\documentclass[11pt]{article}
\usepackage{bbm}
\usepackage{eqnarray,amsmath,amsfonts,amsthm,mathrsfs}
\usepackage{color}
\usepackage{bm}
\usepackage{enumerate}
\usepackage{amssymb}
\usepackage{caption}
\usepackage{graphicx}
\usepackage{epsfig}
\usepackage{epsf}
\usepackage{float}
\usepackage{subfigure}
\usepackage{amsfonts,amsmath,amsthm,amssymb,graphicx,float,fancyhdr,multirow,hyperref}
\usepackage{booktabs,longtable,authblk}
\usepackage{mathrsfs,hhline}

\usepackage{fancyhdr}
\usepackage[top=2.5cm, bottom=2.5cm, left=3cm, right=3cm]{geometry}
\usepackage{graphicx}
\newtheorem{theorem}{Theorem}

\newtheorem{corollary}{Corollary}
\newtheorem{definition}{Definition}

\newtheorem{lemma}{Lemma}
\newtheorem{proposition}{Proposition}
\newtheorem{remark}{Remark}

\usepackage[ruled,linesnumbered]{algorithm2e}

\newcommand{\ppart}[1]{[#1]_+}
\allowdisplaybreaks[4]

\numberwithin{equation}{section}

\title{Weak Physics-Informed Neural Networks for Geometry Compatible Hyperbolic Conservation Laws on Manifolds$^\dag$\footnotetext{\dag~The work of Lei Shi is partially supported by the National Natural Science Foundation of China (Grant No.12171039). The corresponding author is Lei Shi.}}

\author[1]{Hanfei Zhou}
\author[1,2]{Lei Shi}

\affil[1]{School of Mathematical Sciences, \linebreak
Fudan University, Shanghai, 200433, China 
}
\affil[2]{
Shanghai Key Laboratory for Contemporary Applied Mathematics, \linebreak
Fudan University, Shanghai, 200433, China \linebreak
Email:zhouhf23@m.fudan.edu.cn, leishi@fudan.edu.cn
}
\date{}

\begin{document}
	\maketitle
 
\begin{abstract}
Physics-informed neural networks (PINNs) provide a mesh-free approach to solving high-dimensional PDEs on complex geometries, but their theoretical foundations on manifolds remain limited. Moreover, conventional PINN analyses typically rely on solution smoothness, while PINNs may perform poorly for low-regularity solutions arising from nonlinear hyperbolic equations. In this paper, we develop a weak PINN (wPINN) framework for approximating entropy solutions of geometry-compatible hyperbolic conservation laws on Riemannian manifolds $\mathcal{M}^d$. Building on the well-posedness theory, we establish a localized $L_1$-stability estimate that converts localized entropy residuals into terminal error bounds and leads to a rigorous convergence analysis of the proposed method. We then derive approximation guarantees for time-dependent entropy solutions on manifolds, revealing how approximation errors accumulate over long time horizons. For the quadrature error, we develop a problem-adapted localization complexity analysis and show that, for a fixed adversarial test-network architecture, the solution-network contribution achieves the fast rate $\mathrm{VC}_{\mathcal F}/n$, up to logarithmic factors. The resulting algebraic network-complexity exponent depends only on the intrinsic dimension $d$, rather than the ambient dimension. For fixed localization scales, and up to logarithmic factors and the localization bias, the solution-network statistical exponent matches the corresponding minimax exponent in $d$-dimensional Euclidean Sobolev approximation. Numerical experiments illustrate that the proposed wPINN framework accurately approximates entropy solutions on manifold geometries.
\end{abstract}


\section{Introduction}\label{Section: Introduction}
The remarkable success of machine learning techniques—particularly deep learning~\cite{LeCun2015DeepL}—in industrial applications has spurred their growing adoption in scientific computing. Deep neural networks (DNNs), known for their powerful nonlinear approximation capabilities, have been extensively employed to construct surrogate models for approximating solutions of various partial differential equations (PDEs). While classical numerical methods—such as finite difference, finite element, and spectral methods—possess well-established theoretical foundations and deliver strong computational performance for low-dimensional PDEs, they encounter significant challenges in high-dimensional settings due to the curse of dimensionality, which severely limits their efficiency. In recent years, the integration of deep learning with data-driven architectures has led to notable advances in solving high-dimensional PDEs.

In general, algorithms that leverage DNNs to solve PDEs can be broadly classified into two categories: data-driven methods and physics-driven methods. Data-driven methods, which follow the supervised learning paradigm, aim to approximate the solution operator of PDEs, that is, the mapping from input parameters, boundary conditions, or initial conditions to the corresponding solution space. These approaches typically require large volumes of training data, which may be obtained from experimental measurements or precomputed high-fidelity numerical simulations. Notable examples include neural operators~\cite{Li2020NeuralOG, li2021fourier, Kovachki2023NeuralOL,li2024physics} and DeepONets~\cite{Lu2019LearningNO}. In contrast, physics-driven methods adopt an unsupervised learning framework, utilizing the expressive power of DNNs to directly fit the underlying physical laws, which are governed by PDEs. Representative algorithms in this category include physics-informed neural networks (PINNs)~\cite{Raissi2019PhysicsinformedNN, Chen2021PhysicsinformedML} and the deep Ritz method~\cite{Weinan2017TheDR}, which minimize the residuals derived from the classical and variational formulations of PDEs, respectively. For a comprehensive review of deep learning approaches for solving PDEs, we refer the reader to~\cite{de2024numerical}. Despite significant theoretical advances in the analysis of PINNs in recent years~\cite{shin2020convergence,Mishra2020EstimatesOT, biswas2022error, Ryck2022ErrorEF, de2022error}, the majority of existing studies have focused on problems posed in Euclidean domains. As a result, the theoretical underpinnings of PINNs for solving PDEs on manifolds remain largely unexplored. However, PDEs defined on manifolds arise naturally in a wide range of applications, including texture synthesis~\cite{Turk1991GeneratingTO,song2003topics}, image inpainting~\cite{bertalmio2001navier,grossauer2004combined}, and geophysical modeling~\cite{bonito2020divergence,bachini2021intrinsic,jankuhn2021error,haltiner1980numerical,garcia2019shallow}. Although recent efforts have begun extending the PINN framework to PDEs on manifolds~\cite{costabal2024delta, Bastek2023PhysicsInformed}, a rigorous and comprehensive convergence analysis for such methods is still lacking.

Theoretical analyses reveal that the success of PINNs strongly depends on both the well-posedness of the underlying PDEs and the expressive power of neural networks (cf.~\cite{de2024numerical}). More specifically, this dependence can be summarized as follows:  

\begin{itemize}
    \item Existence and uniqueness of solutions: These fundamental properties of PDEs are critical for ensuring the numerical stability of PINNs. Without uniqueness, PINNs may fail to converge or oscillate between multiple potential solutions, resulting in unstable training.
    
    \item Approximation capacity of neural networks: The training objective of PINNs typically involves high-order derivative terms. Ensuring that these terms converge to zero requires both the regularity of the underlying solution and the simultaneous approximation capability of neural networks for high-order derivatives.

    \item Stability of PDEs: Stability allows rigorous control over the total approximation error via the residual term. Analyses in our recent work~\cite{lei2025solving} demonstrate that this can be ensured for elliptic PDEs defined on spheres. In~\cite{lei2025solving}, this property is referred to as the strong convexity condition. Further analysis in the current work shows that this condition improves the error estimates and relaxes constraints on network parameters.
\end{itemize}

Studies in~\cite{Mishra2020EstimatesOT, Ryck2022ErrorEF, Ryck2022GenericBO} have shown that traditional PINNs struggle to efficiently approximate non-smooth solutions of PDEs. This limitation has been confirmed both experimentally and theoretically. To address this issue, several algorithms tailored to low-regularity solutions have been proposed, including weak adversarial networks~\cite{Zang2019WeakAN} and weak PINNs (wPINNs)~\cite{deryck2022wPINN}. The work in~\cite{deryck2022wPINN} introduces wPINNs designed to approximate weak solutions of scalar conservation laws, with a particular focus on their well-posed entropy solutions. The authors develop a theoretical analysis for wPINNs applied to one-dimensional scalar conservation laws and provide empirical evaluations of their performance. Motivated by~\cite{deryck2022wPINN}, we propose a novel wPINN framework for geometry-compatible hyperbolic conservation laws, grounded in the well-posedness theory for conservation laws on Riemannian manifolds~\cite{BenArtzi2006WellposednessTF}. To formulate wPINNs on manifolds, we utilize the measure-valued solutions to hyperbolic conservation laws and establish a more general contraction estimate (the $L_1$ contraction property) to analyze the residuals for wPINNs. Our theoretical analysis proceeds via error decomposition, providing bounds on both the approximation error and quadrature error, the latter of which arises from discretizing the continuous  PDE residuals using numerical quadrature rules~\cite{de2021approximation,longo2021higher, Mishra2020EstimatesOT}. The main contributions of this work are summarized below:

\begin{itemize}
    \item We are the first to establish a wPINN framework for geometry-compatible hyperbolic conservation laws on Riemannian manifolds, extending the structure of wPINNs previously developed for scalar conservation laws in one-dimensional Euclidean space~\cite{deryck2022wPINN}.  Leveraging the well-posedness theory on $d$-dimensional Riemannian manifolds, we derive a more general contraction estimate and further develop an error analysis framework for wPINNs.

    \item Our method imposes less restrictive requirements on the underlying PDEs. By targeting the approximation of weak solutions, our approach alleviates the stringent regularity assumptions (e.g., solutions residing in high-order Sobolev or H{\"o}lder spaces) commonly required in previous studies~\cite{shin2020convergence, Ryck2022ErrorEF, lei2025solving}. We utilize neural networks with ReLU activation functions to approximate weak solutions and validate the effectiveness of this approach numerically. This contrasts with conventional PINNs, which often rely on smooth activation functions such as $\sin$ or $\tanh$, both of which are also employed in the wPINNs framework~\cite{deryck2022wPINN}.  Meanwhile, we introduce a localized ReQU test network class, constructed from extrinsic radial kernels in the ambient coordinates of \(\mathcal M^d\). This structured neural test class is adapted to the embedded manifold geometry and provides the approximate identity, localization, and gradient-cancellation properties required for the localized \(L_1\)-stability estimate for wPINNs, which is essential for the convergence analysis.
    
   \item Regarding approximation error analysis, we are the first to establish approximation rates of neural networks applied to approximate Sobolev functions on manifolds, based on their classical definitions induced by the manifold structure. Previous works have primarily focused on H{\"o}lder classes formulated through ambient coordinates or the ambient Euclidean metric~\cite{Chen2019efficient, Chen2019NonparametricRO, Labate2024LowDA, SchmidtHieber2019DeepRN}, which can be less straightforward to verify intrinsically. The inherently low regularity of weak solutions leads to temporal accumulation of error. Our analytical framework employs a novel neural network architecture that mitigates this accumulation by leveraging the properties of entropy solutions and classical spline interpolation theory. Moreover, for fixed localization scales, and up to logarithmic factors and the localization bias, the statistical exponents of the solution-network contribution match the corresponding minimax exponents for Sobolev functions in $d$-dimensional Euclidean space. This demonstrates that our algorithm can alleviate the curse of dimensionality by exploiting the intrinsic low-dimensional structure of Riemannian manifolds.

    \item Concerning quadrature error analysis, we develop a localized complexity analysis based on the localized \(L_1\)-stability estimate established above. In contrast to the global empirical-process estimates used in~\cite{deryck2022wPINN}, our analysis exploits residual sublevel sets and the stability structure of the wPINN loss. For a fixed adversarial test-network architecture, the part of the quadrature error generated by the growing solution network is bounded by \(C_\lambda^2\frac{\mathrm{VC}_{\mathcal F}\log n}{n}\), up to confidence terms. The remaining fluctuation is attached only to the fixed test family and is separated explicitly in the proof. This makes the dependence on the solution-network VC complexity transparent and relaxes the bounded-parameter requirements used in global covering arguments.
\end{itemize}

The structure of this paper is organized as follows. In~\autoref{Section: Preliminaries and Main Results}, we introduce the well-posedness theory of geometry-compatible conservation laws on manifolds, formulate the wPINN framework, and state the main convergence result. In~\autoref{Section: Approximation Error}, we derive the approximation rates for entropy solutions on manifolds. In~\autoref{Section: Localized L1 contraction}, we establish a localized \(L_1\)-stability estimate that links entropy residuals to the terminal error. In~\autoref{Section: Quadrature Error}, we analyze the quadrature error arising from the empirical residual. In~\autoref{Section: Proof of Main Result}, we combine the preceding estimates to prove the overall convergence rate. In~\autoref{Section: Numerical Experiment}, we present numerical experiments demonstrating the effectiveness of wPINNs on manifolds. All proofs are provided in the \hyperref[allapp]{Appendix}.

\section{Preliminaries and Main Results}\label{Section: Preliminaries and Main Results}
\subsection{Geometry-Compatible Hyperbolic Conservation Laws}\label{Subsection: Geometry compatible hyperbolic conservation laws}

In this subsection, we introduce geometry-compatible hyperbolic conservation laws on manifolds. We first recall some basic notations from~\cite{BenArtzi2006WellposednessTF} for further statement. Let $(\mathcal{M}^d,g)$ be a smooth, closed (that is, compact and without boundary), oriented $d$-dimensional Riemannian manifold. Throughout this paper, we specify that $d \geq 2$. We have tried to make the paper as self-contained as possible; some background and notations are provided in~\autoref{Appendix: Notations} for readability. Recall that Kruzkov’s theory~\cite{kruvzkov1970first} concerns scalar conservation laws in the Euclidean space $\mathbb{R}^n$, establishing the existence, uniqueness, and stability of $L_1 \cap L_\infty$ entropy solutions. In the case of conservation laws with constant flux—that is, when the flux is independent of time or space—the solutions satisfy the maximum principle, $L_1$ contraction, and have bounded total variation. However, it is important to note that a constant flux generally does not exist on a Riemannian manifold. To address this issue, the notion of a geometry-compatible flux was introduced in~\cite{BenArtzi2006WellposednessTF}, leading to hyperbolic conservation laws on $(\mathcal{M}^d, g)$ and establishing corresponding well-posedness results.

For convenience, we follow the notations used in~\cite{BenArtzi2006WellposednessTF}. A flux on the manifold $\mathcal{M}^d$ is a vector field $f = f_x(\bar{u})$ that depends on the parameter $\bar{u}$, with smooth dependence on both variables. Consider the following hyperbolic conservation laws with flux $f$ on $\mathcal{M}^d$:
\begin{equation}\label{equation: hyperbolic conservation laws}
    \partial_t u + \mathbf{div}_g(f_x(u)) = 0.
\end{equation} Here, for each fixed time $t$, the divergence operator $\mathbf{div}_g$ is applied to the vector field $x \to f_x(u(x, t)) \in T_x\mathcal{M}^d$. A flux is said to be geometry-compatible if it satisfies the condition
\begin{equation}\label{equation: constant flux}
    \mathbf{div}_g(f_x(\bar{u}))  = 0,  \quad  \bar{u} \in \mathbb{R}, x \in \mathcal{M}^d.
\end{equation}
Lemma 3.2 in~\cite{BenArtzi2006WellposednessTF} shows that, under the geometry-compatible condition, the conservation law on the manifold takes a non-conservative form when expressed in local coordinates. Given an initial value $u_0 \in L_\infty(\mathcal{M}^d)$, we consider the initial condition:
\begin{equation}\label{equation: initial condition}
    u(x,  0) = u_0(x),  \quad x \in\mathcal{M}^d.
\end{equation} Following~\cite{BenArtzi2006WellposednessTF}, we introduce the definition of weak solutions.
\begin{definition}[cf. Section 3 in~\cite{BenArtzi2006WellposednessTF}]
    A function $u \in L_\infty(\mathcal{M}^d \times \mathbb{R}_+ )$ is called a weak solution if it satisfies the geometry-compatible hyperbolic conservation law \eqref{equation: hyperbolic conservation laws} with flux $f$ satisfying \eqref{equation: constant flux}, in the sense of distributions, and satisfies the initial condition \eqref{equation: initial condition}.
\end{definition}

However, the uniqueness of weak solutions is not guaranteed, necessitating an entropy condition to ensure well-posedness. We introduce a convex entropy-flux pair to establish well-posedness to obtain additional conservation laws. A convex entropy-flux pair is a pair $(U,F)$, where $U : \mathbb{R} \to \mathbb{R}$ is a convex function, and $F = F_x( \bar {u}), \bar{u} \in \mathbb{R}$, is a vector field defined by 
\begin{equation*}
    F_x(\bar{u})=\int^{\bar{u}}\partial_{u'} U (u') \partial_{u'}f_x(u') du', \quad \bar{u} \in \mathbb{R},  x \in \mathcal{M}^d.
\end{equation*}
The definition of an entropy solution is given below.
\begin{definition}[Definition 3.3 in~\cite{BenArtzi2006WellposednessTF}]\label{definition: entropy solution}
    Let $f = f_x(\bar{u})$ be a geometry-compatible flux on $(\mathcal{M}^d, g)$.
    Given an initial value $u_0 \in L_\infty(\mathcal{M}^d)$, $u\in L_\infty(\mathcal{M}^d \times \mathbb{R}_+)$ is called an entropy solution of the initial value problem satisfying \eqref{equation: hyperbolic conservation laws}, \eqref{equation: constant flux}, and \eqref{equation: initial condition}, if for every convex entropy-flux pair $(U,F)$ and every nonnegative smooth test function $\xi\in C_c^\infty(\mathcal M^d\times [0,\infty))$,
    \begin{equation*}
    \begin{aligned}
        &\int_0^\infty\!\int_{\mathcal{M}^d}
        \Bigl(U(u(x,t))\,\partial_t\xi(x,t)
        +g_x\bigl(F_x(u(x,t)),\nabla_g\xi(x,t)\bigr)\Bigr)
        \,dV_g(x)\,dt
        \\
        &\qquad+
        \int_{\mathcal M^d}U(u_0(x))\xi(x,0)\,dV_g(x)
        \geq0.
    \end{aligned}
    \end{equation*}
    Equivalently, if $\xi$ is compactly supported in $\mathcal M^d\times (0,\infty)$, the initial trace term vanishes and the entropy inequality reduces to the space-time integral only.
\end{definition}

\begin{lemma}
\label{lemma: admissible Sobolev entropy tests}
Let $u$ be an entropy solution and let $0<T<\infty$. The interior entropy inequality remains valid for every nonnegative \(\xi\in W_0^{1,\infty}(\mathcal M^d\times(0,T)).\) 
\end{lemma}

\begin{proof}
Since $\xi$ has zero trace on the temporal boundary and compact support in the interior in the Sobolev sense, there exist nonnegative functions $\xi_k\in C_c^\infty(\mathcal M^d\times(0,T))$ such that
\[
\xi_k\longrightarrow\xi
\qquad\text{in }W^{1,1}(\mathcal M^d\times(0,T)).
\]
This follows by using a finite atlas, a partition of unity, and nonnegative Euclidean mollifiers. The state $u$ is bounded, and the flux is smooth on the compact interval containing its essential range. Hence $U(u)$ and $F_x(u)$ are uniformly bounded. Applying the entropy inequality to $\xi_k$ and passing to the limit in the two first-derivative terms proves the claim.
\end{proof}

Now, we are in a position to state the well-posedness theory. We define the following bounded variation function space, which plays a crucial role in characterizing the regularity of the initial condition $u_0$:
\begin{equation*}
\begin{aligned}
\mathrm{BV}(\mathcal M^d;dV_g)
:=\Bigl\{h\in L_1(\mathcal M^d;dV_g):\;&\mathrm{TV}(h)<\infty\Bigr\},
\\[-1mm]
\mathrm{TV}(h)
:=\sup_{\substack{\psi\in C^1(\mathcal M^d;T\mathcal M^d)\\
\|\psi\|_{L_\infty(\mathcal M^d)}\le1}}
&\int_{\mathcal M^d}h\,\mathbf{div}_g\psi\,dV_g.
\end{aligned}
\end{equation*}
where $\psi$ describes all $C^1$ vector fields on the manifold. We refer to $\mathrm{TV}(h)$ as the total variation of the function $h$.
\begin{proposition}[Theorem 4.4 in~\cite{BenArtzi2006WellposednessTF}]\label{proposition: well-posedness theory of Geometry compatible hyperbolic conservation laws}
    Let $u_0 \in L_\infty(\mathcal{M}^d) \cap \mathrm{BV} (\mathcal{M}^d;dV_g)$, then there exists an entropy solution $u$ to the initial value problem satisfying \eqref{equation: hyperbolic conservation laws}, \eqref{equation: constant flux}, and \eqref{equation: initial condition}. Moreover, this solution satisfies the following properties: \begin{equation}\label{equation: L1 contraction1 with bounded variation function}
        \|u(t)\|_{L_p(\mathcal{M}^d;dV_g)} \leq \|u_0\|_{L_p(\mathcal{M}^d;dV_g)}, \quad t \in \mathbb{R}_+,  1 \leq p \leq \infty,
    \end{equation} and
    \begin{equation}\label{equation: bounded variation inequality}
        \|u(t) - u(t')\|_{L_1(\mathcal{M}^d;dV_g)} \leq \mathrm{TV} (u_0) |t-t'|, \quad 0 \leq t' \leq  t.
    \end{equation}
    And there exists a constant $C_1 > 0$ depending only on $\|u_0\|_{L_\infty(\mathcal{M}^d)}$ and $\mathcal{M}^d$, such that
    \begin{equation}\label{equation: bounds on TV}
        \mathrm{TV}(u(t)) \leq e^{C_1t}  (1 + \mathrm{TV} (u_0)), \quad t \in \mathbb{R}_+.
    \end{equation}
    Furthermore, if $u$ and $v$ are entropy solutions corresponding to the initial values $u_0$ and $v_0$, respectively, then the following contraction inequality holds:
    \begin{equation*}
        \|v(t) - u(t)\|_{L_1(\mathcal{M}^d;dV_g)} \leq \|v_0 - u_0\|_{L_1(\mathcal{M}^d;dV_g)}, \quad t \in \mathbb{R}_+.
    \end{equation*}
\end{proposition}
It is worth noting that the assumption of bounded total variation for the initial condition here can be relaxed. In~\cite{BenArtzi2006WellposednessTF}, the authors proved a more general well-posedness theory based on the entropy measure-valued solution. However, in that setting, the total variation estimate is no longer valid.
\subsection{Mathematical Definition of Neural Networks}\label{Subsubsection: Mathematical Definition of Neural Networks}

Given $\sigma:\mathbb{R}\to\mathbb{R}$, $\mathbf b=(b_1,\ldots,b_r)^T$, and $\mathbf y=(y_1,\ldots,y_r)^T$, define
\[
\sigma_{\mathbf b}(\mathbf y)
=
(\sigma(y_1-b_1),\ldots,\sigma(y_r-b_r))^T.
\]
For an input dimension $m\in\mathbb N_+$, a real-valued feedforward neural network is a function
\begin{equation}\label{equation: Neural Network}
 f:\mathbb R^m\to\mathbb R,
 \qquad
 x\longmapsto
 W_L\sigma_{\mathbf v_L}W_{L-1}\sigma_{\mathbf v_{L-1}}\cdots
 W_1\sigma_{\mathbf v_1}W_0x.
\end{equation}
We denote the number of hidden layers by $L$ and the width of the $i$-th hidden layer by $p_i$. Thus $W_i$ has size $p_{i+1}\times p_i$, where $p_0=m$ and $p_{L+1}=1$, while the bias vector satisfies $\mathbf v_i\in\mathbb R^{p_i}$. Throughout the theoretical approximation analysis, the activation function is the ReLU function $\sigma(x)=\max\{x,0\}=(x)_+$.

Let $(Q,G)\in[0,\infty)^2$ and $(S,B,F)\in[0,\infty]^3$. We define
\begin{equation}\label{equation: Assumption for hypothesis space}
\mathcal F_m(Q,G,S,B,F)
:=
\left\{
 f:\mathbb R^m\to\mathbb R
 \left\vert
 \begin{array}{l}
 f\text{ has the form \eqref{equation: Neural Network}},\\
 L\le Q,\quad \max_{1\le i\le L}p_i\le G,\\
 \displaystyle\sum_{i=0}^L\|W_i\|_0+
 \sum_{i=1}^L\|\mathbf v_i\|_0\le S,\\
 \displaystyle\max_{0\le i\le L}\|W_i\|_\infty
 \vee
 \max_{1\le i\le L}\|\mathbf v_i\|_\infty\le B,\\
 \|f\|_{L_\infty}\le F
 \end{array}
\right.
\right\}.
\end{equation}
When the input dimension is clear, we write $\mathcal F(Q,G,S,B,F)$. In particular, the spatial approximation networks have input dimension $D$, whereas the time-dependent solution networks have input dimension $D+1$. The parameters $B$ and $F$ in~\eqref{equation: Assumption for hypothesis space} may equal $\infty$, indicating that no parameter or output bound is imposed; similarly, $S=\infty$ removes the sparsity constraint. Some standard fixed-depth approximation constructions, including those in~\cite{yarotsky2017error,guhring2020error}, use parameter bounds that grow as the target accuracy increases. Our bounded-weight duplication construction instead encodes large coefficients through logarithmic depth and therefore also covers the case $B=1$. In the quadrature analysis, uniform parameter-space covering numbers are replaced by empirical $L_2(P_n)$-covering numbers controlled by pseudo-dimension. This permits both $B=1$ and $B=\infty$ and separates the growing solution-network complexity from the fixed localized test architecture; see~\autoref{Appendix: Analysis of Quadrature Error}. It also differs from global bounded-parameter analyses such as~\cite{lu2022machine}.
\subsection{Construction of wPINNs on Manifolds}\label{Subsubsection: Manifolds Weak PINNs}

According to the Nash embedding theorem~\cite{nash1956imbedding}, any compact manifold can be isometrically embedded into a higher-dimensional Euclidean space. Accordingly, we consistently assume that the manifold $(\mathcal{M}^d,g)$ is embedded into the ambient space $\mathbb{R}^D$, where $d \leq D$.

Due to the limitations of classical PINNs in approximating weak solutions of hyperbolic conservation laws, we focus on approximating entropy solutions over the time interval $[0,T]$ using wPINNs. Specifically, we formulate Kruzkov-type entropy residuals for geometry-compatible hyperbolic conservation laws. We use the convex entropy-flux pairs
\begin{equation*}
    \widetilde{U}(\bar{u}, \bar{v}) = |\bar{v} - \bar{u}|, \quad \widetilde{F}(\bar{u}, \bar{v}) = (f (\bar{v})-f(\bar{u})) \mathrm{sgn}(\bar{v}-\bar{u}).
\end{equation*}
We introduce the signed entropy residual first.
\begin{definition}\label{definition: entropy residual}
For \(v \in (L_\infty \cap L_1)(\mathcal{M}^d \times [0, T])\),  
\(\xi \in W_0^{1, \infty}(\mathcal{M}^d \times [0, T])\), and \(c \in \mathbb{R}\), 
the signed Kruzkov entropy residual is defined as
\begin{equation*}
    \mathcal{R}^{int}(v, \xi, c)
    =
    \int\int_{\mathcal{M}^d \times [0, T]}
    r^{int}(v, \xi, c)\,dV_g(x)dt.
\end{equation*}
Here the internal residual density is defined as 
\begin{equation*}
    r^{int}(v, \xi, c) 
    =
    -\widetilde{U}(v, c)\partial_t\xi(x, t)
    -
    g_x\bigl(\widetilde{F}(v, c), \nabla_g\xi(x, t)\bigr).
\end{equation*}
\end{definition}

The entropy inequality is one-sided and is imposed only on nonnegative test functions. Therefore the theoretical wPINN objective is defined by the positive entropy violation rather than by an unconstrained signed supremum. Let \(\Xi\subset W_0^{1,\infty}(\mathcal M^d\times(0,T))\) be a nonnegative and uniformly bounded test class, and let \(\mathcal C\subset\mathbb R\) be a compact interval containing the essential range of the entropy solution. Define
\begin{equation}\label{equation: positive internal loss}
    \mathcal L_\Xi^+(v)
    :=
    \sup_{\xi\in\Xi,\ c\in\mathcal C}
    \ppart{\mathcal R^{int}(v,\xi,c)} ,
\end{equation}
and
\begin{equation*}
    \mathcal J_\Xi(v)
    :=
    \mathcal R^{tb}(v)
    +
    T\operatorname{Vol}_g(\mathcal M^d)\mathcal L_\Xi^+(v),
    \qquad
    \mathcal R^{tb}(v)
    :=
    \int_{\mathcal M^d}|v(x,0)-u_0(x)|\,dV_g(x).
\end{equation*}
The positive part in \eqref{equation: positive internal loss} is essential. Without a uniform normalization of the test functions, a signed supremum over \(W_0^{1,\infty}\) is not finite unless the residual distribution is identically zero, because the test function can be multiplied by an arbitrary positive scalar. Moreover, a negative signed residual does not violate the entropy inequality.

For the quadrature analysis below, the test class is the localized ReQU family
\(\overline\Phi_\lambda\) introduced in \autoref{theorem: L1 contraction for total error}. We regard this family as a fixed-architecture adversarial test network. Its center, time-location, and horizon parameters are optimized in the inner maximization, but the architecture and the number of test parameters are fixed independently of the quadrature sample size and independently of the size of the solution network. This is the theoretical counterpart of the adversarial test update in Algorithm~\ref{algorithm: training algorithm}. The fixed test-side parameter entropy is absorbed into constants depending on \(\lambda\); the only network complexity that grows with \(n\) is therefore \(\mathrm{VC}_{\mathcal F}:=SL\log S.\) A precise covering formulation of this fixed architecture is given in \autoref{Appendix: Analysis of Quadrature Error}.

In practice, the integrals are approximated by numerical quadrature rules. In our implementation, we adopt Monte Carlo quadrature. Let \(P\) denote the uniform probability distribution on \(\mathcal{M}^d \times[0, T]\). By sampling \(\{(X_i,t_i)\}_{i=1}^n\) independently and identically from \(P\), define
\begin{equation*}
     \mathcal{R}^{int}_n(u, \xi, c)
     =
     \frac{T\operatorname{Vol}_g(\mathcal M^d)}{n}\sum_{i=1}^n
     r^{int}(u(X_i,t_i), \xi(X_i,t_i), c).
\end{equation*}
For the initial-time term, let \(\{X_i^{0}\}_{i=1}^{n_0}\) be i.i.d. samples from the uniform distribution on \(\mathcal M^d\). We write
\begin{equation*}
     \mathcal{R}_n^{tb}(u)
     =
     \frac{\operatorname{Vol}_g(\mathcal M^d)}{n_0}\sum_{i=1}^{n_0}
     |u(X_i^0,0)-u_0(X_i^0)|.
\end{equation*}
For notational simplicity, and as in the sequel, we take \(n_0=n\). The empirical positive wPINN objective associated with the test class \(\Xi\) is
\begin{equation}\label{equation: positive empirical objective}
    \mathcal J_{\Xi,n}(u)
    :=
    \mathcal R_n^{tb}(u)
    +
    T\operatorname{Vol}_g(\mathcal M^d)
    \sup_{\xi\in\Xi,\ c\in\mathcal C}
    \ppart{\mathcal R_n^{int}(u,\xi,c)} .
\end{equation}
The theoretical empirical estimator is any exact minimizer
\begin{equation}\label{equation: empirical minimizer positive}
    u_n\in \arg\min_{u\in\mathcal F}\mathcal J_{\Xi,n}(u).
\end{equation}
The numerical implementation in \autoref{Section: Numerical Experiment} uses a normalized adversarial loss and squares the positive violation for stability, but the convergence proof is stated for \eqref{equation: positive empirical objective} because the map \(z\mapsto \ppart z\) is \(1\)-Lipschitz and preserves the localized empirical-process estimates used in \autoref{Section: Quadrature Error}.

\subsection{Main Results}\label{subsection: Main Results}
In this subsection, we state the convergence results in three levels. We first give a bound with the class-wise localization bias left explicit. We then restrict the empirical minimization to the temporally controlled approximation class supplied by \autoref{corollary: TV controlled temporal modulus}, which converts the abstract bias into an explicit deterministic term. Finally, we state the improved rate under uniform-in-time Sobolev regularity.

For two positive sequences \(\{A_n\}_{n\ge1}\) and \(\{B_n\}_{n\ge1}\), we write \(A_n\lesssim B_n\) if there exists a constant \(C>0\), independent of \(n\), such that \(A_n\le CB_n\). We write \(A_n\asymp B_n\) if both \(A_n\lesssim B_n\) and \(B_n\lesssim A_n\) hold. Let \(u^*\) be the entropy solution in \autoref{definition: entropy solution}.

The localized ReQU test functions used below are defined in \autoref{theorem: L1 contraction for total error}. We write
\[
  \overline\Phi_\lambda=\bigcup_{\tau\in[4h,T]}\Phi_{\lambda,\tau},
  \qquad \lambda=(\alpha,\eta,h),
\]
and
\[
  \mathfrak L_\lambda^+(v)
  :=
  \sup_{\xi\in\overline\Phi_\lambda,\ c\in\mathcal C}
  \ppart{\mathcal R^{int}(v,\xi,c)}.
\]
The family \(\overline\Phi_\lambda\) is optimized in the adversarial step, but it has a fixed ReQU architecture and a compact finite-dimensional parameter set. Its metric-entropy constants may depend on \(\lambda\), but do not grow with the quadrature sample size or with the size of the solution network. For the solution class \(\mathcal F(L,W,S,B,M)\), set $\mathrm{VC}_{\mathcal F}:=SL\log S.$ The reduction of the joint residual entropy to this quantity is proved in \autoref{lemma: empirical residual covering}.

For any trial class \(\mathcal F\subset L_\infty(\mathcal M^d\times[0,T])\cap C([0,T];L_1(\mathcal M^d))\), define the uniform localization bias
\begin{equation*}
  b_{\lambda,T}(\mathcal F)
  :=
  \sup_{v\in\mathcal F}C_T\left[
  \omega_t(v;C\eta)+\omega_t(v;Ch)+ (\eta+h)TV(u_0)
  +\alpha(1+TV(u_0))\right].
\end{equation*}

\begin{theorem}\label{theorem: main result}
Let
\[
M\ge \max\{\|u^*\|_{L_\infty(\mathcal M^d\times[0,T])},\|u_0\|_{L_\infty(\mathcal M^d)}\},
\qquad
a=1/(d+2).
\]
Fix a compact interval \(\mathcal C\subset[-M,M]\) containing the essential range of \(u^*\).
Assume that \(u_0\in L_\infty(\mathcal M^d)\cap BV(\mathcal M^d)\), that \(\mathcal M^d\) admits smooth local coordinates in the sense of \autoref{definition: smooth coordinates}, and that the flux is geometry-compatible and smooth enough so that, for \(|u|\le M\),
\[
   \|\partial_u f(\cdot,u)\|_{L_\infty}
   +\|\nabla_x f(\cdot,u)\|_{L_\infty}
   +\|\nabla_x\partial_u f(\cdot,u)\|_{L_\infty}<\infty.
\]
Assume also that the adversarial test architecture is fixed in the sense described above and in \autoref{Appendix: Analysis of Quadrature Error}. For either \(B=1\) or \(B=\infty\), choose \(\mathcal F=\mathcal F(L,W,S,B,M)\) with the network sizes specified below, and assume \(b_{\lambda,T}(\mathcal F)<\infty\). Let \(u_n\) be any exact minimizer of the localized empirical positive wPINN objective \(\mathcal J_{\overline\Phi_\lambda,n}\) defined in \eqref{equation: positive empirical objective} over \(\mathcal F\). Then, for all sufficiently large \(n\), with probability at least \(1-\exp\left(-c\,n^{1-2a}(\log n)^{8a}\right),\) there holds
\begin{equation}\label{equation: main theorem bias rate}
\|u_n(\cdot,T)-u^*(\cdot,T)\|_{L_1(\mathcal M^d)}
\lesssim
C_\lambda^2 n^{-a}(\log n)^{4a}
+b_{\lambda,T}(\mathcal F).
\end{equation}
The network sizes can be chosen so that
\begin{align*}
L&\lesssim \log\bigl(n(\log n)^{-4}\bigr),\\
W&\lesssim n^{a(d+1)}(\log n)^{-4a(d+1)},\\
S&\lesssim n^{a(d+1)}(\log n)^{-4a(d+1)}
       \log\bigl(n(\log n)^{-4}\bigr).
\end{align*}
The hidden constants are independent of \(n\). They may depend on the fixed embedding of \(\mathcal M^d\), and hence on \(D\), as well as on the geometric and flux constants, \(M,T\), and the fixed test architecture. The algebraic exponents in \(n\) depend only on the intrinsic dimension \(d\).
\end{theorem}

\begin{remark}
In our proof, we built upon the well-posedness theory of hyperbolic conservation laws, specifically leveraging the temporal Lipschitz continuity of entropy solutions, as established in~\eqref{equation: bounded variation inequality}. Utilizing the bounded total variation property of entropy solutions,  we construct a $\mathcal{W}^1_1$ function approximation, and subsequently apply spline interpolation to build neural networks that approximate the solution at discrete time points. Notably, the overall convergence rate depends only on the intrinsic dimension $d$ of the underlying manifold, which may be significantly smaller than the ambient dimension $D$. For fixed localization scales, the solution-network contribution achieves the classical minimax exponent for approximating functions in $W^1_1$ over $d$-dimensional Euclidean domains, up to logarithmic factors and the localization bias. This result underscores the power of exploiting low-dimensional manifold structures to mitigate the curse of dimensionality.
\end{remark}

\begin{remark}
    A key finding from our analysis is that the constant term in~\eqref{equation: main theorem bias rate} depends on the total variation of the entropy solution. Furthermore, the estimate provided in~\eqref{equation: bounds on TV} suggests that this term scales exponentially with respect to the final time $T$, approximately as $e^{C_1T}$. This highlights the significant impact of long time horizons on the overall algorithmic error. Such time-accumulated error has been observed experimentally in~\cite{deryck2022wPINN}, where approaches like total variation regularization have been proposed to mitigate its effects. Motivated by both the practical benefits of regularization and the theoretical challenge posed by exponential time scaling, our future work will explore a synergistic strategy. Specifically, we aim to address this limitation by employing accurate approximations over shorter, discretized time intervals, combined with appropriate regularization techniques, to enable robust and reliable long-term predictions.
\end{remark}

\begin{remark}
The localization bias \(b_{\lambda,T}(\mathcal F)\) in~\autoref{theorem: main result} appears because the localized \(L_1\)-stability estimate involves the temporal modulus $\omega_t(v;r):=\sup_{\substack{s,t\in[0,T], |t-s|\le r}}\|v(\cdot,t)-v(\cdot,s)\|_{L_1(\mathcal M^d)}$ of trial functions. For a general trial class \(\mathcal F\), this modulus is not explicitly controlled. However, by~\autoref{corollary: TV controlled temporal modulus}, the approximation network constructed in~\autoref{theorem: approximation error for entropy solutions} can be chosen to satisfy the TV-controlled temporal modulus $\omega_t(v;r)\le C(1+TV(u_0))r .$ This observation suggests restricting the empirical minimization to a temporally controlled subclass. Such a restriction preserves the approximation rate, because the constructed approximation network remains admissible, and it does not enlarge the quadrature error since the subclass is contained in the original trial class. It also turns the abstract localization bias into the explicit term $C_T(1+TV(u_0))(\alpha+\eta+h).$ The following corollary records this explicit-bias version of~\autoref{theorem: main result}.
\end{remark}

The abstract bias in \eqref{equation: main theorem bias rate} becomes explicit on the temporally controlled subclass suggested by the approximation construction. Define
\begin{equation*}
\mathcal F_{\rm tm}(L,W,S,B,M)
:=
\left\{
 v\in\mathcal F(L,W,S,B,M)
 \left\vert
 \begin{array}{l}
 \omega_t(v;r)\le C_{\rm tm}(1+TV(u_0))r,\\
  0\le r\le T
 \end{array}
\right.
\right\},
\end{equation*}
where \(C_{\rm tm}\) is chosen large enough that the approximation network in \autoref{corollary: TV controlled temporal modulus} belongs to this class.

\begin{corollary}\label{corollary: explicit bias convergence}
Under the assumptions of \autoref{theorem: main result}, let \(u_n\) be an exact minimizer of \(\mathcal J_{\overline\Phi_\lambda,n}\) over \(\mathcal F_{\rm tm}(L,W,S,B,M)\). Then the same network-size bounds and confidence level hold, and
\begin{equation}\label{equation: explicit bias convergence rate}
\|u_n(\cdot,T)-u^*(\cdot,T)\|_{L_1(\mathcal M^d)}
\lesssim
C_\lambda^2 n^{-a}(\log n)^{4a}
+C_T(1+TV(u_0))(\alpha+\eta+h).
\end{equation}
\end{corollary}

If the entropy solution has additional Sobolev regularity, then the nodewise spatial approximation is more efficient.

\begin{corollary}\label{corollary: sobolev convergence rate}
Assume, in addition to the hypotheses of \autoref{theorem: main result}, that for some \(s\in\mathbb N_+\), the entropy solution admits a representative satisfying
\[
u^*(\cdot,t)\in\mathcal W_1^s(\mathcal M^d)
\quad\text{for every }t\in[0,T],
\qquad
\sup_{0\le t\le T}\|u^*(\cdot,t)\|_{\mathcal W_1^s(\mathcal M^d)}<\infty.
\]
Set \(a_s:=s/(2s+d).\) Then, for either \(B=1\) or \(B=\infty\), there exists a choice of \(L,W,S\) and an exact minimizer \(u_n\in\mathcal F(L,W,S,B,M)\) of \(\mathcal J_{\overline\Phi_\lambda,n}\) such that, for all sufficiently large \(n\), with probability at least \(1-\exp\left(-c\,n^{1-2a_s}(\log n)^{8a_s}\right),\) there holds
\begin{equation}\label{equation: sobolev rate abstract bias}
\|u_n(\cdot,T)-u^*(\cdot,T)\|_{L_1(\mathcal M^d)}
\lesssim
C_\lambda^2 n^{-a_s}(\log n)^{4a_s}
+b_{\lambda,T}(\mathcal F).
\end{equation}
The network sizes can be chosen as
\begin{align*}
L&\lesssim \log\bigl(n(\log n)^{-4}\bigr),\\
W&\lesssim n^{a_s(d/s+1)}(\log n)^{-4a_s(d/s+1)},\\
S&\lesssim n^{a_s(d/s+1)}(\log n)^{-4a_s(d/s+1)}
       \log\bigl(n(\log n)^{-4}\bigr).
\end{align*}
If the empirical minimization is instead performed over \(\mathcal F_{\rm tm}(L,W,S,B,M)\), then \eqref{equation: sobolev rate abstract bias} improves to
\begin{equation}\label{equation: sobolev rate explicit bias}
\|u_n(\cdot,T)-u^*(\cdot,T)\|_{L_1(\mathcal M^d)}
\lesssim
C_\lambda^2 n^{-a_s}(\log n)^{4a_s}
+C_T(1+TV(u_0))(\alpha+\eta+h).
\end{equation}
\end{corollary}
\begin{remark}
   Under the pointwise-in-time uniform Sobolev assumption, the spatial-localization terms are controlled directly by $TV(u^*(t))\le C\sup_{0\le r\le T}\|u^*(\cdot,r)\|_{\mathcal W_1^s}$, so the corresponding constant does not use the exponential $e^{C_1T}$ BV-growth factor. For fixed localization scales, and up to logarithmic factors and the localization bias, the solution-network contribution has the minimax exponent $n^{-s/(2s+d)}$ associated with Sobolev approximation over $d$-dimensional Euclidean domains. The algebraic exponent remains independent of the ambient dimension $D$, further highlighting that our approach effectively mitigates the curse of dimensionality. However, to the best of our knowledge, the specific minimax optimal rates for solving PDE directly on the $d$-dimensional manifold $\mathcal{M}^d$, while leveraging its intrinsic structure, have yet to be established. A rigorous investigation into these manifold-specific minimax rates is left for future work. 
\end{remark}

\paragraph{Error decomposition.}
We briefly display the decomposition used in the proof. The detailed version, including all scale factors, is given in \autoref{Appendix: Error Decomposition}. Let
\[
   \mathfrak Q_\lambda(v)
   :=\mathcal R^{tb}(v)
   +T\operatorname{Vol}_g(\mathcal M^d)\mathfrak L_\lambda^+(v)
   +b_{\lambda,T}(\mathcal F),
\]
and let \(\mathfrak Q_{\lambda,n}\) be its empirical analogue. Since the bias is constant on \(\mathcal F\), the empirical minimizer of \(\mathfrak Q_{\lambda,n}\) is the same as that of the localized objective \(\mathcal J_{\overline\Phi_\lambda,n}\) in \eqref{equation: positive empirical objective}. If $u_{\mathcal F}\in\arg\min_{v\in\mathcal F}\mathfrak Q_\lambda(v),$ then empirical minimality gives the basic inequality
\begin{equation}\label{equation: main basic error decomposition}
\mathfrak Q_\lambda(u_n)-\mathfrak Q_\lambda(u_{\mathcal F})
\le
\bigl(\mathfrak Q_\lambda-\mathfrak Q_{\lambda,n}\bigr)(u_n)
+
\bigl(\mathfrak Q_{\lambda,n}-\mathfrak Q_\lambda\bigr)(u_{\mathcal F}).
\end{equation}
For \(\theta=(\xi,c)\), define
\[
  h_{v,\theta}:=r^{int}(v,\xi,c)-r^{int}(u^*,\xi,c).
\]
The positive-part contraction yields
\[
\left|\mathfrak L_\lambda^+(v)-\mathfrak L_{\lambda,n}^+(v)\right|
\lesssim
\sup_\theta |(P-P_n)h_{v,\theta}|
+
\sup_\theta |(P-P_n)r^{int}(u^*,\xi,c)|.
\]
Thus \eqref{equation: main basic error decomposition} consists of the population approximation term \(\mathfrak Q_\lambda(u_{\mathcal F})\), localized solution-dependent quadrature processes at \(u_n\) and \(u_{\mathcal F}\), and a fixed calibration process at \(u^*\). The key Bernstein relation is \(P h_{v,\theta}^2 \le 2M C_\lambda^2P|v-u^*| \lesssim C_\lambda^2\mathfrak Q_\lambda(v).\) Consequently, on the sublevel set \(\{v:\mathfrak Q_\lambda(v)\le\tau\}\), the variance is of order \(\tau\). Dudley's entropy integral then gives a local fluctuation of order \(C_\lambda\sqrt{\frac{\mathrm{VC}_{\mathcal F}\tau\log n}{n}},\) and the local Rademacher deviation estimate, critical-radius calculation, and dyadic peeling argument convert this quantity into the fast term \(C_\lambda^2\frac{\mathrm{VC}_{\mathcal F}\log n}{n}.\) The complete derivation is given in \autoref{Appendix: Analysis of Quadrature Error}.

\section{Approximation Error Analysis}\label{Section: Approximation Error}
This section focuses on estimating the approximation error. We begin by introducing the definition of Sobolev functions on manifolds. Classical Sobolev spaces on manifolds admit several definitions. Under suitable manifold assumptions, Theorem 3 in~\cite{devito2019reproducingkernelhilbertspaces} establishes the equivalence of these definitions. Let $1\le p\le\infty$ and $s\in\mathbb N$. For a chart $(V_j,\psi_j)$ and a partition-of-unity function $\rho_j$, let $\widetilde{(\rho_jf)\circ\psi_j^{-1}}$ denote the zero extension from $\psi_j(V_j)$ to $\mathbb R^d$. We define
\begin{equation*}
\mathcal W_p^s(\mathcal M^d)
:=
\left\{f:\mathcal M^d\to\mathbb R
\left \vert
\|f\|_{\mathcal W_p^s(\mathcal M^d)}:=
\sum_{j=1}^k
\left\|
\widetilde{(\rho_jf)\circ\psi_j^{-1}}
\right\|_{W_p^s(\mathbb R^d)}<\infty
\right.
\right\}.
\end{equation*}
In the above definition, $\{\rho_j\}_{j=1}^k$ is a partition of unity on $\mathcal{M}^d$, and $\{\psi_j\}_{j=1}^k$ represents a collection of local coordinate mappings. Detailed definitions and properties of these mappings are provided in the~\autoref{Appendix: Notations}. It is worth noting that different choices of partition of unity and coordinate mappings lead to equivalent definitions of Sobolev spaces on manifolds. Prior work on neural network approximation theory for functions on manifolds has primarily focused on H{\"o}lder functions~\cite{Chen2019efficient, Chen2019NonparametricRO, Labate2024LowDA, SchmidtHieber2019DeepRN}. Several of these results formulate the H{\"o}lder class through ambient coordinates or the ambient Euclidean metric, which can make intrinsic verification less direct. For instance,~\cite{Labate2024LowDA} employs the following ambient formulation of the H{\"o}lder norm:
\begin{equation*}
    \|f\|_{\mathcal{H}^s_D(\mathcal{M}^d)} =  \max_{\alpha: \|\alpha\|_1 \le \lfloor s \rfloor} \sup_{x \in \mathcal{M}^d} |\partial^\alpha f(x)|\vee \max_{\alpha: \|\alpha\|_1 = \lfloor s\rfloor} \sup_{\substack{x, y \in \mathcal{M}^d \\ x \neq y}} \frac{|\partial^\alpha f(x) - \partial^\alpha f(y)|}{\|x - y\|_2^{s - \lfloor s\rfloor}} .
\end{equation*}
In contrast, the Sobolev norm we adopt is independent of the metric on the ambient Euclidean space. The following theorem illustrates the ability of DNNs to approximate Sobolev functions on manifolds.
\begin{theorem}\label{theorem: approximation error for sobolev functions on manifolds}
For \(1\le p<\infty,0< \varepsilon< \frac{1}{2} \) and  \(s\in\mathbb N_+\), let \(f\in\mathcal W_p^s(\mathcal M^d)\). Corresponding to \(B=1\) or \(B=\infty\), we can construct a neural network
\[
g\in \mathcal F(L,W,S,B,\infty)
\]
such that
\[
\|f-g\|_{L_p(\mathcal M^d)}\le \varepsilon.
\]
Furthermore,
\[
L\le C\log(\varepsilon^{-1}),\qquad
W\le C\varepsilon^{-d/s},\qquad
S\le C\varepsilon^{-d/s}\log(\varepsilon^{-1}),
\]
where \(C\) is independent of \(\varepsilon\). If \(f\in L_\infty(\mathcal M^d)\) and an output bound \(M_f\ge\|f\|_{L_\infty(\mathcal M^d)}\) is required, then the clipped network
\[
\mathrm{clip}_{M_f}(g)=\min\{M_f,\max\{-M_f,g\}\}
\]
belongs to an output-bounded ReLU class and satisfies the same \(L_p\)-error bound up to the same constant, since clipping is \(L_p\)-nonexpansive.
\end{theorem}

Notably, the algebraic exponent of the approximation tolerance depends only on the intrinsic dimension $d$. The constants and the input-layer size may depend on the fixed embedding into $\mathbb R^D$. Our result demonstrates that, within the manifold approximation setting, it is possible to achieve high approximation accuracy while avoiding an ambient-dimensional exponent in the approximation rate. Next, we establish a neural network approximation theorem tailored to entropy solutions.
\begin{theorem}\label{theorem: approximation error for entropy solutions}
Let \(u^*\) be the entropy solution in \autoref{definition: entropy solution}, and let \(M\ge \max\{\|u^*\|_{L_\infty},\|u_0\|_{L_\infty}\}\). Let \(\Xi\subset W^{1,\infty}_0(\mathcal M^d\times(0,T))\) be a nonnegative test class satisfying
\[
C_\Xi
:=
1+
\sup_{\xi\in\Xi}\left(
\|\partial_t\xi\|_{L_\infty(\mathcal M^d\times[0,T])}
+
3d\mathrm{Lip}_f\|\nabla_g\xi\|_{L_\infty(\mathcal M^d\times[0,T])}
\right)<\infty .
\]
Assume that the flux \(f\) is component-wise Lipschitz continuous with constant \(\mathrm{Lip}_f\). Then, for every \(0< \varepsilon< \frac{1}{2}\), corresponding to \(B=1\) or \(B=\infty\), there exists a neural network \(u\in\mathcal F(L,W,S,B,M)\) such that
\[
\|u^*-u\|_{L_1(\mathcal M^d\times[0,T])}\le \varepsilon,\qquad
\|u^*(\cdot,0)-u(\cdot,0)\|_{L_1(\mathcal M^d)}\le \varepsilon,
\]
and
\[
\mathcal R^{tb}(u)+
T\operatorname{Vol}_g(\mathcal M^d)
\sup_{\xi\in\Xi,\ c\in\mathcal C}\ppart{\mathcal R^{int}(u,\xi,c)}
\lesssim
C_\Xi\varepsilon .
\]
The network size satisfies
\[
L \le C\log(\varepsilon^{-1}),\qquad
W \le C\varepsilon^{-(d+1)},\qquad
S \le C\varepsilon^{-(d+1)}\log(\varepsilon^{-1}).
\]
The constant \(C\) is independent of \(\varepsilon\). Moreover, the construction can be chosen to preserve the TV-controlled temporal modulus of the entropy solution, as stated in \autoref{corollary: TV controlled temporal modulus}.
\end{theorem}

The explicit expressions for the constants are provided and verified in the proof in~\autoref{Appendix: Analysis of Approximation Error}. Moreover, the above construction can be refined to preserve the TV-controlled temporal modulus of the entropy solution, which is the regularity needed in the localization bias estimate.
\begin{corollary}\label{corollary: TV controlled temporal modulus}
    Under the assumptions of~\autoref{theorem: approximation error for entropy solutions}, the network constructed in the proof of~\autoref{theorem: approximation error for entropy solutions} can be chosen such that
    \begin{equation*}
        \omega_t(u;r):=\sup_{\substack{s,t\in[0,T], |t-s|\le r}}\|u(\cdot,t)-u(\cdot,s)\|_{L_1(\mathcal M^d)}\le C(1+TV(u_0))r,\qquad 0\le r\le T .
    \end{equation*}
    Consequently, for the localization bias appearing in~\autoref{theorem: L1 contraction for total error}, the constructed network satisfies
    \begin{equation}\label{equation: constructed approximation bias}
        b_{\lambda,T}(u)\le C_T(1+TV(u_0))(\alpha+\eta+h).
    \end{equation}
    The network size bounds remain
    \begin{equation*}
        L\le C\log(\varepsilon^{-1}),\qquad W\le C\varepsilon^{-(d+1)},\qquad S\le C\varepsilon^{-(d+1)}\log(\varepsilon^{-1}),
    \end{equation*}
    with a different constant \(C\) independent of \(\varepsilon\).
\end{corollary}

\section{\texorpdfstring{Localized \(L_1\)-Stability Analysis}{Localized L1-Stability Analysis}}\label{Section: Localized L1 contraction}
In this section, we identify the structural stability and coercivity mechanism underlying the wPINN loss on manifolds. This mechanism is induced by the entropy structure and the \(L_1\) contraction property of geometry-compatible conservation laws (see \eqref{equation: L1 contraction1 with bounded variation function}). Although the wPINN objective is not convex with respect to the neural network parameters, it is coercive in the solution space: small temporal boundary mismatch and small localized entropy violation force the trial function to be close to the entropy solution in terminal \(L_1\)-distance. Thus the localized entropy residual provides coercive control in the natural \(L_1\)-metric around the entropy solution.

This stability and coercivity property has two roles in our analysis. First, it provides the deterministic stability estimate that converts residual control into a terminal \(L_1\)-error bound. Second, it supplies the localization structure required for the fast-rate generalization analysis. More precisely, after decomposing the hypothesis class into residual sublevel sets, the stochastic fluctuation on each localized shell has variance proportional to the shell radius. This Bernstein-type localization is the mechanism that yields a complexity term of order \(\mathrm{VC}_{\mathcal F}/n\), rather than the slower square-root rate.

\begin{theorem}\label{theorem: L1 contraction for total error}
Let \(u\) be the entropy solution to \eqref{equation: hyperbolic conservation laws}--\eqref{equation: initial condition}, with initial data \(u_0\in L_\infty(\mathcal M^d)\cap BV(\mathcal M^d)\). Let \(v\in L_\infty(\mathcal M^d\times[0,T])\cap C([0,T];L_1(\mathcal M^d))\), and let \(\mathcal C\subset\mathbb R\) be a compact interval containing the essential range of the entropy solution \(u\).

Let \(\sigma_2(r):=(r_+)^2\). Define
\[
\rho(r):=\frac{15}{16}\sigma_2(1-r^2)=\frac{15}{16}(1-r^2)_+^2,\qquad
\rho_\eta(r):=\eta^{-1}\rho(r/\eta).
\]
Then \(\rho_\eta\ge0\), \(\rho_\eta\) is even, \(\operatorname{supp}\rho_\eta\subset[-\eta,\eta]\), and \(\int_{\mathbb R}\rho_\eta(r)\,dr=1\). Define
\[
\zeta(r):=2\sigma_2(r-1)-4\sigma_2\left(r-\frac32\right)+2\sigma_2(r-2).
\]
Then \(0\le \zeta\le1\), \(\zeta(r)=0\) for \(r\le1\), and \(\zeta(r)=1\) for \(r\ge2\). For \(0<h<T/4\) and each horizon \(\tau\in[4h,T]\), set
\[
\chi_{h,\tau}(s):=\zeta(s/h)\zeta((\tau-s)/h).
\]
Then \(\chi_{h,\tau}=1\) on \([2h,\tau-2h]\), \(\chi_{h,\tau}=0\) near \(0\) and \(\tau\), \(\operatorname{supp}\chi_{h,\tau}'\subset[h,2h]\cup[\tau-2h,\tau-h]\), and \(\|\chi_{h,\tau}'\|_{L_1(0,\tau)}\le C\).

Let \(\iota:\mathcal M^d\hookrightarrow\mathbb R^D\) be the fixed isometric embedding. Let \(\omega_d\) denote the Lebesgue measure of the unit ball in \(\mathbb R^d\), and define
\[
c_d:=\left(\int_{\mathbb R^d}(1-|z|^2)_+^2\,dz\right)^{-1}
=\frac{(d+2)(d+4)}{8\omega_d}.
\]
For \(0<\alpha\le\alpha_0\), with \(\alpha_0>0\) sufficiently small, define
\[
K_\alpha^{\mathrm{ext}}(x,y)
:=
c_d\alpha^{-d}
\sigma_2\left(1-\frac{|\iota(x)-\iota(y)|^2}{\alpha^2}\right).
\]
The required estimates for \(K_\alpha^{\mathrm{ext}}\) are given in \autoref{lemma: extrinsic ReQU kernel estimates}.

Let \(\lambda=(\alpha,\eta,h)\), \(0<\alpha\le\alpha_0\), \(0<\eta<h/4\). For each \(\tau\in[4h,T]\), define
\[
\Phi_{\lambda,\tau}
:=
\left\{
\Phi_{\lambda,\tau}^{y,s}(x,t)
:=
\chi_{h,\tau}\left(\frac{t+s}{2}\right)
\rho_\eta(t-s)
K_\alpha^{\mathrm{ext}}(x,y)
\,\middle|\,
(y,s)\in\mathcal M^d\times[0,\tau]
\right\},
\]
and
\[
\overline\Phi_\lambda:=\bigcup_{\tau\in[4h,T]}\Phi_{\lambda,\tau}.
\]
For $s=0$ or $s=\tau$, the support condition $|t-s|\le\eta<h/4$ forces the time cutoff to vanish, so the corresponding test function is identically zero. Thus adjoining the endpoints does not change either supremum below, while it makes the center--horizon parameter set compact.
Define the horizon-indexed and global positive localized entropy violations by
\[
\mathcal L_{\lambda,\tau}^+(v)
:=
\sup_{\xi\in\Phi_{\lambda,\tau},\ c\in\mathcal C}
\ppart{\mathcal R^{int}_{[0,\tau]}(v,\xi,c)},
\qquad
\mathfrak L_\lambda^+(v)
:=
\sup_{\xi\in\overline\Phi_\lambda,\ c\in\mathcal C}
\ppart{\mathcal R^{int}(v,\xi,c)}.
\]
Here \(\mathcal R^{int}_{[0,\tau]}\) denotes the signed residual integrated over \(\mathcal M^d\times[0,\tau]\).

Then for every \(\tau\in[4h,T]\),
\begin{equation}\label{equation: horizon terminal stability}
\|v(\cdot,\tau)-u(\cdot,\tau)\|_{L_1(\mathcal M^d)}
\le
\|v(\cdot,0)-u_0\|_{L_1(\mathcal M^d)}
+
\tau\operatorname{Vol}_g(\mathcal M^d)\mathcal L_{\lambda,\tau}^+(v)
+
b_{\lambda,\tau}(v),
\end{equation}
where
\[
b_{\lambda,\tau}(v)
\le
C_T\left[
\omega_t(v;C\eta)+\omega_t(v;Ch)+(\eta+h)TV(u_0)+\alpha(1+TV(u_0))
\right].
\]
Consequently,
\begin{equation}\label{equation: terminal stability global positive}
\|v(\cdot,T)-u(\cdot,T)\|_{L_1(\mathcal M^d)}
\le
\mathcal R^{tb}(v)
+
T\operatorname{Vol}_g(\mathcal M^d)\mathfrak L_\lambda^+(v)
+
b_{\lambda,T}(v).
\end{equation}
Moreover, the following space-time stability estimate holds:
\begin{equation}\label{equation: spacetime stability global positive}
\begin{aligned}
&\frac{1}{T\operatorname{Vol}_g(\mathcal M^d)}
\int_0^T
\|v(\cdot,t)-u(\cdot,t)\|_{L_1(\mathcal M^d)}\,dt
\\
&\qquad\le
C_T\left[
\mathcal R^{tb}(v)
+T\operatorname{Vol}_g(\mathcal M^d)\mathfrak L_\lambda^+(v)
+b_{\lambda,T}(v)
\right].
\end{aligned}
\end{equation}
In particular, if
\[
\omega_t(v;r)\le C(1+TV(u_0))r,\qquad 0\le r\le T,
\]
then
\[
b_{\lambda,T}(v)\le C_T(1+TV(u_0))(\alpha+\eta+h).
\]
Here \(C_T\) may depend on \(T\), the fixed geometry and embedding, uniform \(L_\infty\)-bounds for \(u\) and \(v\), the compact interval \(\mathcal C\), and the relevant flux norms, but it is independent of \(\alpha,\eta,h\).
\end{theorem}

The proof is provided in~\autoref{Appendix: L1 Contraction for Total Error}.

\section{Quadrature Error Analysis}\label{Section: Quadrature Error}

This section states the fast-rate quadrature theorem. The full proof is deferred to \autoref{Appendix: Analysis of Quadrature Error}, where the local Rademacher deviation estimate, the empirical \(L_2\)-localization, the Dudley entropy integral, the critical radius, and the dyadic shell argument are written out in detail.
Throughout this section, all empirical-process classes are assumed to be pointwise measurable. Equivalently, all probabilities and expectations may be interpreted in the outer sense.

Let \(\mathcal Z:=\mathcal M^d\times[0,T], A_T:=T\operatorname{Vol}_g(\mathcal M^d),A_0:=\operatorname{Vol}_g(\mathcal M^d).\) Let \(P\) be the uniform probability measure on \(\mathcal Z\), \(P_n\) the empirical measure generated by the interior quadrature points, and let \(P_0,P_{0,n}\) be the corresponding population and empirical measures on \(\mathcal M^d\) at the initial time. Then
\[
   \mathcal R^{int}(v,\xi,c)=A_TP r^{int}(v,\xi,c),
   \qquad
   \mathcal R_n^{int}(v,\xi,c)=A_TP_n r^{int}(v,\xi,c),
\]
and
\[
   \mathcal R^{tb}(v)=A_0P_0|v(\cdot,0)-u_0|,
   \qquad
   \mathcal R_n^{tb}(v)=A_0P_{0,n}|v(\cdot,0)-u_0|.
\]

Write the adversarial localized ReQU family as
\[
   \overline\Phi_\lambda=\{\xi_\vartheta:\vartheta\in\Theta_\lambda\},
   \qquad
   \Theta_\lambda
   :=
   \left\{(y,s,\tau):
   y\in\mathcal M^d,\ 4h\le\tau\le T,\ 0\le s\le\tau
   \right\}.
\]
The parameter set \(\Theta_\lambda\) is compact and finite-dimensional. The endpoint values \(s=0\) and \(s=\tau\) produce the zero test function, so including them does not change the adversarial supremum. For \(\vartheta,\vartheta'\in\Theta_\lambda\), define
\begin{align*}
 d_\lambda(\vartheta,\vartheta')
 :={}
 \|\partial_t(\xi_\vartheta-\xi_{\vartheta'})\|_{L_\infty(\mathcal Z)}
 +3d\mathrm{Lip}_f
 \|\nabla_g(\xi_\vartheta-\xi_{\vartheta'})\|_{L_\infty(\mathcal Z)}.
\end{align*}
The localization constant is enlarged, if necessary, so that \(C_\lambda \ge 1+\sup_{\vartheta\in\Theta_\lambda}(\|\partial_t\xi_\vartheta\|_{L_\infty(\mathcal Z)}+3d\mathrm{Lip}_f\|\nabla_g\xi_\vartheta\|_{L_\infty(\mathcal Z)}),\) and so that it also dominates the fixed test-parameter Lipschitz constants used below. We assume the fixed-architecture metric-entropy estimate
\begin{equation}\label{equation: section5 fixed test entropy}
   \log N(\delta,\Theta_\lambda,d_\lambda)
   \le q_\lambda\log\left(\frac{A_\lambda}{\delta}\right),
   \qquad 0<\delta\le1,
\end{equation}
where \(q_\lambda,A_\lambda\) do not depend on \(n,L,W,S\). For the explicit localized ReQU family in \autoref{theorem: L1 contraction for total error}, this follows from the compact finite-dimensional parametrization by the spatial center, time center, and terminal horizon. The test parameters are therefore optimized adversarially, but their complexity does not grow with the sample size. The only growing network complexity is \(\mathrm{VC}_{\mathcal F}:=SL\log S.\)

Let \(B_{\lambda,\mathcal F}:=b_{\lambda,T}(\mathcal F),\) and define the population and empirical stability functionals
\begin{align}
\mathfrak Q_\lambda(v)
&:=
\mathcal R^{tb}(v)
+A_T\mathfrak L_\lambda^+(v)
+B_{\lambda,\mathcal F},
\label{equation: section5 population Q}
\\
\mathfrak Q_{\lambda,n}(v)
&:=
\mathcal R_n^{tb}(v)
+A_T\mathfrak L_{\lambda,n}^+(v)
+B_{\lambda,\mathcal F},
\label{equation: section5 empirical Q}
\end{align}
where \(\mathfrak L_{\lambda,n}^+(v) := \sup_{\xi\in\overline\Phi_\lambda,\ c\in\mathcal C} \ppart{\mathcal R_n^{int}(v,\xi,c)}.\) Since \(B_{\lambda,\mathcal F}\) is constant on the class, it does not alter the empirical minimizer. Let
\[
   u_{\mathcal F}\in\arg\min_{v\in\mathcal F}\mathfrak Q_\lambda(v),
   \qquad
   u_n\in\arg\min_{v\in\mathcal F}\mathfrak Q_{\lambda,n}(v).
\]
For \(\tau>0\), set \(\mathcal F_\tau:=\{v\in\mathcal F:\mathfrak Q_\lambda(v)\le\tau\}.\) The space-time stability estimate \eqref{equation: spacetime stability global positive} implies \(P|v-u^*|\lesssim\tau,v\in\mathcal F_\tau.\) For \(h_{v,\xi,c}:=r^{int}(v,\xi,c)-r^{int}(u^*,\xi,c),\)
the pointwise Lipschitz estimate gives the Bernstein condition \(P h_{v,\xi,c}^2 \le2MC_\lambda^2P|v-u^*| \lesssim C_\lambda^2\tau.\) The joint residual covering estimate and Dudley's entropy integral then yield
\[
   \mathrm{Rad}_n(\mathcal H_\tau)
   \lesssim
   C_\lambda
   \sqrt{\frac{\mathrm{VC}_{\mathcal F}\tau\log n}{n}},
\]
where \(\mathcal H_\tau\) is the localized joint residual-increment class. The local Rademacher deviation estimate of~\cite{bartlett2005local,mendelson2003few} converts this local complexity bound into a deviation inequality with variance proportional to \(\tau\). A critical-radius calculation and dyadic peeling then produce the \(\mathrm{VC}_{\mathcal F}/n\) rate.

For the final fast-rate statement, define \(\Gamma_n:=\mathrm{VC}_{\mathcal F}\log n.\) The lower bound in the following theorem ensures that the fixed exact-solution calibration process is dominated by the localized solution-dependent term; the upper bound ensures that the confidence parameter remains of order at most \(\Gamma_n\).

\begin{theorem}\label{theorem: oracle inequality}
Let \(\mathcal F=\mathcal F(L,W,S,B,M)\), assume \(B_{\lambda,\mathcal F}<\infty\), and suppose that the fixed adversarial test family satisfies \eqref{equation: section5 fixed test entropy}. Assume that, for constants \(c_0,c_1>0\),
\begin{equation}\label{equation: fast regime Gamma}
   c_0\sqrt{n\log n}
   \le \Gamma_n
   \le c_1n,
\end{equation}
and that \(n\ge C\mathrm{VC}_{\mathcal F}\) for a sufficiently large constant \(C\). Then there exist constants \(c,C'>0\), independent of \(n,L,W,S\), such that with probability at least
\begin{equation}\label{equation: quadrature fast probability}
   1-\exp\left(-c\frac{\Gamma_n^2}{n}\right),
\end{equation}
there holds
\begin{equation}\label{equation: theorem6 VCF oracle fast}
\mathfrak Q_\lambda(u_n)
\le
C'\left[
\mathfrak Q_\lambda(u_{\mathcal F})
+C_\lambda^2\frac{\Gamma_n}{n}
\right].
\end{equation}
The hidden constant may depend on \(M,T,\mathcal M^d\), the flux bounds, and the fixed adversarial test architecture, but not on \(n,L,W,S\). The same conclusion, with the same upper bound \(\mathrm{VC}_{\mathcal F}=SL\log S\), holds when \(\mathcal F\) is replaced throughout by any subclass \(\mathcal G\subseteq\mathcal F(L,W,S,B,M)\), with \(B_{\lambda,\mathcal G}\) and the corresponding population and empirical minimizers taken over \(\mathcal G\).
\end{theorem}

The proof first derives a general localized oracle inequality with an arbitrary confidence parameter. It then chooses \(t_n\asymp \frac{\Gamma_n^2}{n}\) and uses \eqref{equation: fast regime Gamma} to absorb every remaining stochastic contribution into \(C_\lambda^2\Gamma_n/n\). Thus the final quadrature estimate contains only the fast-rate term displayed in \eqref{equation: theorem6 VCF oracle fast}.

\section{Proof of Main Results}\label{Section: Proof of Main Result}

\begin{proof}[Proof of \autoref{theorem: main result}]
Fix an approximation tolerance \(0<\varepsilon<1/2\). By \autoref{theorem: approximation error for entropy solutions}, applied to the bounded nonnegative test class \(\overline\Phi_\lambda\), there exists a network \(v_\varepsilon\in\mathcal F(L,W,S,B,M)\) such that
\begin{equation}\label{equation: main proof approximation J abstract}
   \mathcal R^{tb}(v_\varepsilon)
   +A_T\mathfrak L_\lambda^+(v_\varepsilon)
   \lesssim C_\lambda\varepsilon,
\end{equation}
and
\begin{equation}\label{equation: main proof network size abstract}
L\lesssim\log(\varepsilon^{-1}),\qquad
W\lesssim\varepsilon^{-(d+1)},\qquad
S\lesssim\varepsilon^{-(d+1)}\log(\varepsilon^{-1}).
\end{equation}
Let \(u_{\mathcal F}\) be a population minimizer of \(\mathfrak Q_\lambda\) over \(\mathcal F\). Since the class-wise bias is constant over \(\mathcal F\), \eqref{equation: main proof approximation J abstract} gives
\begin{equation}\label{equation: main proof approximation Q abstract}
\mathfrak Q_\lambda(u_{\mathcal F})
\le \mathfrak Q_\lambda(v_\varepsilon)
\lesssim C_\lambda\varepsilon+b_{\lambda,T}(\mathcal F).
\end{equation}

Choose constants \(c_L,c_W,c_S>0\), depending only on the approximation construction, sufficiently large that \(v_\varepsilon\) belongs to the class with budgets
\begin{equation}\label{equation: main proof exact architecture budgets}
\begin{aligned}
L&=\left\lceil c_L\log(\varepsilon^{-1})\right\rceil,\\
W&=\left\lceil c_W\varepsilon^{-(d+1)}\right\rceil,\\
S&=\left\lceil c_S\varepsilon^{-(d+1)}\log(\varepsilon^{-1})\right\rceil.
\end{aligned}
\end{equation}
Unused units and parameters can be set to zero, so enlarging these budgets does not alter the approximation network. Consequently,
\[
\mathrm{VC}_{\mathcal F}=SL\log S
\asymp
\varepsilon^{-(d+1)}\left(\log\frac1\varepsilon\right)^3.
\]
Choose
\begin{equation}\label{equation: main proof epsilon choice}
   \varepsilon=n^{-a}(\log n)^{4a},
   \qquad a=1/(d+2).
\end{equation}
Then
\begin{equation*}
   \Gamma_n=\mathrm{VC}_{\mathcal F}\log n
   \asymp n^{1-a}(\log n)^{4a}
\end{equation*}
up to constants depending only on the approximation construction. In particular, for all sufficiently large \(n\), \(c_0\sqrt{n\log n}\le\Gamma_n\le c_1n,\) and hence the assumptions of \autoref{theorem: oracle inequality} are satisfied. Moreover,
\begin{equation}\label{equation: main proof fast term balance}
   \frac{\Gamma_n}{n}
   \lesssim n^{-a}(\log n)^{4a}.
\end{equation}
Therefore, with probability at least
\[
1-\exp\left(-c\frac{\Gamma_n^2}{n}\right)
\ge
1-\exp\left(-c'n^{1-2a}(\log n)^{8a}\right),
\]
\autoref{theorem: oracle inequality}, \eqref{equation: main proof approximation Q abstract}, and \eqref{equation: main proof fast term balance} yield
\begin{align*}
\mathfrak Q_\lambda(u_n)
&\lesssim
\mathfrak Q_\lambda(u_{\mathcal F})
+C_\lambda^2\frac{\Gamma_n}{n}
\\
&\lesssim
C_\lambda\varepsilon
+b_{\lambda,T}(\mathcal F)
+C_\lambda^2n^{-a}(\log n)^{4a}
\\
&\lesssim
C_\lambda^2n^{-a}(\log n)^{4a}
+b_{\lambda,T}(\mathcal F),
\end{align*}
where we used \(C_\lambda\ge1\). Finally, by \eqref{equation: terminal stability global positive} and the definition of the class-wise bias,
\[
\|u_n(\cdot,T)-u^*(\cdot,T)\|_{L_1(\mathcal M^d)}
\le
\mathcal R^{tb}(u_n)+A_T\mathfrak L_\lambda^+(u_n)+b_{\lambda,T}(u_n)
\le
\mathfrak Q_\lambda(u_n).
\]
This proves \eqref{equation: main theorem bias rate}. Substituting \eqref{equation: main proof epsilon choice} into the explicit budget choices in \eqref{equation: main proof exact architecture budgets} gives the stated network-size bounds.
\end{proof}

\begin{proof}[Proof of \autoref{corollary: explicit bias convergence}]
Apply the preceding proof with
\[
   \mathcal F=\mathcal F_{\rm tm}(L,W,S,B,M).
\]
By \autoref{corollary: TV controlled temporal modulus}, the approximation network \(v_\varepsilon\) belongs to this restricted class, so the approximation estimate and network-size bounds are unchanged. By the definition of \(\mathcal F_{\rm tm}\),
\begin{align*}
b_{\lambda,T}(\mathcal F_{\rm tm})
&\le C_T\Bigl[
C_{\rm tm}(1+TV(u_0))(\eta+h)
\\[-1mm]
&\hspace{34mm}
+(\eta+h)TV(u_0)+\alpha(1+TV(u_0))
\Bigr]
\\
&\lesssim C_T(1+TV(u_0))(\alpha+\eta+h).
\end{align*}
Substitution into \eqref{equation: main theorem bias rate} proves \eqref{equation: explicit bias convergence rate}.
\end{proof}

\begin{proof}[Proof of \autoref{corollary: sobolev convergence rate}]
Under the pointwise-in-time Sobolev hypothesis in \autoref{corollary: sobolev convergence rate}, the nodewise spatial approximation in the construction of \autoref{theorem: approximation error for entropy solutions} uses \autoref{theorem: approximation error for sobolev functions on manifolds} with order \(s\). The uniform bound makes all nodewise approximation constants independent of the time node. It also gives
\[
TV(u^*(\cdot,t))
\le C\|u^*(\cdot,t)\|_{\mathcal W_1^1(\mathcal M^d)}
\le C\sup_{0\le r\le T}\|u^*(\cdot,r)\|_{\mathcal W_1^s(\mathcal M^d)},
\qquad 0\le t\le T.
\]
Thus, in the spatial-localization estimates in the proof of \autoref{theorem: L1 contraction for total error}, the exponential BV-growth bound \eqref{equation: bounds on TV} is replaced by this uniform Sobolev bound. For every \(0<\varepsilon<1/2\), the construction yields a network \(v_\varepsilon\) satisfying
\[
   \mathcal R^{tb}(v_\varepsilon)+A_T\mathfrak L_\lambda^+(v_\varepsilon)
   \lesssim C_\lambda\varepsilon,
\]
with
\begin{equation}\label{equation: sobolev proof network sizes}
L\lesssim\log(\varepsilon^{-1}),\qquad
W\lesssim\varepsilon^{-(d/s+1)},\qquad
S\lesssim\varepsilon^{-(d/s+1)}\log(\varepsilon^{-1}).
\end{equation}
Choose constants \(c_{L,s},c_{W,s},c_{S,s}>0\) sufficiently large and take
\[
L=\left\lceil c_{L,s}\log(\varepsilon^{-1})\right\rceil,
\qquad
W=\left\lceil c_{W,s}\varepsilon^{-(d/s+1)}\right\rceil,
\qquad
S=\left\lceil c_{S,s}\varepsilon^{-(d/s+1)}\log(\varepsilon^{-1})\right\rceil.
\]
As above, the approximation network is embedded in this class by setting unused units and parameters to zero, and hence
\[
\mathrm{VC}_{\mathcal F}
\asymp
\varepsilon^{-(d/s+1)}\left(\log\frac1\varepsilon\right)^3.
\]
Set
\[
   \varepsilon=n^{-a_s}(\log n)^{4a_s},
   \qquad a_s=s/(2s+d).
\]
Then
\[
   \Gamma_n=\mathrm{VC}_{\mathcal F}\log n
   \asymp n^{1-a_s}(\log n)^{4a_s},
   \qquad
   \frac{\Gamma_n}{n}
   \lesssim n^{-a_s}(\log n)^{4a_s}.
\]
Since \(0<a_s<1/2\), the fast-regime condition \eqref{equation: fast regime Gamma} holds for all sufficiently large \(n\). Applying \autoref{theorem: oracle inequality} exactly as in the proof of \autoref{theorem: main result} gives
\[
\mathfrak Q_\lambda(u_n)
\lesssim
C_\lambda^2n^{-a_s}(\log n)^{4a_s}
+b_{\lambda,T}(\mathcal F)
\]
with probability at least \(1-\exp\left(-c n^{1-2a_s}(\log n)^{8a_s}\right).\) The terminal stability estimate proves \eqref{equation: sobolev rate abstract bias}, and the explicit Sobolev budget choices above give the stated architecture.

If the empirical minimization is performed over \(\mathcal F_{\rm tm}\), the same approximation network remains admissible and
\[
   b_{\lambda,T}(\mathcal F_{\rm tm})
   \lesssim C_T(1+TV(u_0))(\alpha+\eta+h).
\]
This proves \eqref{equation: sobolev rate explicit bias}.
\end{proof}

\section{Numerical Experiments}\label{Section: Numerical Experiment}
This section presents the simulation results and implementation details of the wPINNs algorithm (see Algorithm~\ref{algorithm: training algorithm}). The wPINNs algorithm is implemented using Python scripts powered by the PyTorch framework \href{https://pytorch.org/}{https://pytorch.org/}. The scripts can be downloaded from \href{https://github.com/hanfei27/wPINNonmaifolds}{https://github.com/hanfei27/wPINNonmaifolds}. While the preceding theoretical analysis focuses on boundaryless manifolds, the algorithm can be extended to manifolds with boundaries, albeit with unresolved theoretical challenges.  We briefly describe the algorithm and present numerical experiments specifically for the case of the Dirichlet boundary condition, leaving convergence analysis for the boundary case to future work. The core of the algorithm lies in the loss function, which consists of three previously introduced components: the internal loss, the temporal boundary loss, and the spatial boundary loss. These loss terms are reformulated below, with the key distinction that the test functions are now represented by neural networks parameterized by $\eta$, constrained via multiplication with a compactly supported function. The specific form of loss functions is described as
\begin{equation*}
    L^{int}(u_\theta,\xi_\eta,c) = \int_{\mathcal{M}^d \times (0,T)}-\widetilde{U}(u_\theta, c)  \partial_t\xi_\eta(t,  x) - g_x\bigl(\widetilde{F}(u_\theta, c),  \nabla_g\xi_\eta(t,  x)\bigr) dV_g(x)dt,
\end{equation*}
\begin{equation*}
    L^{tb}(u_\theta) = \int_{\mathcal{M}^d}|u_\theta(x, 0)- u(x, 0)|dV_g(x),
\end{equation*}
\begin{equation*}
    L^{sb}(u_\theta)  = \int_{\partial \mathcal{M}^d \times [0,T]} |u_\theta(x,t)-u(x,t)|\,dV_{\partial g}(x)dt.
\end{equation*}
Our theoretical analysis demonstrates that a judicious choice of the test function space can substantially enhance training efficiency. Guided by this insight, define
\begin{align*}
\mathcal N(u_\theta,\xi_\eta,c)
&:=\int_{\mathcal M^d\times(0,T)}
\Bigl[-\widetilde U(u_\theta,c)\,\partial_t\xi_\eta
-g_x\bigl(\widetilde F(u_\theta,c),\nabla_g\xi_\eta\bigr)\Bigr]
\,dV_g\,dt,
\\
\mathcal D(\xi_\eta)
&:=\varepsilon_{\rm den}
+\int_{\mathcal M^d\times(0,T)}
\Bigl(\xi_\eta^2+\|\nabla_g\xi_\eta\|_g^2\Bigr)
\,dV_g\,dt,
\end{align*}
and use the normalized internal loss
\begin{equation*}
L^{int}(u_\theta,\xi_\eta,c)
:=\frac{[\mathcal N(u_\theta,\xi_\eta,c)]_+^2}
{\mathcal D(\xi_\eta)}.
\end{equation*}
Here, $\varepsilon_{\rm den}>0$ is a fixed denominator stabilizer, and squaring the positive part of the numerator enhances numerical stability during optimization.  For weak adversarial networks, the choice of the test function space is crucial. As demonstrated in~\cite{Chen2020ACS}, applying soft or hard constraints to enforce zero boundary conditions can significantly accelerate training and improve accuracy. We implement such constraints using a cutoff function—specifically, a compactly supported function $\omega : \mathcal{M}^d \times (0,T) \to \mathbb{R}_+$ that vanishes on $\partial \mathcal{M}^d$ while remaining positive in the interior of $\mathcal{M}^d$. Consequently, the test functions are defined as $\xi_\eta(x,t) = \omega(x,t)\tilde{\xi}_\eta(x,t)$, where $\tilde{\xi}_\eta$ is a neural network parameterized by $\eta$.

Additionally, the work in~\cite{chaumet2022efficient} proposed solving a dual problem to select the test functions, achieving promising numerical performance. A hyperparameter $\varrho$ is used to balance the internal and boundary loss terms. The total loss is defined as
\begin{equation*}
    L(u_\theta,\xi_\eta,c) =   L^{int} + \varrho (L^{tb} + L^{sb}).
\end{equation*}
In practice, these integrals are approximated using numerical quadrature methods, such as Monte Carlo integration. For the parameter $c$, points $\{c_i\}_{i=1}^{N_c}$ are sampled uniformly from the interval $[c_{\min}, c_{\max}]$, and the maximum loss value $\max_{c_i, i=1,\ldots,N_c} L(u_\theta,\xi_\eta,c_i)$ is used to enforce the entropy condition.

\begin{algorithm}[htb]
\caption{wPINN algorithm on manifolds}\label{algorithm: training algorithm}
\KwIn{Initial data $u_0$, boundary data $g$, hyperparameter: $N_{int},N_{sb},N_{ini},N_{\min},N_{\max},N_c,N_{ep},\tau_{\min},\tau_{\max},\varrho,r$, neural network $u_\theta, \xi_\eta$ }
\KwOut{Best networks $u_\theta^b$}
Initialize the networks  $u_\theta, \xi_\eta$ \;
Generate interior points $S^{int} = \{(x_i,t_i)\}_{i=1}^{N_{int}} \subset \mathcal{M}^d\times(0,T)  $, boundary points $S^{sb} = \{(x_i,t_i)\}_{i=1}^{N_{sb}}\subset\partial \mathcal{M}^d  \times[0,T]  $, and initial points $S^{ini} = \{x_i\}_{i=1}^{N_{ini}} \subset \mathcal{M}^d $\;
Initialize performance loss $L^b \leftarrow \infty $\\
\For{$i=1,...,N_{ep}$}{
    \If{$i \% r = 0$ }{Randomly initialize $\eta$\;}
    \For{$N=1,...,N_{\max}$}{
    Compute $L_{\max}(u_\theta,\xi_\eta) = \max_{i=1}^{N_c} L(u_\theta,\xi_\eta,c_i)$\;
    Update $\eta \leftarrow \eta + \tau_{\max}\nabla_\eta L_{\max}(u_\theta,\xi_\eta)$\;
    }
    \For{$N=1,...,N_{\min}$}{
    Compute $L_{\max}(u_\theta,\xi_\eta) = \max_{i=1}^{N_c} L(u_\theta,\xi_\eta,c_i)$\;
    Update $\theta \leftarrow \theta - \tau_{\min}\nabla_\theta L_{\max}(u_\theta,\xi_\eta)$\;
    }
    \If{$L_{\max} < L^b$}{
    $L^b\leftarrow L_{\max}$\;
    Save best networks $u_\theta^b$;
    }
}
\end{algorithm}
We performed ensemble training as described in~\cite{lye2020deep}. A hyperparameter $r$ is introduced to periodically reinitialize the test function, which enhances the training stability. After retraining the model $N$ times, we obtain a set of networks $\{u^b_{\theta_i}\}_{i=1}^N$. The final output function is defined as the average of these networks:
\begin{equation}\label{equation: retrain}
    u_{\text{avg}}(x,t) = \frac{1}{N}\sum_{i=1}^N u^b_{\theta_i}(x,t)  .
\end{equation}

This ensemble averaging approach improves both the accuracy and stability of the resulting output. To assess accuracy, we compute the relative space-time $L_1$ test error:
\begin{equation*}
    \mathcal{E}_{\mathrm{st}}(u_\theta) = \frac{\int_{ \mathcal{M}^d\times [0,T]} |u_\theta(x, t) - u(x, t)| \, dV_g(x) \, dt}{\int_{\mathcal{M}^d\times [0,T]} |u(x, t)| \, dV_g(x) \, dt}.
\end{equation*}

In the next section, we present numerical experiments conducted on the sphere. The experimental setup is detailed there.

\subsection{Experimental Setup}\label{Subsection: Experimental Setup}
Numerical experiments are conducted on the unit sphere $\mathbb{S}^2$; the three nonperiodic cases are implemented on longitude--latitude patches with imposed boundary data, whereas the periodic case covers the full longitude range. This choice of domain is not only mathematically illustrative but also practically significant—particularly in the context of solving the shallow water equations on the sphere. These equations play a crucial role in geophysical modeling~\cite{bauer2015quiet} and pose considerable computational challenges due to their nonlinear hyperbolic nature and the curvature of the spherical domain. Achieving accurate numerical solutions requires careful consideration of the sphere’s global geometry and the potential emergence of discontinuities in the flow. The unit sphere is parameterized by the normalized longitude $\lambda\in[-1,1]$ and latitude $\phi\in[-1,1]$ as
\begin{equation*}
\begin{cases}
 x_1=\cos\left(\frac{\pi}{2}\phi\right)\cos(\pi\lambda),\\
 x_2=\cos\left(\frac{\pi}{2}\phi\right)\sin(\pi\lambda),\\
 x_3=\sin\left(\frac{\pi}{2}\phi\right).
\end{cases}
\end{equation*}
The corresponding unit tangent vectors are
\begin{equation*}
\begin{aligned}
\mathbf i_\lambda
&=-\sin(\pi\lambda)\mathbf i_1+\cos(\pi\lambda)\mathbf i_2,\\
\mathbf i_\phi
&=-\sin\left(\frac{\pi}{2}\phi\right)\cos(\pi\lambda)\mathbf i_1
  -\sin\left(\frac{\pi}{2}\phi\right)\sin(\pi\lambda)\mathbf i_2
  +\cos\left(\frac{\pi}{2}\phi\right)\mathbf i_3.
\end{aligned}
\end{equation*}
Here $\mathbf i_1,\mathbf i_2,\mathbf i_3$ are the standard basis vectors in $\mathbb R^3$. Writing \(f=f_\lambda\mathbf i_\lambda+f_\phi\mathbf i_\phi,\) the surface divergence and gradient are
\begin{equation*}
\mathbf{div}_g f
=
\frac{1}{\cos(\frac{\pi}{2}\phi)}
\left[
\frac{1}{\pi}\frac{\partial f_\lambda}{\partial\lambda}
+
\frac{2}{\pi}\frac{\partial}{\partial\phi}
\left(f_\phi\cos\left(\frac{\pi}{2}\phi\right)\right)
\right],
\end{equation*}
and
\begin{equation*}
\nabla_g\xi
=
\frac{1}{\pi\cos(\frac{\pi}{2}\phi)}
\frac{\partial\xi}{\partial\lambda}\mathbf i_\lambda
+
\frac{2}{\pi}\frac{\partial\xi}{\partial\phi}\mathbf i_\phi.
\end{equation*}
Let
\[
\widetilde F(u_\theta,c)
=
\widetilde F_\lambda(u_\theta,c)\mathbf i_\lambda
+
\widetilde F_\phi(u_\theta,c)\mathbf i_\phi.
\]
Then the internal residual used in the experiments is
\begin{equation*}
\begin{aligned}
R^{int}(u_\theta,\xi_\eta,c)
={}&\int_0^T\int_{\mathbb S^2}
-\widetilde U(u_\theta,c)\partial_t\xi_\eta
-
\frac{1}{\pi\cos(\frac{\pi}{2}\phi)}
\partial_\lambda\xi_\eta\,\widetilde F_\lambda(u_\theta,c)
\\
&
-
\frac{2}{\pi}
\partial_\phi\xi_\eta\,\widetilde F_\phi(u_\theta,c)
\,dSdt.
\end{aligned}
\end{equation*}

As shown in~\cite{Amorim2005HYPERBOLICCL}, every smooth vector field on $\mathbb S^2$ can be represented as
\begin{equation*}
f(x,u)=n(x)\times\Phi(x,u).
\end{equation*}
If $\Phi(x,u)=\Phi(u)$ is independent of $x$, then the flux is geometry-compatible~\cite{BenArtzi2008HyperbolicCL}. Writing
\[
\Phi(u)=f_1(u)\mathbf i_1+f_2(u)\mathbf i_2+f_3(u)\mathbf i_3,
\]
one obtains
\begin{equation*}
\begin{aligned}
f_\lambda(\lambda,\phi,u)
&=
f_1(u)\sin\left(\frac{\pi}{2}\phi\right)\cos(\pi\lambda)
+f_2(u)\sin\left(\frac{\pi}{2}\phi\right)\sin(\pi\lambda)
-f_3(u)\cos\left(\frac{\pi}{2}\phi\right),\\
f_\phi(\lambda,\phi,u)
&=-f_1(u)\sin(\pi\lambda)+f_2(u)\cos(\pi\lambda).
\end{aligned}
\end{equation*}
For $f_1=f_2\equiv0$, the conservation law reduces on each fixed latitude to
\begin{equation}\label{equation: HCL-simple}
\frac{\partial u}{\partial t}
-
\frac{1}{\pi}\frac{\partial}{\partial\lambda}f_3(u)
=0.
\end{equation}
We choose \(f_3(u)=-\frac{\pi}{2}u^2,\) so that \eqref{equation: HCL-simple} becomes the inviscid Burgers equation \(\partial_tu+\partial_\lambda\left(\frac{u^2}{2}\right)=0.\) Let $\widetilde u(\lambda,t)$ be a periodic entropy solution of this one-dimensional equation and set
\[
\widehat u(\phi):=\mathbf 1_{[-1/2,1/2]}(\phi),
\qquad
u(\lambda,\phi,t):=\widetilde u(\lambda,t)\widehat u(\phi).
\]
Because $\widehat u$ takes only the values $0$ and $1$, $f_3(0)=0$, and the flux has no $\phi$-component, the equation is satisfied independently on almost every latitude. Thus $u$ gives the band-supported spherical solution used in the experiments.

\subsection{Simulation Results}\label{Subsection: Simulation Result}
We will conduct numerical experiments for four scenarios: standing shock, moving shock, rarefaction wave, and sine wave. The sine wave scenario will employ periodic boundary conditions, requiring the simulation to cover the entire longitudinal range of the sphere. In contrast, the other three scenarios (standing shock, moving shock, and rarefaction wave) will be simulated on a defined longitudinal segment of the sphere, without periodic boundary conditions. In each experiment, we use $\widehat u(\phi)=\mathbf 1_{[-1/2,1/2]}(\phi)$. Two entropy functions are tested: the Kruzkov entropy $\widetilde U_1(u,c)=|u-c|$ and the quadratic entropy $\widetilde U_2(u,c)=\frac12(u-c)^2$. Their corresponding scalar Burgers entropy fluxes are $\widetilde F_1(u,c)=\mathrm{sgn}(u-c)(\frac{u^2}{2}-\frac{c^2}{2})$ and $\widetilde F_2(u,c)=\frac13u^3-\frac12cu^2$, respectively, up to an irrelevant additive constant in $c$. Each experiment will perform a grid search over a subset of hyperparameters to identify the optimal settings. Retraining will be carried out according to equation~\eqref{equation: retrain} for different network initializations to ensure the robustness of the model.~\autoref{table:loss compare} compares the space-time test errors of the average predicted solutions obtained using the two entropy losses and lists the number of collocation points used in each of the four experiments.
\begin{table}[htbp]
    \centering
    \caption{Configurations and relative space-time test errors for the four experiments. $\mathcal{E}_{\mathrm{st}}$-K and $\mathcal{E}_{\mathrm{st}}$-S denote the relative space-time test errors of the average predicted solutions obtained with the Kruzkov and quadratic entropies, respectively.}
    \label{table:loss compare}
    \vspace{5pt}
    \begin{tabular}{l|l|l|l|l|l}
      \textbf{} & $N_{int}$ & $N_{sb}$&$N_{ini}$ & \textbf{$\mathcal{E}_{\mathrm{st}}$}-K &\textbf{$\mathcal{E}_{\mathrm{st}}$}-S\\
      \hline
      Standing shock & 32768 & 16384 & 16384 & 0.0026 & 0.0038\\
      Moving shock & 32768 & 16384 & 16384 & 0.008 & 0.016\\
      Rarefaction wave & 32768 & 16384 & 16384 & 0.026 & 0.035\\
      Sine wave & 32768 & 16384 & 16384 & 0.051 & 0.064\\
      \hline
    \end{tabular}
\end{table}

\subsubsection{Standing Shock and Moving Shock}

The first two numerical experiments considered are the standing shock and moving shock, where the initial condition $\widetilde{u}(\lambda,0)$ is given by

\begin{equation*}
    \widetilde{u}(\lambda,0) =\begin{cases}
        1,&\lambda <0,\\
        -1,&\lambda \geq 0, \\
    \end{cases} \quad
    \widetilde{u}(\lambda,0) =\begin{cases}
        1,&\lambda <0,\\
        0,&\lambda \geq 0.\\
    \end{cases}
\end{equation*}
These initial conditions correspond to a shock wave stationary at $\lambda = 0$ and a shock wave moving at a speed of $0.5$, respectively:

\begin{equation*}
    \widetilde{u}(\lambda,t) =\begin{cases}
        1,&\lambda <0,\\
        -1,&\lambda \geq 0, \\
    \end{cases} \quad
    \widetilde{u}(\lambda,t) =\begin{cases}
        1,&\lambda <\frac{t}{2},\\
        0,&\lambda \geq \frac{t}{2}. \\
    \end{cases}
\end{equation*}
\begin{figure}[!ht]
\centering
\subfigure[Kruzkov entropy]{\includegraphics[height=4cm,width=6cm]{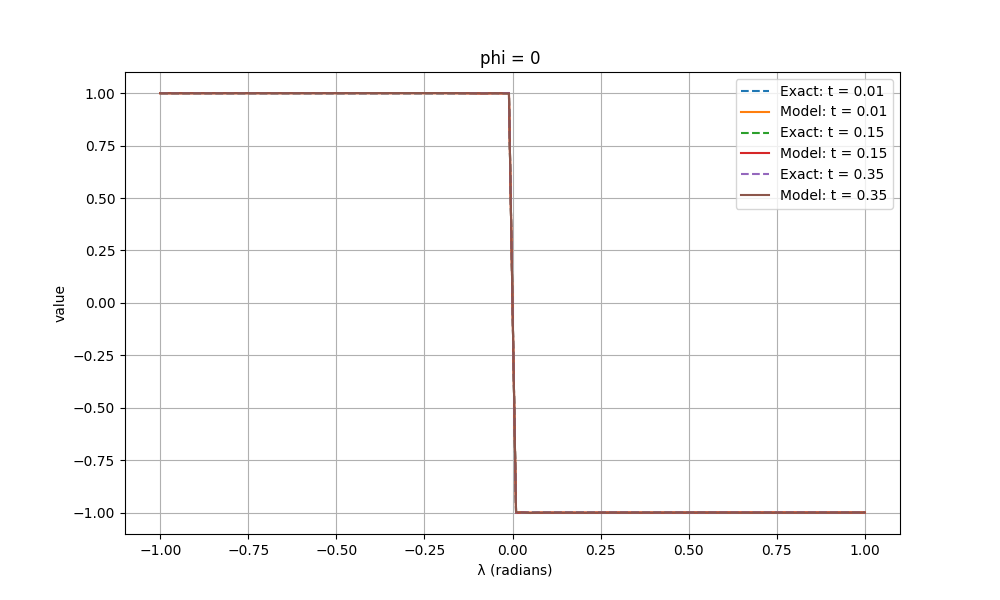}}
\subfigure[Quadratic entropy]{\includegraphics[height=4cm,width=6cm]{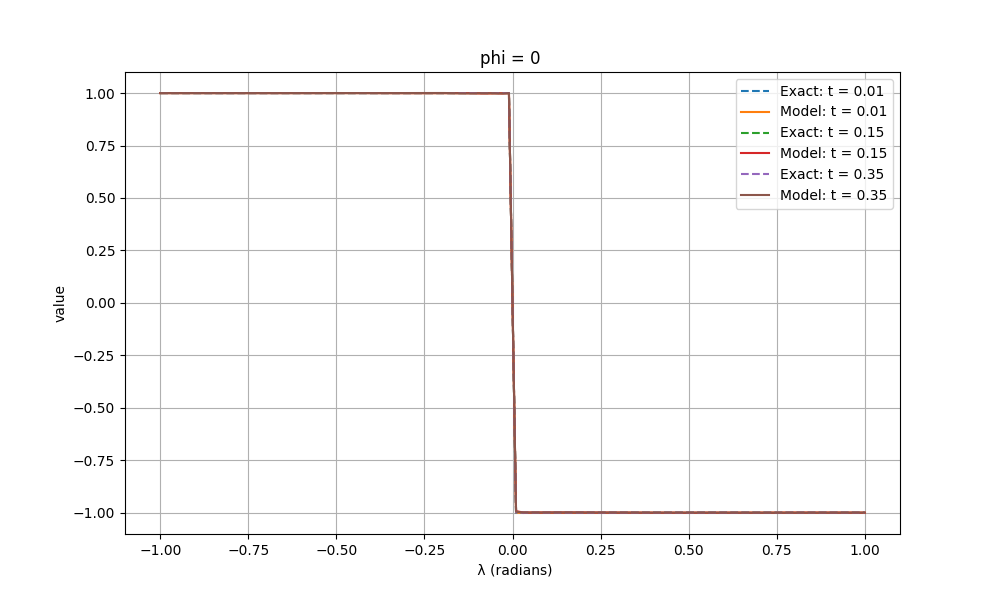}}
\captionsetup{justification=centering,singlelinecheck=false}
\caption{Standing shock}
\label{fig:Standing shock}
\end{figure}
For both experiments, the following hyperparameters were fixed: $N_{int} = 32768, N_{sb} = 16384, N_{ini} = 16384, N_{\min} = 1, N_{\max} = 6, N_{ep} = 5000, \tau_{\max} = 0.015$, and $\tau_{\min} = 0.01$. The solution network consists of 20 neurons per layer, and the test function network has 10 neurons per layer. A hyperparameter grid search is conducted for activation functions (ReLU, tanh, sin), solution network depths ($L_{\theta} = 4, 6$), test function network depths ($L_{\eta} = 2, 4$), reset frequencies ($r = 100, 200, 500, 1000$), and penalty parameters ($\varrho = 1, 5, 10$). The number of retrains is set to 10. Numerical results indicate that the choice of activation function (ReLU, tanh, sin) has only a minor impact, with the space-time test errors for all three remaining within the same order of magnitude. Since the theoretical analysis presented in this paper is specifically developed for the ReLU case, the figures below display the numerical results obtained using the ReLU activation function.

\begin{figure}[!ht]
\centering
\subfigure[Kruzkov entropy]{\includegraphics[height=4cm,width=6cm]{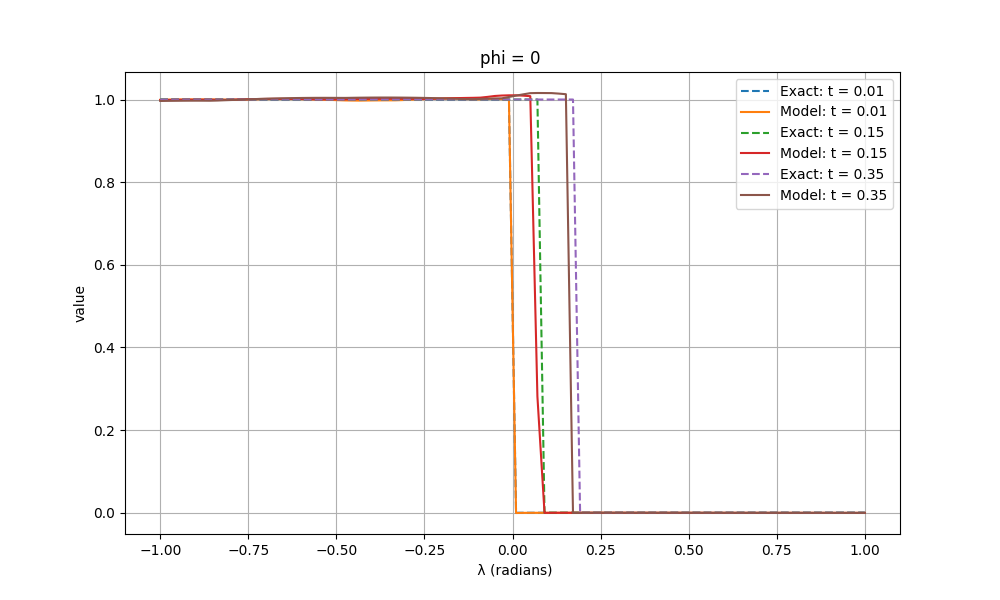}}
\subfigure[Quadratic entropy]{\includegraphics[height=4cm,width=6cm]{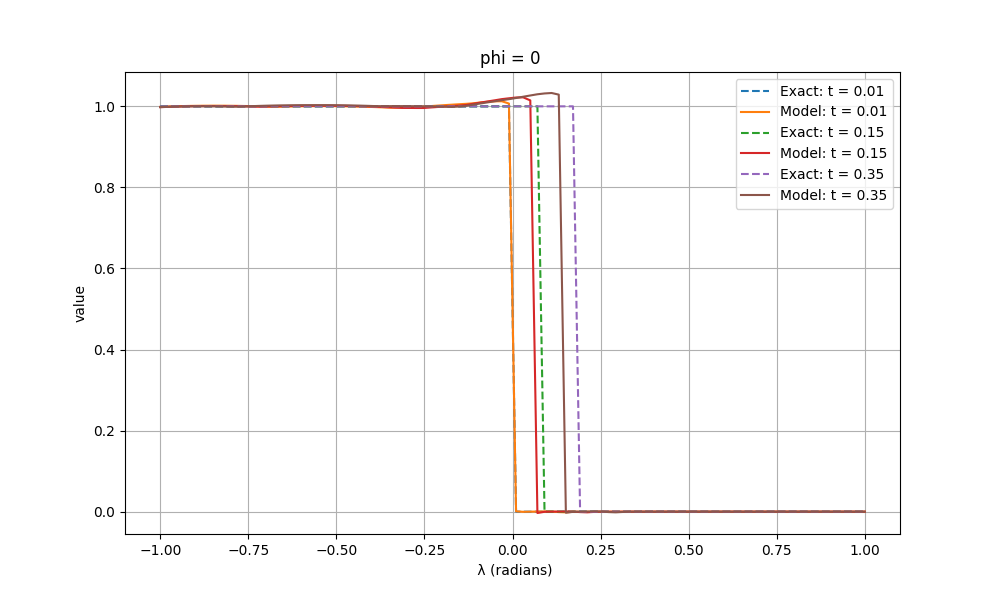}}
\caption{Moving shock}
\label{fig:Moving shock}
\end{figure}

~\autoref{fig:Standing shock} and~\autoref{fig:Moving shock} show the predicted solutions at $\phi = 0$ for the two experiments. For both types of parametric entropy, the model efficiently predicts the location of the shock wave. Notably, in the case of the moving shock, the solution predicted by the quadratic entropy exhibits slight oscillations. This issue could potentially be addressed by incorporating a total variation penalty, as suggested by~\cite{deryck2022wPINN}.

\begin{figure}[H]
    \centering
    \includegraphics[width=0.75\linewidth]{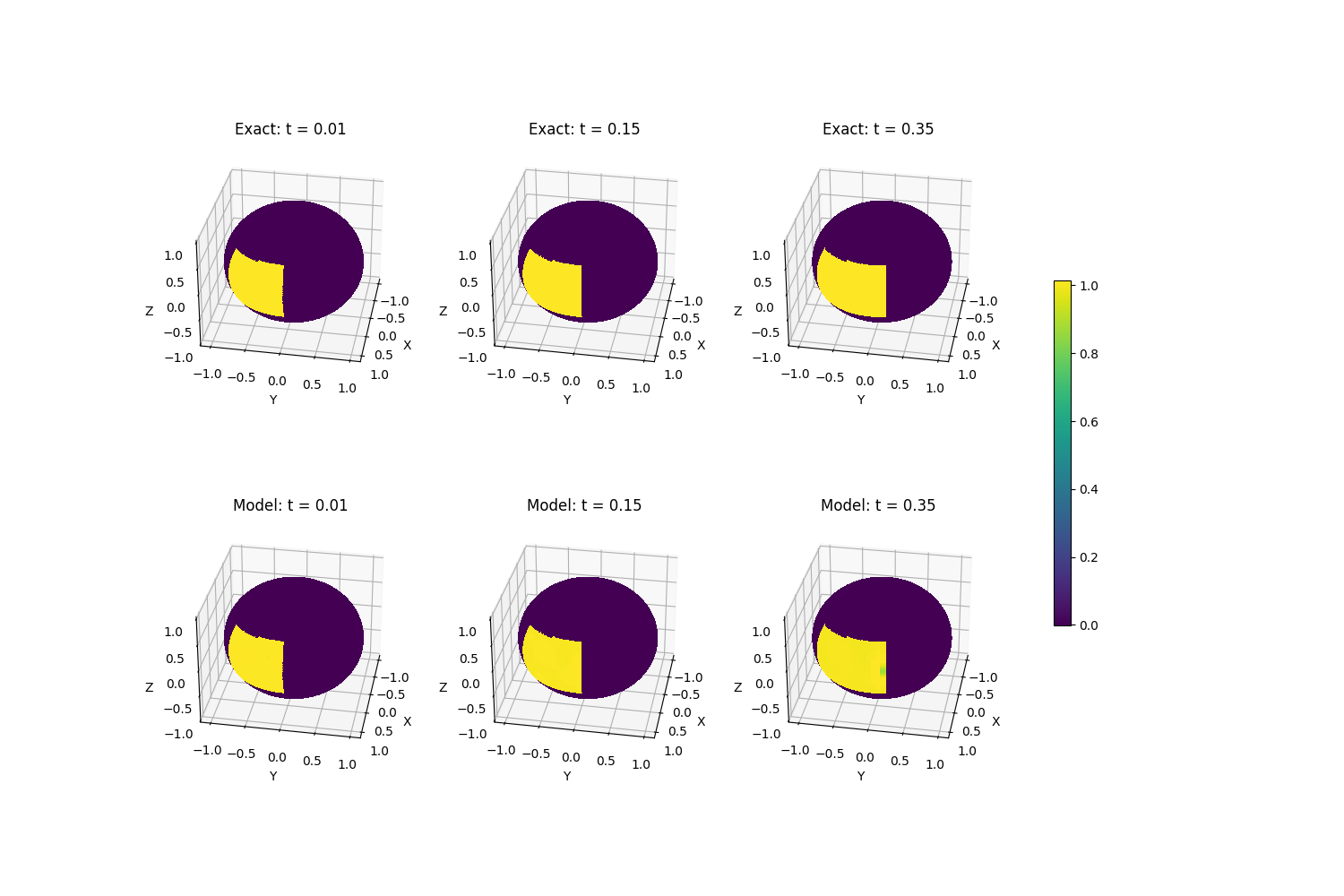}
    \caption{Moving shock solution predicted on the sphere.}
    \label{fig:spheres image of Moving shock}
\end{figure}

Although the space-time test errors for both types of parametric entropy are of the same order of magnitude, as shown in~\autoref{table:loss compare}, the space-time test errors for the Kruzkov entropy loss ($\mathcal{E}_{\mathrm{st}} = 0.26\%,\ 0.8\%$) are smaller than those for the quadratic entropy loss ($\mathcal{E}_{\mathrm{st}} = 0.38\%,\ 1.6\%$), which is also evident from~\autoref{fig:Moving shock}. 
\begin{figure}[H]
    \centering
    \includegraphics[width=0.75\linewidth]{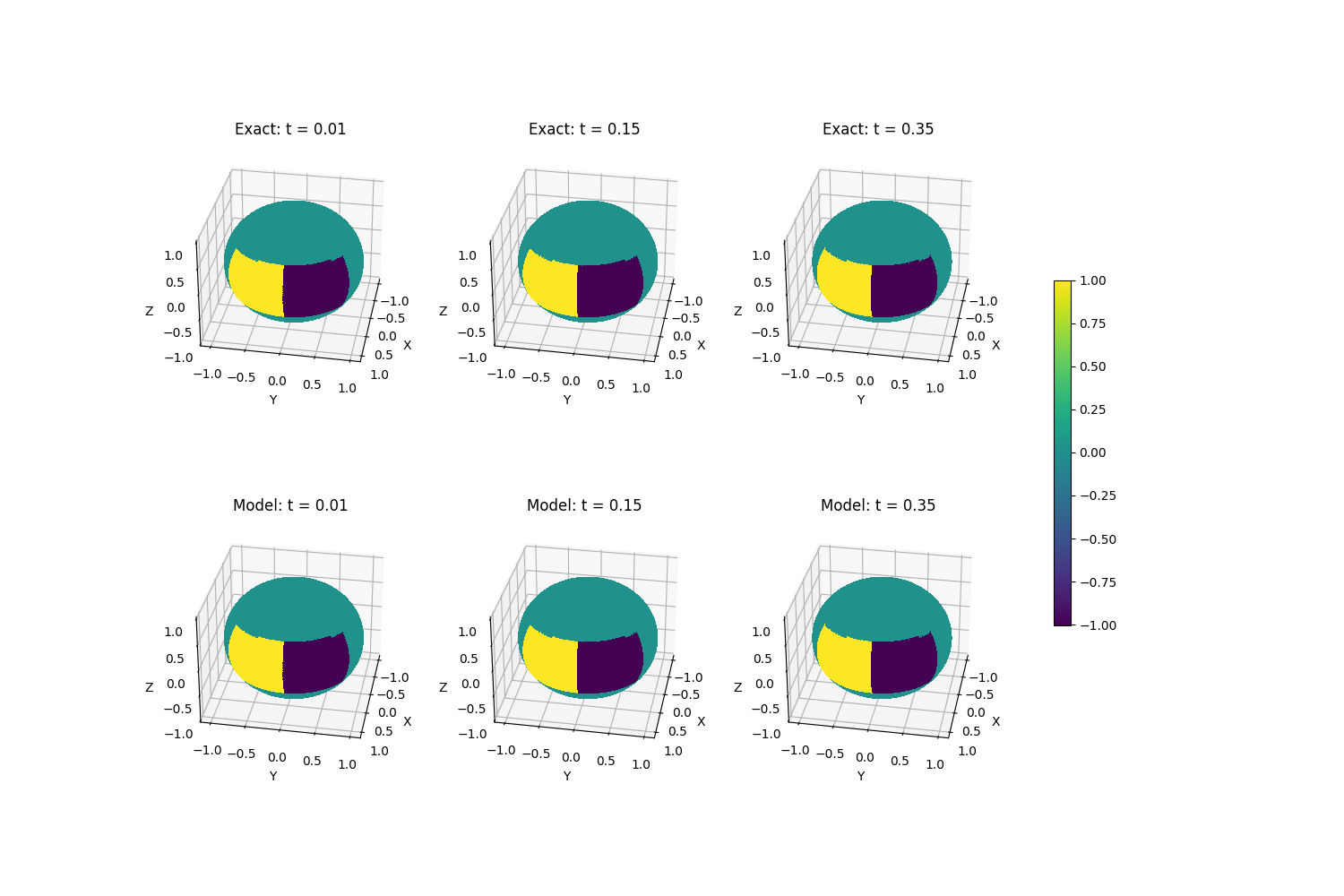}
    \caption{Standing shock solution predicted on the sphere.}
    \label{fig:spheres image of Standing shock}
\end{figure}
This difference might be attributed to the smoother nature of the quadratic entropy loss, which could pose challenges when approximating functions with sharp local variations. Subsequently, the predicted solutions on the sphere obtained using the Kruzkov entropy loss are shown in~\autoref{fig:spheres image of Standing shock} and~\autoref{fig:spheres image of Moving shock}.

The figures clearly illustrate the formation location and propagation speed of the shock. Essentially, these two numerical experiments require the neural network to learn a linear transformation, which is simpler than the transformations required in the other two experiments.  As a result, the space-time test errors are relatively smaller. Furthermore, subsequent experiments reveal only minor differences between the solution images generated using the Kruzkov entropy loss and those generated using the quadratic entropy loss. Given this similarity, the figures presented below are based on the solutions predicted with the Kruzkov entropy loss.

\subsubsection{Rarefaction Wave and Sine Wave}
We first consider the rarefaction wave generated by the following initial condition:
\begin{equation*}
    \widetilde{u}(\lambda,0) =\begin{cases}
        -1,&\lambda <0,\\
        1,&\lambda \geq 0.\\
    \end{cases}
\end{equation*}
This initial condition leads to a region $-t < \lambda < t$ where no characteristic lines pass through, resulting in the formation of a rarefaction wave. A rarefaction wave is a region where the velocity transitions smoothly from low to high:
\begin{equation*}
    \widetilde{u}(\lambda,t) =\begin{cases}
        -1,&\lambda <-t,\\
        \frac{\lambda}{t},&-t<\lambda \leq t, \\
        1,& \lambda \geq t. \\
    \end{cases}
\end{equation*}

\begin{figure}[H]
    \centering
    \includegraphics[width=0.75\linewidth]{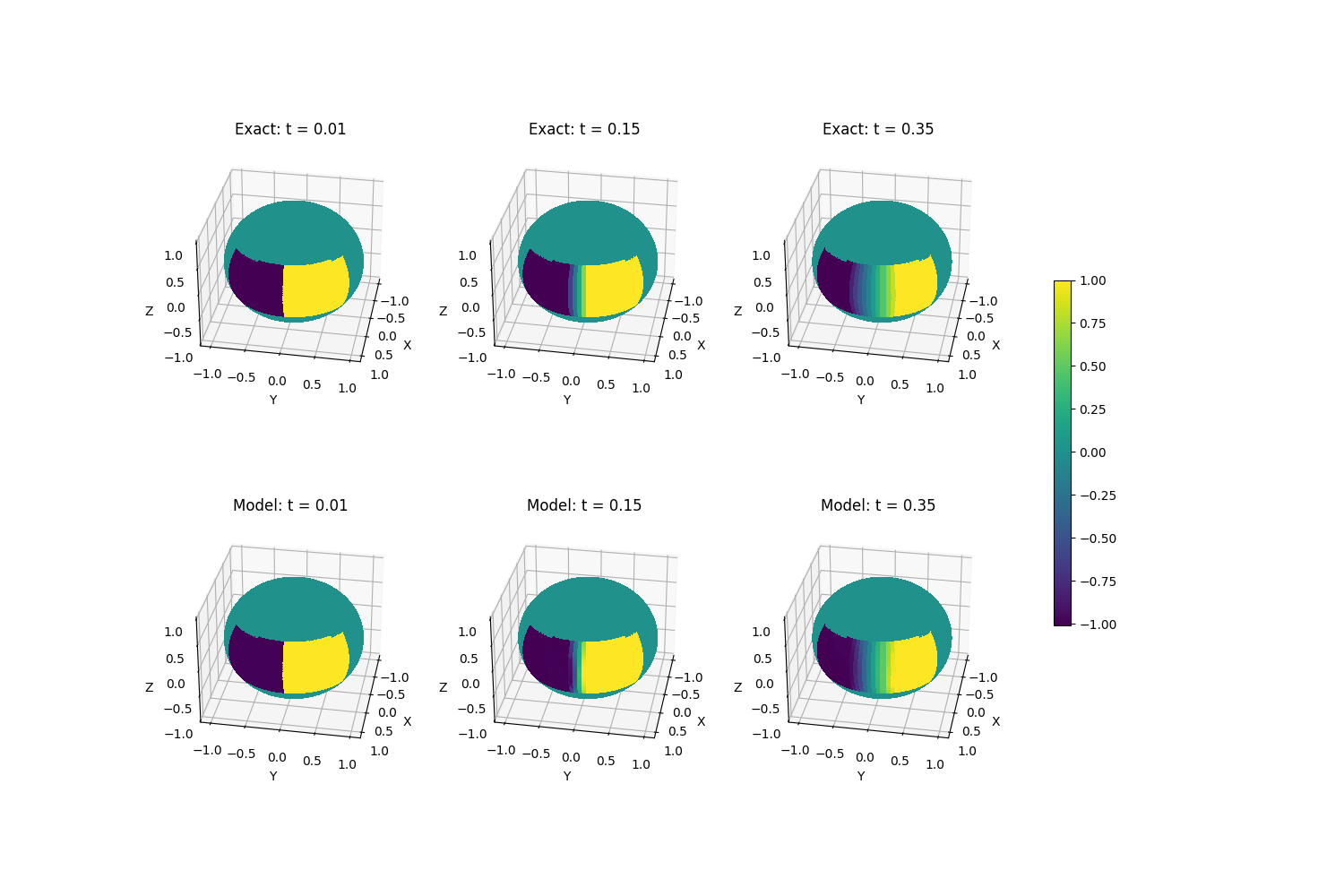}
    \caption{Rarefaction wave solution predicted on the sphere.}
    \label{fig:spheres image of Rarefaction wave}
\end{figure}

It is worth noting that the weak solution corresponding to this initial condition is not unique; the standing shock is also a weak solution. However, the rarefaction wave represents the unique entropy solution.~\autoref{table:loss compare} shows the space-time test errors for both types of parametric entropy. The space-time test error of the model trained with the Kruzkov entropy loss ($\mathcal{E}_{\mathrm{st}} = 2.6\%$) is slightly lower than that of the model trained with the quadratic entropy loss ($\mathcal{E}_{\mathrm{st}} = 3.5\%$). 
\begin{figure}[htbp]
\centering
\subfigure[Rarefaction wave]{\includegraphics[height=4cm,width=6cm]{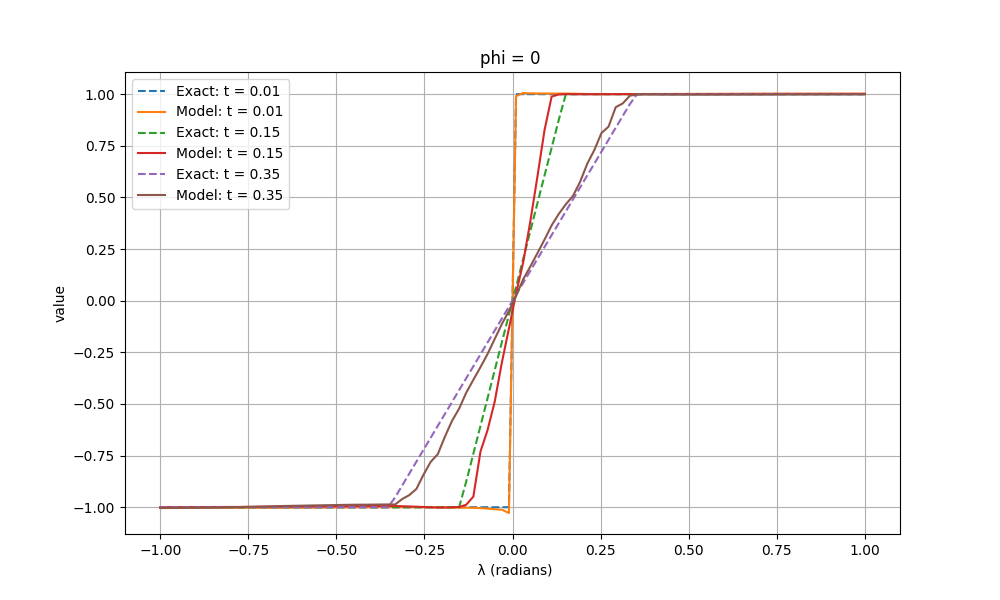}}
\subfigure[Sine wave]{\includegraphics[height=4cm,width=6cm]{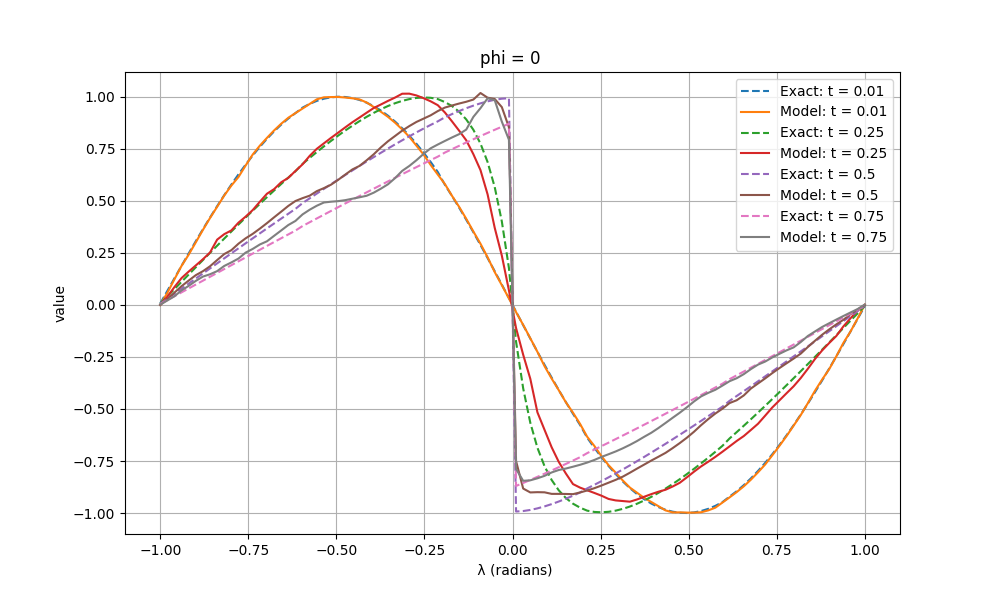}}
\caption{ Rarefaction wave and sine wave }
\label{fig:Rarefaction wave and Sine wave}
\end{figure}
Due to the more complex structure of the characteristic lines associated with this initial condition, the space-time test errors are larger than those in the standing shock and moving shock scenarios.~\autoref{fig:Rarefaction wave and Sine wave}(a) illustrates the predicted solution at $\phi = 0$.~\autoref{fig:spheres image of Rarefaction wave} shows the predicted solution on the sphere.

Finally, the sine wave case is presented as the last numerical experiment. Consider the following initial condition:
\begin{equation*}
    \widetilde{u}(\lambda,0) = -\sin(\pi \lambda),
\end{equation*}
with periodic boundary conditions.  Since the characteristic lines have varying slopes at different initial positions, they intersect at later times. As a result, the solution to the Burgers' equation with a sine wave initial condition deforms over time and eventually develops a shock wave. This scenario is more challenging than the previous experiments, as it requires the neural network to learn discontinuous behavior from smooth initial data. To improve training efficiency, a refined Sobol low-discrepancy sequence sampling method is employed~\cite{mishra2021enhancing}.

\begin{figure}[!ht]
    \centering
    \includegraphics[width=0.75\linewidth]{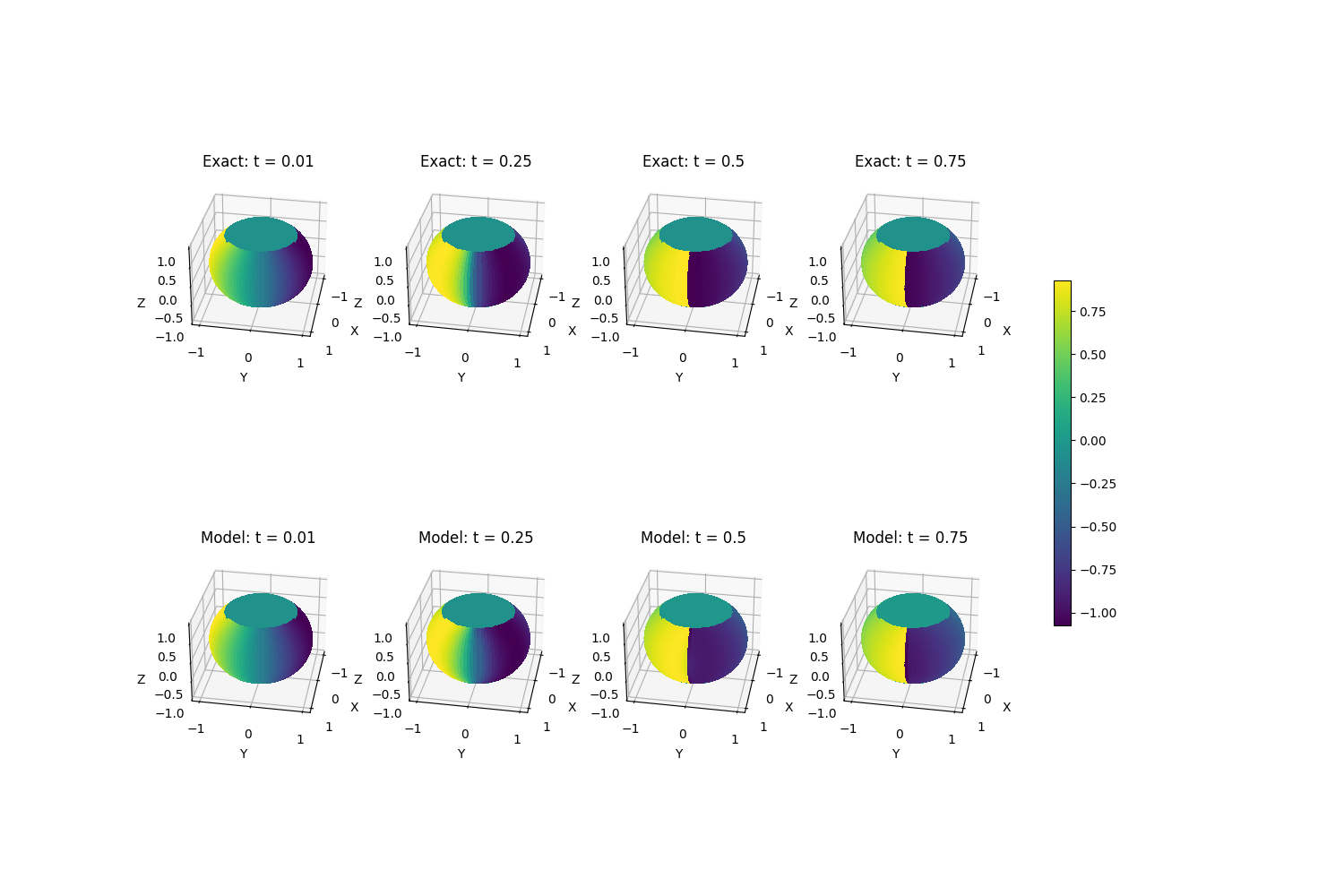}
    \caption{Sine wave solution predicted on the sphere.}
    \label{fig:spheres image of Sine wave}
\end{figure}
Adaptive sampling strategies could also be considered~\cite{wu2023comprehensive}. For this experiment, the following hyperparameters are fixed: $N_{int} = 32768, N_{sb} = 16384, N_{ini} = 16384, N_{\min} = 1, N_{ep} = 80000, \tau_{\max} = 0.015, \tau_{\min} = 0.01$. A hyperparameter grid search is conducted over solution network depths ($L_{\theta} = 6, 8$), test function network depths ($L_{\eta} = 4, 6$), reset frequencies ($r = 500, 1000, 2000, 5000, 10000$), $N_{\max}=6, 8$, and penalty parameters ($\varrho = 1, 5, 10$). The number of retrains is set to 15.~\autoref{table:loss compare} reports the space-time test errors for the average predicted solutions corresponding to the two types of entropy. The space-time test errors for the sine wave case ($\mathcal{E}_{\mathrm{st}}$-K $= 5.1\%$, $\mathcal{E}_{\mathrm{st}}$-S $= 6.4\%$) are larger than those in the previous three experiments, but remain within the same order of magnitude as the results reported in~\cite{deryck2022wPINN}. During experimentation, it is observed that increasing the number of collocation points could further reduce the space-time test errors; however, the optimal number of sampling points is not explored in this study.~\autoref{fig:Rarefaction wave and Sine wave}(b) displays the predicted solution at $\phi = 0$. Notable oscillations are present near the shock location, which constitutes the primary source of prediction error.~\autoref{fig:spheres image of Sine wave} shows the predicted solution on the sphere. The network effectively captures the shock location and its temporal evolution. It is worth noting that for problems sensitive to time dynamics, a pre-training strategy—dividing the temporal domain into shorter intervals for separate training—can produce solutions with sharper temporal resolution~\cite{guo2023pre}. However, this approach is not adopted here to maintain comparability with the original wPINN framework.

\section*{Appendix}
\addcontentsline{toc}{section}{Appendix}
\label{allapp}
The appendix is organized as follows. We first collect the notation used throughout this paper in~\autoref{Appendix: Notations}. We then present the detailed proof of the error decomposition in~\autoref{Appendix: Error Decomposition}. The proof of the localized \(L_1\)-stability estimate for the total error is given in~\autoref{Appendix: L1 Contraction for Total Error}, where we also define the ReQU localized test network family and prove the auxiliary kernel estimates needed for the stability argument. Finally,~\autoref{Appendix: Analysis of Approximation Error} and~\autoref{Appendix: Analysis of Quadrature Error} provide rigorous analyses of the approximation and quadrature errors, respectively.
\appendix
\section{Notations}\label{Appendix: Notations}
\smallskip
\noindent\textbf{\textit{Manifolds.}}\quad In this paper, we consider a smooth, closed, oriented $d$-dimensional Riemannian manifold $(\mathcal{M}^{d},g)$, where the metric tensor $ g $ is a positive definite, second-order covariant tensor field. For any point $ x \in \mathcal{M}^{d} $, $ T_x\mathcal{M}^{d} $ the tangent space at $x$, and the tangent bundle $ T\mathcal{M}^{d} = \bigcup_{x \in \mathcal{M}^{d}} T_x\mathcal{M}^{d} $, together with the cotangent bundle $ T^*\mathcal{M}^{d}$, forms the fundamental geometric structures on the manifold. In a local coordinate system $(x^i)$, the natural basis for the tangent space is given by $\partial_j = \frac{\partial}{\partial x^j}$, and the corresponding dual basis for the cotangent space is $dx^i$. Using the Einstein summation convention, the metric tensor is written as $g = g_{ij} \, dx^i \, dx^j$, where $g_{ij} = g(\partial_i, \partial_j)$. The inner product of any tangent vectors $X_x, Y_x \in T_x\mathcal{M}^{d}$ is defined as $\langle X_x, Y_x \rangle_{g_x} = g_x(X_x, Y_x)$, and the inverse of the metric matrix $(g_{ij})$ is denoted by $(g^{ij})$. 

For a smooth function $ f: \mathcal{M}^{d} \rightarrow \mathbb{R} $, its differential satisfies $ df_x(X_x) = X_x(f) $ for any tangent vector  $X_x$. The volume element on the manifold is denoted by $ dV_g $, and in local coordinates, it takes the form $ dV_g = \sqrt{|g|} \, dx^1 \wedge \dots \wedge dx^d $, where $ |g| = \det(g_{ij}) > 0 $.  The gradient of a function $ h $, denoted by $ \mathrm{grad}_g h $, is defined via the isomorphism induced by the Riemannian metric $ T\mathcal{M}^{d} \cong T^*\mathcal{M}^{d} $, with local expression $ \mathrm{grad}_g h = g^{ij} \frac{\partial h}{\partial x^i} \, \partial_j $. The Levi-Civita connection $ \nabla $, which is compatible with the metric, satisfies $ \nabla g = 0 $. Specifically, for a vector field $ X = X^j \partial_j $, its covariant derivative is given by $ \nabla_j X^i = \partial_j X^i + \Gamma^i_{kj} X^k $, where the Christoffel symbols are defined as $ \Gamma^i_{kj} = \frac{1}{2} g^{il} \left( \partial_k g_{lj} + \partial_j g_{lk} - \partial_l g_{kj} \right) $. The divergence of a vector field $ X $ is defined as $ \mathbf{div}_g X = \nabla_i X^i = \partial_i X^i + \Gamma^i_{ki} X^k $. 

On a compact manifold $ \mathcal{M}^d $, a partition of unity $ \{\rho_i\}_{i=1}^k $ is a collection of non-negative smooth functions satisfying $ \sum_{i=1}^k \rho_i = 1 $, with each $ \rho_i $ supported in a single element of an open cover. The associated coordinate maps $ \{\psi_i\}_{i=1}^k :  \mathcal{M}^d\to U_i \subset \mathbb{R}^d$ are smooth homeomorphisms onto their images, providing local coordinate systems on the manifold.

\smallskip
\noindent\textbf{\textit{Norms.}}\quad Let $U$ be a bounded domain. We introduce the function spaces $ L_p(U) $ and the Sobolev spaces $ W_p^k(U) $. For $ 1 \leq p \leq \infty $, the $ L_p $ norm is defined as
\begin{equation*}
    \|f\|_{L_p} = \begin{cases}
    \left( \int_{U} |f(x)|^p \, dx \right)^{1/p}, & 1 \leq p < \infty, \\
    \operatorname*{ess\,sup}_{x \in U} |f(x)|, & p = \infty.
    \end{cases} 
\end{equation*}
The Sobolev norm is defined as
\begin{equation*}
    \|f\|_{W_p^k} = \begin{cases}
    \left( \sum_{|\alpha| \leq k} \|D^\alpha f\|_{L_p}^p \right)^{1/p},  &1 \leq p < \infty, \\
    \max_{|\alpha| \leq k} \|D^\alpha f\|_{L_\infty}, & p = \infty,
    \end{cases}
\end{equation*}
where $ D^\alpha f $ denotes the $ \alpha $-th weak derivative of $ f $. The spaces $ L_p(U) $ and $ W_p^k(U) $ consist of functions with finite norms. Additionally, we define the seminorm
\begin{equation*}
    |f|_{W_p^k} = 
    \begin{cases}
        \left( \sum_{|\alpha| = k} \|D^\alpha f\|_{L_p}^p \right)^{1/p} ,  &1 \leq p < \infty,\\
         \max_{|\alpha| = k} \|D^\alpha f\|_{L_\infty}, & p = \infty.
    \end{cases}
\end{equation*}
 In particular, we consider the Sobolev space $W_0^{1,\infty}(U)$, which consists of functions in $W_\infty^1(U)$ whose trace vanishes on the boundary $\partial U$.

Throughout this appendix, the notation
\[
C(\texttt{parameter}_1,\ldots,\texttt{parameter}_m)
\]
denotes a positive constant depending only on the displayed parameters. Its value may change from line to line, and the same symbol may denote different constants in different proofs.

\section{Error Decomposition}\label{Appendix: Error Decomposition}

This appendix makes the error decomposition precise and records the scalar and local Rademacher deviation estimates used in the quadrature analysis. The empirical objective contains the positive part of the entropy residual, but the stochastic localization is performed on centered signed residual increments. This distinction is essential: the increment is signed, so one cannot use a bound of the form \(Ph^2\lesssim Ph\). The correct Bernstein relation is instead obtained from the stability-controlled quantity \(P|v-u^*|\).

Let
\[
   \mathcal Z:=\mathcal M^d\times[0,T],
   \qquad
   A_T:=T\operatorname{Vol}_g(\mathcal M^d),
   \qquad
   A_0:=\operatorname{Vol}_g(\mathcal M^d).
\]
Let \(P\) be the uniform probability measure on \(\mathcal Z\), let \(P_n\) be the interior empirical measure, and let \(P_0,P_{0,n}\) denote the population and empirical measures on \(\mathcal M^d\). Thus
\[
\mathcal R^{int}(v,\xi,c)=A_TP r^{int}(v,\xi,c),
\qquad
\mathcal R_n^{int}(v,\xi,c)=A_TP_n r^{int}(v,\xi,c),
\]
and
\[
\mathcal R^{tb}(v)=A_0P_0b_v,
\qquad
\mathcal R_n^{tb}(v)=A_0P_{0,n}b_v,
\qquad
b_v(x):=|v(x,0)-u_0(x)|.
\]

\begin{lemma}\label{lemma: bounds for Entropy Residual}
Let \(p,q>1\) satisfy \(1/p+1/q=1\), or let \(p=\infty\) and \(q=1\). Let \(u^*\) denote the entropy solution, and let \(u,v\in\mathcal F(L,W,S,B,M)\). Suppose
\[
M\ge
\max\{\|u^*\|_{L_\infty(\mathcal Z)},\|u\|_{L_\infty(\mathcal Z)},\|v\|_{L_\infty(\mathcal Z)}\}.
\]
For \(\xi\in W_0^{1,\infty}(\mathcal Z)\), define
\[
C_\xi:=
\|\partial_t\xi\|_{L_\infty(\mathcal Z)}
+3d\mathrm{Lip}_f\|\nabla_g\xi\|_{L_\infty(\mathcal Z)}.
\]
Then, pointwise on \(\mathcal Z\),
\begin{equation}\label{equation: pointwise residual lipschitz}
|r^{int}(u,\xi,c)-r^{int}(v,\xi,c)|
\le C_\xi|u-v|.
\end{equation}
Consequently,
\begin{align*}
|\mathcal R^{int}(u,\xi,c)-\mathcal R^{int}(v,\xi,c)|
&\le C_\xi\|u-v\|_{L_1(\mathcal Z)},
\\
|\mathcal R^{int}(u,\xi,c)-\mathcal R^{int}(v,\xi,c)|
&\le A_T^{1/p}C_\xi\|u-v\|_{L_q(\mathcal Z)}.
\end{align*}
If \(\xi\in\overline\Phi_\lambda\), then \(C_\xi\le C_\lambda\).
\end{lemma}

\begin{proof}
The Kruzkov entropy satisfies \(\bigl||u-c|-|v-c|\bigr|\le|u-v|.\) Moreover, Lemma 3.1 in~\cite{deryck2022wPINN} gives, componentwise, \(|\widetilde F(u,c)-\widetilde F(v,c)| \le 3d\mathrm{Lip}_f|u-v|.\) Substitution into the definition of \(r^{int}\) proves \eqref{equation: pointwise residual lipschitz}. Integration gives the \(L_1\)-estimate, and H{\"o}lder's inequality on the finite-measure set \(\mathcal Z\) gives the factor \(A_T^{1/p}\) in the \(L_q\)-estimate.
\end{proof}

\begin{lemma}\label{lemma: positive part contraction}
Let \(\Theta\) be any index set and let \(\{A_\theta\}_{\theta\in\Theta}\), \(\{B_\theta\}_{\theta\in\Theta}\subset\mathbb R\). Then
\[
\left|
\sup_{\theta\in\Theta}\ppart{A_\theta}
-
\sup_{\theta\in\Theta}\ppart{B_\theta}
\right|
\le
\sup_{\theta\in\Theta}|A_\theta-B_\theta|.
\]
\end{lemma}

\begin{proof}
The map \(z\mapsto[z]_+\) is one-Lipschitz. Hence, for every \(\theta\),
\[
[A_\theta]_+
\le [B_\theta]_++|A_\theta-B_\theta|
\le \sup_{\vartheta}[B_\vartheta]_+
   +\sup_{\vartheta}|A_\vartheta-B_\vartheta|.
\]
Taking the supremum over \(\theta\), and then interchanging \(A\) and \(B\), proves the assertion.
\end{proof}

Let
\[
   \Theta_\lambda^c:=\Theta_\lambda\times\mathcal C,
   \qquad
   \theta=(\vartheta,c),
   \qquad
   \xi_\theta:=\xi_\vartheta.
\]
Define
\[
   r_{v,\theta}:=r^{int}(v,\xi_\vartheta,c),
   \qquad
   h_{v,\theta}:=r_{v,\theta}-r_{u^*,\theta}.
\]
Also set
\begin{align*}
\Delta_n^{tb}(v)
&:=|(P_0-P_{0,n})b_v|,
\\
\Delta_n^{int}(v)
&:=\sup_{\theta\in\Theta_\lambda^c}|(P-P_n)h_{v,\theta}|,
\\
\Delta_{*,n}
&:=\sup_{\theta\in\Theta_\lambda^c}|(P-P_n)r_{u^*,\theta}|.
\end{align*}

\begin{proposition}\label{proposition: basic oracle decomposition}
Let \(\mathcal F\) be a trial class with finite class-wise bias \(B_{\lambda,\mathcal F}:=b_{\lambda,T}(\mathcal F)\). Define \(\mathfrak Q_\lambda\) and \(\mathfrak Q_{\lambda,n}\) by \eqref{equation: section5 population Q}--\eqref{equation: section5 empirical Q}, and let
\[
   u_{\mathcal F}\in\arg\min_{v\in\mathcal F}\mathfrak Q_\lambda(v),
   \qquad
   u_n\in\arg\min_{v\in\mathcal F}\mathfrak Q_{\lambda,n}(v).
\]
Then
\begin{align}
\mathfrak Q_\lambda(u_n)
\le{}&
\mathfrak Q_\lambda(u_{\mathcal F})
+A_0\bigl(\Delta_n^{tb}(u_n)+\Delta_n^{tb}(u_{\mathcal F})\bigr)
\notag\\
&+A_T^2\bigl(\Delta_n^{int}(u_n)+\Delta_n^{int}(u_{\mathcal F})+2\Delta_{*,n}\bigr).
\label{equation: exact oracle decomposition}
\end{align}
\end{proposition}

\begin{proof}
Since the bias is constant on \(\mathcal F\), empirical minimality gives
\[
\mathfrak Q_{\lambda,n}(u_n)
\le \mathfrak Q_{\lambda,n}(u_{\mathcal F}).
\]
Therefore,
\begin{align*}
\mathfrak Q_\lambda(u_n)-\mathfrak Q_\lambda(u_{\mathcal F})
\le{}&
\bigl(\mathfrak Q_\lambda-\mathfrak Q_{\lambda,n}\bigr)(u_n)
+
\bigl(\mathfrak Q_{\lambda,n}-\mathfrak Q_\lambda\bigr)(u_{\mathcal F}).
\end{align*}
The initial-time terms are bounded by \(A_0\Delta_n^{tb}\). For the internal loss, \autoref{lemma: positive part contraction} gives
\begin{align*}
|\mathfrak L_\lambda^+(v)-\mathfrak L_{\lambda,n}^+(v)|
&\le
\sup_{\theta}
|\mathcal R^{int}(v,\theta)-\mathcal R_n^{int}(v,\theta)|
\\
&=
A_T\sup_{\theta}|(P-P_n)r_{v,\theta}|
\\
&\le
A_T\Delta_n^{int}(v)+A_T\Delta_{*,n}.
\end{align*}
The outer factor \(A_T\) in \(\mathfrak Q_\lambda\) produces \(A_T^2\), and \eqref{equation: exact oracle decomposition} follows.
\end{proof}

We next record the local Rademacher deviation estimate used below; see~\cite{bartlett2005local,mendelson2003few}.

\begin{lemma}
\label{lemma: local Rademacher deviation}
Let $\mathcal G$ be a class of measurable functions satisfying $\|g\|_\infty\le b$ and $Pg^2\le\sigma^2$ for every $g\in\mathcal G$. Then the standard localized Rademacher deviation estimate gives, for every $t>0$, with probability at least $1-2e^{-t}$,
\begin{equation}\label{equation: class local Rademacher deviation}
\sup_{g\in\mathcal G}|(P_n-P)g|
\lesssim
\mathrm{Rad}_n(\mathcal G)
+
\sqrt{\frac{\sigma^2t}{n}}
+
\frac{bt}{n},
\end{equation}
where
\[
\mathrm{Rad}_n(\mathcal G)
:=
\mathbb E_{Z,\varepsilon}
\sup_{g\in\mathcal G}
\left|\frac1n\sum_{i=1}^n\varepsilon_i g(Z_i)\right|.
\]
The same estimate holds with the empirical Rademacher complexity on the right-hand side, after changing the universal constant. These forms are standard consequences of symmetrization and localization; see~\cite{bartlett2005local,mendelson2003few}.
\end{lemma}

The localized variance condition follows directly from \autoref{lemma: bounds for Entropy Residual}. Indeed,
\begin{equation*}
|h_{v,\theta}|
\le C_\lambda|v-u^*|
\le2MC_\lambda,
\end{equation*}
and
\begin{equation*}
P h_{v,\theta}^2
\le C_\lambda^2P|v-u^*|^2
\le2MC_\lambda^2P|v-u^*|.
\end{equation*}
Thus, if a localized stability sublevel set satisfies \(P|v-u^*|\lesssim\tau\), then \(h_{v,\theta}\) has envelope \(2MC_\lambda\) and variance \(O(C_\lambda^2\tau)\). Applying \eqref{equation: class local Rademacher deviation} to the corresponding localized residual class gives
\[
\sup_{v,\theta}|(P_n-P)h_{v,\theta}|
\lesssim
\mathrm{Rad}_n(\mathcal H_\tau)
+C_\lambda\sqrt{\frac{\tau t}{n}}
+MC_\lambda\frac{t}{n}.
\]
This is the localized Rademacher step that converts the local \(L_1\)-stability estimate into a fast quadrature rate. The remaining task is to bound \(\mathrm{Rad}_n(\mathcal H_\tau)\) by Dudley's entropy integral and to solve the resulting critical-radius inequality; this is carried out in \autoref{Appendix: Analysis of Quadrature Error}.

\section{\texorpdfstring{$L_1$ Contraction for Total Error}{L1 Contraction for Total Error}}\label{Appendix: L1 Contraction for Total Error}
In this section, we prove~\autoref{theorem: L1 contraction for total error}. The proof is based on a localized doubling-of-variables argument. We first record the properties of the localized test functions used in the proof of~\autoref{theorem: L1 contraction for total error}. Note that the functions \(\rho_\eta\), \(\chi_h\), and \(K_\alpha^{\mathrm{ext}}\) are exactly representable by ReQU networks in the variables \((\iota(x),t)\), with \((\iota(y),s)\) treated as parameters. Indeed, ReQU networks exactly
represent the square function through $r^2=\sigma_2(r)+\sigma_2(-r),$
and hence exactly represent squared distances and products via $ab=\frac{(a+b)^2-(a-b)^2}{4}.$
Therefore, the localized test functions
\[
    \Phi_{\lambda,\tau}^{y,s}(x,t)
    =
    \chi_{h,\tau}\left(\frac{t+s}{2}\right)
    \rho_\eta(t-s)
    K_\alpha^{\mathrm{ext}}(x,y)
\]
form a fixed ReQU localized test neural network. In the proof below we first fix the terminal horizon $T$ and use the abbreviations
\[
\chi_h:=\chi_{h,T},
\qquad
\Phi_\lambda^{y,s}:=\Phi_{\lambda,T}^{y,s}.
\]
The horizon-indexed statement follows by replacing $T$ everywhere by $\tau$, $\chi_h$ by $\chi_{h,\tau}$, and $\Phi_\lambda$ by $\Phi_{\lambda,\tau}$. The trial function is denoted by $v$, as in the theorem statement.
\begin{lemma}
\label{lemma: extrinsic ReQU kernel estimates}
There exist constants \(C>0\) and \(\alpha_0>0\), depending only on
\((\mathcal M^d,g)\) and the embedding \(\iota\), such that for all
\(0<\alpha\le\alpha_0\), the following estimates hold.

First,
\begin{equation}\label{eq:requ-kernel-support}
    0\le K_\alpha^{\mathrm{ext}}(x,y)
    \le
    C\alpha^{-d}
    \mathbf 1_{\{d_g(x,y)\le C\alpha\}}.
\end{equation}
Second,
\begin{equation}\label{eq:requ-kernel-gradient}
    |\nabla_xK_\alpha^{\mathrm{ext}}(x,y)|_g
    +
    |\nabla_yK_\alpha^{\mathrm{ext}}(x,y)|_g
    \le
    C\alpha^{-d-1}
    \mathbf 1_{\{d_g(x,y)\le C\alpha\}}.
\end{equation}
Third, if $Z_\alpha(x):=\int_{\mathcal M^d}K_\alpha^{\mathrm{ext}}(x,y)\,dV_g(y),$ then
\begin{align}
    |Z_\alpha(x)-1|
    &\le C\alpha^2,
    \label{eq:requ-kernel-normalization}
    \\
    |\nabla_gZ_\alpha(x)|_g
    &\le C\alpha.
    \label{eq:requ-kernel-normalization-gradient}
\end{align}
Fourth, the cancellation estimate holds:
\begin{equation}\label{eq:requ-kernel-cancellation}
    \left|
    \nabla_xK_\alpha^{\mathrm{ext}}(x,y)
    +
    P_{y\to x}\nabla_yK_\alpha^{\mathrm{ext}}(x,y)
    \right|_{g_x}
    \le
    C\alpha^{-d}
    \mathbf 1_{\{d_g(x,y)\le C\alpha\}},
\end{equation}
where \(P_{y\to x}:T_y\mathcal M^d\to T_x\mathcal M^d\) denotes parallel
transport along the unique minimizing geodesic.

Finally, for every \(w\in BV(\mathcal M^d)\),
\begin{equation}\label{eq:requ-indicator-bv-localization}
    \int_{\mathcal M^d}\int_{\mathcal M^d}
    \alpha^{-d}
    \mathbf 1_{\{d_g(x,y)\le C\alpha\}}
    |w(x)-w(y)|
    \,dV_g(x)dV_g(y)
    \le
    C\alpha\,TV(w),
\end{equation}
and consequently
\begin{equation}\label{eq:requ-bv-localization}
    \int_{\mathcal M^d}\int_{\mathcal M^d}
    K_\alpha^{\mathrm{ext}}(x,y)|w(x)-w(y)|
    \,dV_g(x)dV_g(y)
    \le
    C\alpha\,TV(w).
\end{equation}
\end{lemma}

\begin{proof}
The support estimate follows from the definition. If
\(K_\alpha^{\mathrm{ext}}(x,y)>0\), then $|\iota(x)-\iota(y)|<\alpha.$
Since \(\mathcal M^d\) is compact and smoothly embedded, the extrinsic and
intrinsic distances are locally equivalent. Hence, after choosing
\(\alpha_0>0\) sufficiently small, we may require both
\[
    |\iota(x)-\iota(y)|<\alpha
    \quad\Longrightarrow\quad
    d_g(x,y)\le C\alpha
\]
and
\[
    C\alpha_0<\operatorname{inj}(\mathcal M^d).
\]
Thus every pair in the support of the kernel is joined by a unique minimizing geodesic, and the parallel transport used below is well defined.
Therefore,
\[
    0\le K_\alpha^{\mathrm{ext}}(x,y)
    \le
    C\alpha^{-d}
    \mathbf 1_{\{d_g(x,y)\le C\alpha\}},
\]
which proves~\eqref{eq:requ-kernel-support}.

For the gradient estimate, define the ambient kernel
\[
    \overline K_\alpha(X,Y)
    :=
    c_d\alpha^{-d}
    \left(
    1-\frac{|X-Y|^2}{\alpha^2}
    \right)_+^2.
\]
On its support, \(|X-Y|\le\alpha\). A direct differentiation gives
\[
    \nabla_X\overline K_\alpha(X,Y)
    =
    -4c_d\alpha^{-d-2}(X-Y)
    \left(
    1-\frac{|X-Y|^2}{\alpha^2}
    \right)_+,
\]
and similarly for \(\nabla_Y\overline K_\alpha\). Hence,
\[
    |\nabla_X\overline K_\alpha|
    +
    |\nabla_Y\overline K_\alpha|
    \le
    C\alpha^{-d-1}.
\]
Since the Riemannian gradients are the tangent projections of the ambient
gradients, we obtain~\eqref{eq:requ-kernel-gradient}.

We next prove the approximate normalization. Fix \(x\in\mathcal M^d\) and
write \(y=\exp_x z\) in normal coordinates. Since \(\iota\) is an isometric
embedding,
\[
    |\iota(\exp_x z)-\iota(x)|^2
    =
    |z|^2+O(|z|^4),
\]
and the Riemannian volume form satisfies
\[
    dV_g(y)=\bigl(1+O(|z|^2)\bigr)\,dz.
\]
Thus, using \(z=\alpha w\),
\[
\begin{aligned}
    Z_\alpha(x)
    &=
    c_d
    \int_{\mathbb R^d}
    \left(
    1-|w|^2+O(\alpha^2|w|^4)
    \right)_+^2
    \bigl(1+O(\alpha^2|w|^2)\bigr)\,dw  \\
    &=
    c_d
    \int_{\mathbb R^d}(1-|w|^2)_+^2\,dw
    +O(\alpha^2)
    =
    1+O(\alpha^2).
\end{aligned}
\]
This gives~\eqref{eq:requ-kernel-normalization}. By compactness, the normal-coordinate expansions above, together with their first derivatives with respect to the base point $x$, are uniform in $x$. Differentiating under the integral sign, the leading Euclidean term is independent of $x$, whereas the differentiated geometric correction terms are $O(\alpha)$ after the change of variables $z=\alpha w$. Consequently,
\[
    \sup_{x\in\mathcal M^d}|\nabla_gZ_\alpha(x)|_g\le C\alpha,
\]
which proves~\eqref{eq:requ-kernel-normalization-gradient}.

We now prove the cancellation estimate. Let \(P_x:\mathbb R^D\to
T_x\mathcal M^d\) be the orthogonal projection. Since
\[
    \nabla_Y\overline K_\alpha(X,Y)
    =
    -\nabla_X\overline K_\alpha(X,Y),
\]
we have
\[
    \nabla_xK_\alpha^{\mathrm{ext}}
    =
    P_x\nabla_X\overline K_\alpha,
    \qquad
    \nabla_yK_\alpha^{\mathrm{ext}}
    =
    P_y\nabla_Y\overline K_\alpha
    =
    -P_y\nabla_X\overline K_\alpha.
\]
Therefore,
\[
    \nabla_xK_\alpha^{\mathrm{ext}}
    +
    P_{y\to x}\nabla_yK_\alpha^{\mathrm{ext}}
    =
    (P_x-P_{y\to x}P_y)\nabla_X\overline K_\alpha.
\]
Because the tangent bundle is smooth and \(\mathcal M^d\) is compact,
\[
    \|P_x-P_{y\to x}P_y\|
    \le
    C d_g(x,y).
\]
On the support of \(K_\alpha^{\mathrm{ext}}\), \(d_g(x,y)\le C\alpha\).
Together with
\[
    |\nabla_X\overline K_\alpha|
    \le
    C\alpha^{-d-1},
\]
we obtain~\eqref{eq:requ-kernel-cancellation}.

It remains to prove the \(BV\)-localization estimate. We first prove
\eqref{eq:requ-indicator-bv-localization} for
\(w\in C^1(\mathcal M^d)\).
Since \(\mathcal M^d\) is compact, it can be covered by finitely many
normal coordinate charts. For \(0<\alpha\le\alpha_0\), the condition
\(d_g(x,y)\le C\alpha\) ensures that \(x\) and \(y\) lie in a common
coordinate chart after possibly refining the finite covering. In these
charts, the Riemannian distance and volume form are uniformly comparable
with the Euclidean distance and Lebesgue measure. Hence it suffices to
estimate local terms of the form
\[
    \int \int_{|h|\le C\alpha}
    \alpha^{-d}|w(z+h)-w(z)|\,dh\,dz.
\]
By the fundamental theorem of calculus along the line segment in the
coordinate chart, $|w(z+h)-w(z)|\le|h|\int_0^1|\nabla w(z+\theta h)|\,d\theta.$ Therefore,
\[
\begin{aligned}
&\int\!\!\int_{|h|\le C\alpha}
\alpha^{-d}|w(z+h)-w(z)|\,dh\,dz
\\
&\qquad\le
\int_0^1\!\int\!\!\int_{|h|\le C\alpha}
\alpha^{-d}|h|\,|\nabla w(z+\theta h)|
\,dh\,dz\,d\theta .
\end{aligned}
\]
Changing variables \(z+\theta h\mapsto z\), and using the uniform
boundedness of the Jacobians of the coordinate changes, we obtain
\[
\begin{aligned}
&\int\!\!\int_{|h|\le C\alpha}
\alpha^{-d}|w(z+h)-w(z)|\,dh\,dz
\\
&\qquad\le
C\left(
\int_{|h|\le C\alpha}\alpha^{-d}|h|\,dh
\right)
\int_{\mathcal M^d}|\nabla_g w|_g\,dV_g .
\end{aligned}
\]
With the change of variables \(h=\alpha q\), $\int_{|h|\le C\alpha}\alpha^{-d}|h|\,dh=\alpha\int_{|q|\le C}|q|\,dq\le C\alpha.$
Thus,
\[
    \int_{\mathcal M^d}\int_{\mathcal M^d}
    \alpha^{-d}
    \mathbf 1_{\{d_g(x,y)\le C\alpha\}}
    |w(x)-w(y)|
    \,dV_g(x)dV_g(y)
    \le
    C\alpha
    \int_{\mathcal M^d}|\nabla_g w|_g\,dV_g .
\]
Since \(w\in C^1(\mathcal M^d)\), we have $TV(w)=\int_{\mathcal M^d}|\nabla_g w|_g\,dV_g.$
Therefore,
\[
    \int_{\mathcal M^d}\int_{\mathcal M^d}
    \alpha^{-d}
    \mathbf 1_{\{d_g(x,y)\le C\alpha\}}
    |w(x)-w(y)|
    \,dV_g(x)dV_g(y)
    \le
    C\alpha\,TV(w).
\]

For a general \(w\in BV(\mathcal M^d)\), choose a strict
\(BV\)-approximation \(w_k\in C^\infty(\mathcal M^d)\) such that
\[
    w_k\to w
    \quad\text{in }L^1(\mathcal M^d),
    \qquad
    TV(w_k)\to TV(w).
\]
Applying the preceding estimate to \(w_k\), we get
\[
    \int_{\mathcal M^d}\int_{\mathcal M^d}
    \alpha^{-d}
    \mathbf 1_{\{d_g(x,y)\le C\alpha\}}
    |w_k(x)-w_k(y)|
    \,dV_g(x)dV_g(y)
    \le
    C\alpha\,TV(w_k).
\]
Since the kernel $(x,y)\mapsto \alpha^{-d}\mathbf 1_{\{d_g(x,y)\le C\alpha\}}$ is integrable on \(\mathcal M^d\times\mathcal M^d\) for fixed \(\alpha\),
we can pass to the limit \(k\to\infty\). Hence,
\[
    \int_{\mathcal M^d}\int_{\mathcal M^d}
    \alpha^{-d}
    \mathbf 1_{\{d_g(x,y)\le C\alpha\}}
    |w(x)-w(y)|
    \,dV_g(x)dV_g(y)
    \le
    C\alpha\,TV(w),
\]
which proves \eqref{eq:requ-indicator-bv-localization}.

Finally, by the support estimate \eqref{eq:requ-kernel-support}, it holds that
\[
    \int_{\mathcal M^d}\int_{\mathcal M^d}
    K_\alpha^{\mathrm{ext}}(x,y)|w(x)-w(y)|
    \,dV_g(x)dV_g(y)
    \le
    C\alpha\,TV(w),
\]
which proves \eqref{eq:requ-bv-localization}.
\end{proof}

\begin{proof}[Proof of~\autoref{theorem: L1 contraction for total error}]
For \((x,t,y,s)\in
\mathcal M^d\times[0,T]\times\mathcal M^d\times[0,T]\), set
\[
  \Phi(x,t,y,s)
  :=
  \chi_h\left(\frac{t+s}{2}\right)
  \rho_\eta(t-s)
  K_\alpha^{\mathrm{ext}}(x,y).
\]
For brevity, write \(K_\alpha:=K_\alpha^{\mathrm{ext}}\). Since
\(K_\alpha(x,y)=K_\alpha(y,x)\) and \(\rho_\eta\) is even, we have
\[
    \Phi(x,t,y,s)=\Phi(y,s,x,t).
\]
Moreover, \(\Phi\ge0\), and the condition \(0<\eta<h/4\) ensures that the
support of \(\Phi\) stays away from the temporal boundary.

The localized ReQU test functions belong to $W_0^{1,\infty}(\mathcal M^d\times(0,T))$ and are nonnegative. Their use in the entropy inequality is therefore justified by \autoref{lemma: admissible Sobolev entropy tests}.

For fixed \((x,t)\), apply the entropy inequality for \(u\) with the constant $c=v(x,t)$
and with the nonnegative test function \((y,s)\mapsto\Phi(x,t,y,s)\). For almost every fixed $(x,t)$, the value $v(x,t)$ is used only as a real constant in the entropy inequality and is not required to belong to $\mathcal C$. This gives
\[
\begin{aligned}
\int_0^T\int_{\mathcal M^d}
\Big[
&\widetilde U(u(y,s),v(x,t))\partial_s\Phi(x,t,y,s)
\\
&+
g_y\bigl(
\widetilde F_y(u(y,s),v(x,t)),
\nabla_y\Phi(x,t,y,s)
\bigr)
\Big]
\,dV_g(y)\,ds
\ge0.
\end{aligned}
\]
Integrating this inequality with respect to \((x,t)\), and using the symmetry
\[
  \widetilde U(u(y,s),v(x,t))
  =
  \widetilde U(v(x,t),u(y,s)),
\]
\[
  \widetilde F_y(u(y,s),v(x,t))
  =
  \widetilde F_y(v(x,t),u(y,s)),
\]
we obtain
\[
\begin{aligned}
\int_0^T\int_{\mathcal M^d}\int_0^T\int_{\mathcal M^d}
\Big[
&\widetilde U(v(x,t),u(y,s))\partial_s\Phi
\\
&+
g_y\bigl(
\widetilde F_y(v(x,t),u(y,s)),
\nabla_y\Phi
\bigr)
\Big]
\,dV_g(y)\,ds\,dV_g(x)\,dt
\ge0.
\end{aligned}
\]

On the other hand, for fixed \((y,s)\), by the definition of the signed
residual \(\mathcal R^{int}\),
\[
\begin{aligned}
&\mathcal R^{int}(v,\Phi_\lambda^{y,s},u(y,s))
\\
&\quad=
-\int_0^T\int_{\mathcal M^d}
\Big[
\widetilde U(v(x,t),u(y,s))\partial_t\Phi(x,t,y,s)
\\
&\hspace{41mm}+
 g_x\bigl(
\widetilde F_x(v(x,t),u(y,s)),
\nabla_x\Phi(x,t,y,s)
\bigr)
\Big]
\,dV_g(x)\,dt.
\end{aligned}
\]
Since \(u(y,s)\in\mathcal C\) for a.e. \((y,s)\), we have
\[
  \mathcal R^{int}(v,\Phi_\lambda^{y,s},u(y,s))
  \le
  \mathcal L_{\lambda,T}^+(v).
\]
Integrating in \((y,s)\), we get
\[
\begin{aligned}
-\int_0^T\int_{\mathcal M^d}\int_0^T\int_{\mathcal M^d}
\Big[
&\widetilde U(v(x,t),u(y,s))\partial_t\Phi
\\
&+
g_x\bigl(
\widetilde F_x(v(x,t),u(y,s)),
\nabla_x\Phi
\bigr)
\Big]
\,dV_g(x)\,dt\,dV_g(y)\,ds
\\
\le
T\operatorname{Vol}_g(\mathcal M^d)\mathcal L_{\lambda,T}^+(v).
\end{aligned}
\]
Equivalently,
\[
\begin{aligned}
\int_0^T\int_{\mathcal M^d}\int_0^T\int_{\mathcal M^d}
\Big[
&\widetilde U(v(x,t),u(y,s))\partial_t\Phi
\\
&+
g_x\bigl(
\widetilde F_x(v(x,t),u(y,s)),
\nabla_x\Phi
\bigr)
\Big]
\,dV_g(x)\,dt\,dV_g(y)\,ds
\\
\ge
-
T\operatorname{Vol}_g(\mathcal M^d)\mathcal L_{\lambda,T}^+(v).
\end{aligned}
\]
Adding the two estimates gives
\begin{equation}\label{eq:finite-scale-A-B}
  A_\lambda(v)+B_\lambda(v)
  \ge
  -
  T\operatorname{Vol}_g(\mathcal M^d)\mathcal L_{\lambda,T}^+(v),
\end{equation}
where
\[
\begin{aligned}
A_\lambda(v)
:=
\int_0^T\int_{\mathcal M^d}\int_0^T\int_{\mathcal M^d}
&\widetilde U(v(x,t),u(y,s))
(\partial_t+\partial_s)\Phi(x,t,y,s)
\\
&\qquad\qquad
\,dV_g(x)\,dt\,dV_g(y)\,ds,
\end{aligned}
\]
and
\[
\begin{aligned}
B_\lambda(v)
:=\int_0^T\int_{\mathcal M^d}\int_0^T\int_{\mathcal M^d}
\Big[
&g_x\bigl(
\widetilde F_x(v(x,t),u(y,s)),
\nabla_x\Phi(x,t,y,s)
\bigr)
\\
&+
g_y\bigl(
\widetilde F_y(v(x,t),u(y,s)),
\nabla_y\Phi(x,t,y,s)
\bigr)
\Big]
\\[-1mm]
&\hspace{43mm}
\,dV_g(x)\,dt\,dV_g(y)\,ds.
\end{aligned}
\]

Define
\[
    \mathcal E(t)
    :=
    \|v(\cdot,t)-u(\cdot,t)\|_{L_1(\mathcal M^d)}.
\]
Set
\[
  b_{\lambda,T}(v)
  :=
  \left|
  A_\lambda(v)-\bigl(\mathcal E(0)-\mathcal E(T)\bigr)
  \right|
  +
  |B_\lambda(v)|.
\]
Then~\eqref{eq:finite-scale-A-B} implies
\[
  \mathcal E(T)
  \le
  \mathcal E(0)
  +
  T\operatorname{Vol}_g(\mathcal M^d)\mathcal L_{\lambda,T}^+(v)
  +
  b_{\lambda,T}(v).
\]
Since \(u(\cdot,0)=u_0\), this gives
\[
\|v(\cdot,T)-u(\cdot,T)\|_{L_1(\mathcal M^d)}
\le
\|v(\cdot,0)-u_0\|_{L_1(\mathcal M^d)}
+
T\operatorname{Vol}_g(\mathcal M^d)\mathcal L_{\lambda,T}^+(v)
+
b_{\lambda,T}(v).
\]

It remains to estimate \(b_{\lambda,T}(v)\). Since $\Phi(x,t,y,s)=\chi_h\left(\frac{t+s}{2}\right)\rho_\eta(t-s)K_\alpha(x,y),$
we have
\[
  (\partial_t+\partial_s)\Phi(x,t,y,s)
  =
  \chi_h'\left(\frac{t+s}{2}\right)
  \rho_\eta(t-s)
  K_\alpha(x,y).
\]
Introduce $\tau:=\frac{t+s}{2},r:=t-s.$
Then \(t=\tau+r/2\), \(s=\tau-r/2\), and \(dt\,ds=d\tau\,dr\). Because
\(\chi_h(\tau)\rho_\eta(r)\neq0\) implies $\tau\in[h,T-h],|r|\le\eta<h/4,$ the points \(\tau\pm r/2\) remain in \([0,T]\). Therefore,
\[
  A_\lambda(v)
  =
  \int_0^T\chi_h'(\tau)J_{\alpha,\eta}(\tau)\,d\tau,
\]
where
\begin{align*}
J_{\alpha,\eta}(\tau)
:=\int_{\mathbb R}\int_{\mathcal M^d}\int_{\mathcal M^d}
&\rho_\eta(r)K_\alpha(x,y)
\widetilde U\left(
 v\left(x,\tau+\frac r2\right),
 u\left(y,\tau-\frac r2\right)
\right)
\\[-1mm]
&\,dV_g(x)\,dV_g(y)\,dr.
\end{align*}

Define the finite-scale distance
\[
    \mathcal E_\alpha(\tau)
    :=
    \int_{\mathcal M^d}\int_{\mathcal M^d}
    K_\alpha(x,y)
    \widetilde U(v(x,\tau),u(y,\tau))
    \,dV_g(x)dV_g(y).
\]
Since \(\int_{\mathbb R}\rho_\eta(r)\,dr=1\), we can write
\[
\mathcal E_\alpha(\tau)
=\int_{\mathbb R}\int_{\mathcal M^d}\int_{\mathcal M^d}\rho_\eta(r)K_\alpha(x,y)\widetilde U(v(x,\tau),u(y,\tau))dV_g(x)dV_g(y)dr.
\]
Therefore, using \(\widetilde U(a,b)=|a-b|\), we obtain
\[
\begin{aligned}
&|J_{\alpha,\eta}(\tau)-\mathcal E_\alpha(\tau)|
\\
&\le
\int_{\mathbb R}
\int_{\mathcal M^d}
\int_{\mathcal M^d}
\rho_\eta(r)K_\alpha(x,y)
\left|
v\left(x,\tau+\frac r2\right)-v(x,\tau)
\right|
\,dV_g(x)dV_g(y)dr
\\
&\quad+
\int_{\mathbb R}
\int_{\mathcal M^d}
\int_{\mathcal M^d}
\rho_\eta(r)K_\alpha(x,y)
\left|
u\left(y,\tau-\frac r2\right)-u(y,\tau)
\right|
\,dV_g(x)dV_g(y)dr
\\
&=:I_1+I_2.
\end{aligned}
\]

We first estimate \(I_1\). By the definition $Z_\alpha(x) =\int_{\mathcal M^d}K_\alpha(x,y)\,dV_g(y),$ and by \eqref{eq:requ-kernel-normalization}, we have $\sup_{x\in\mathcal M^d}Z_\alpha(x)\le C$ for all \(0<\alpha\le\alpha_0\). Hence,
\[
\begin{aligned}
I_1
&=
\int_{\mathbb R}
\rho_\eta(r)
\int_{\mathcal M^d}
Z_\alpha(x)
\left|
v\left(x,\tau+\frac r2\right)-v(x,\tau)
\right|
\,dV_g(x)dr
\\
&\le
C
\int_{\mathbb R}
\rho_\eta(r)
\left\|
v\left(\cdot,\tau+\frac r2\right)-v(\cdot,\tau)
\right\|_{L_1(\mathcal M^d)}
\,dr .
\end{aligned}
\]
Since \(\rho_\eta\) is supported in \([-\,\eta,\eta]\), we have $\left|\left(\tau+\frac r2\right)-\tau\right| \le \frac{\eta}{2}.$ Consequently,
\[
    \left\|
    v\left(\cdot,\tau+\frac r2\right)-v(\cdot,\tau)
    \right\|_{L_1(\mathcal M^d)}
    \le
    \omega_t(v;\eta/2)
    \le
    \omega_t(v;C\eta).
\]
Since \(\int_{\mathbb R}\rho_\eta(r)\,dr=1\), we get $I_1\le C\omega_t(v;C\eta).$

We next estimate \(I_2\). By the symmetry of \(K_\alpha\),
\[
\int_{\mathcal M^d}K_\alpha(x,y)\,dV_g(x)
=
\int_{\mathcal M^d}K_\alpha(y,x)\,dV_g(x)
=
Z_\alpha(y).
\]
Again by \eqref{eq:requ-kernel-normalization},
\[
\sup_{y\in\mathcal M^d}Z_\alpha(y)\le C.
\]
Therefore,
\[
\begin{aligned}
I_2
&=
\int_{\mathbb R}
\rho_\eta(r)
\int_{\mathcal M^d}
Z_\alpha(y)
\left|
u\left(y,\tau-\frac r2\right)-u(y,\tau)
\right|
\,dV_g(y)dr
\\
&\le
C
\int_{\mathbb R}
\rho_\eta(r)
\left\|
u\left(\cdot,\tau-\frac r2\right)-u(\cdot,\tau)
\right\|_{L_1(\mathcal M^d)}
\,dr .
\end{aligned}
\]
By~\eqref{equation: bounded variation inequality}, we obtain $\left\|u\left(\cdot,\tau-\frac r2\right)-u(\cdot,\tau)\right\|_{L_1(\mathcal M^d)}\le\frac{|r|}{2}TV(u_0).$
Hence,
\[
I_2
\le
C
\int_{\mathbb R}
\rho_\eta(r)\frac{|r|}{2}TV(u_0)\,dr\le
C\eta\,TV(u_0),
\]
where we used \(\operatorname{supp}\rho_\eta\subset[-\eta,\eta]\) and
\(\int_{\mathbb R}\rho_\eta(r)\,dr=1\).

Combining the estimates for \(I_1\) and \(I_2\), we obtain
\begin{equation}\label{eq:J-Ealpha-estimate}
    |J_{\alpha,\eta}(\tau)-\mathcal E_\alpha(\tau)|
    \le
    C\omega_t(v;C\eta)
    +
    C\eta\,TV(u_0).
\end{equation}

It remains to compare \(\mathcal E_\alpha(\tau)\) with the true error $ \mathcal E(\tau)=\|v(\cdot,\tau)-u(\cdot,\tau)\|_{L_1(\mathcal M^d)}.$
Using the triangle inequality and the definition of \(Z_\alpha\), we have
\[
\begin{aligned}
|\mathcal E_\alpha(\tau)-\mathcal E(\tau)|
&\le
\int_{\mathcal M^d}\int_{\mathcal M^d}
K_\alpha(x,y)|u(y,\tau)-u(x,\tau)|
\,dV_g(x)dV_g(y)
\\
&\quad+
\int_{\mathcal M^d}
|v(x,\tau)-u(x,\tau)|\,|Z_\alpha(x)-1|
\,dV_g(x).
\end{aligned}
\]
The first term is controlled by the \(BV\)-localization estimate
\eqref{eq:requ-bv-localization}:
\[
    \int_{\mathcal M^d}\int_{\mathcal M^d}
    K_\alpha(x,y)|u(y,\tau)-u(x,\tau)|
    \,dV_g(x)dV_g(y)
    \le
    C\alpha\,TV(u(\cdot,\tau)).
\]
Using~\eqref{equation: bounds on TV}, we get $C\alpha\,TV(u(\cdot,\tau))\le C_T\alpha(1+TV(u_0)).$
For the second term, \eqref{eq:requ-kernel-normalization} gives $|Z_\alpha(x)-1|\le C\alpha^2.$ Since \(v,u\in L_\infty(\mathcal M^d\times[0,T])\), it follows that $\int_{\mathcal M^d}|v(x,\tau)-u(x,\tau)|\,|Z_\alpha(x)-1|dV_g(x)\le C\alpha^2\le C\alpha.$
Therefore,
\begin{equation}\label{eq:Ealpha-E-estimate}
    |\mathcal E_\alpha(\tau)-\mathcal E(\tau)|
    \le
    C_T\alpha(1+TV(u_0)).
\end{equation}
Combining \eqref{eq:J-Ealpha-estimate} and \eqref{eq:Ealpha-E-estimate}, we obtain
\begin{equation*}
\begin{aligned}
|J_{\alpha,\eta}(\tau)-\mathcal E(\tau)|
\le
C_T\Big[
    \omega_t(v;C\eta)
    +
    \eta\,TV(u_0)
    +
    \alpha(1+TV(u_0))
\Big].
\end{aligned}
\end{equation*}
Since \(\|\chi_h'\|_{L_1(0,T)}\le C\), it follows that
\[
\left|
A_\lambda(v)
-
\int_0^T\chi_h'(\tau)\mathcal E(\tau)\,d\tau
\right|
\le
C_T
\left[
\omega_t(v;C\eta)
+
\eta\,TV(u_0)
+
\alpha(1+TV(u_0))
\right].
\]

Next compare $\int_0^T\chi_h'(\tau)\mathcal E(\tau)\,d\tau$
with \(\mathcal E(0)-\mathcal E(T)\). Since $\chi_h=0$ near $0,T,\chi_h=1$ on $[2h,T-2h],$ and $\operatorname{supp}\chi_h'
\subset [h,2h]\cup[T-2h,T-h],$
we have $\int_h^{2h}\chi_h'(\tau)\,d\tau=1, \int_{T-2h}^{T-h}\chi_h'(\tau)\,d\tau=-1.$
Hence,
\[
\begin{aligned}
&\left|
\int_0^T\chi_h'(\tau)\mathcal E(\tau)\,d\tau
-
\bigl(\mathcal E(0)-\mathcal E(T)\bigr)
\right|
\\
&\qquad\le
C
\sup_{\tau\in[h,2h]}|\mathcal E(\tau)-\mathcal E(0)|
+
C
\sup_{\tau\in[T-2h,T-h]}|\mathcal E(\tau)-\mathcal E(T)|.
\end{aligned}
\]
For any \(t,s\in[0,T]\),
\[
\begin{aligned}
|\mathcal E(t)-\mathcal E(s)|
&\le
\|v(\cdot,t)-v(\cdot,s)\|_{L_1(\mathcal M^d)}
+
\|u(\cdot,t)-u(\cdot,s)\|_{L_1(\mathcal M^d)}
\\
&\le
\omega_t(v;|t-s|)
+
TV(u_0)|t-s|.
\end{aligned}
\]
Therefore,
\[
\left|
\int_0^T\chi_h'(\tau)\mathcal E(\tau)\,d\tau
-
\bigl(\mathcal E(0)-\mathcal E(T)\bigr)
\right|
\le
C\omega_t(v;Ch)
+
Ch\,TV(u_0).
\]
Combining the preceding estimates gives
\begin{equation}\label{eq:A-lambda-estimate}
\begin{aligned}
\left|
A_\lambda(v)-\bigl(\mathcal E(0)-\mathcal E(T)\bigr)
\right|
\le
C_T
\Big[
&\omega_t(v;C\eta)
+
\omega_t(v;Ch)
\\
&+
(\eta+h)TV(u_0)
+
\alpha(1+TV(u_0))
\Big].
\end{aligned}
\end{equation}

It remains to estimate \(B_\lambda(v)\). Since $\chi_h\left(\frac{t+s}{2}\right)\rho_\eta(t-s)$
does not depend on \(x,y\), we can write
\[
  B_\lambda(v)
  =
  \int_0^T\int_0^T
  \chi_h\left(\frac{t+s}{2}\right)
  \rho_\eta(t-s)
  \mathcal B_\alpha(t,s)
  \,dt\,ds,
\]
where
\[
\begin{aligned}
\mathcal B_\alpha(t,s)
:=
\int_{\mathcal M^d}\int_{\mathcal M^d}
\Big[
&g_x\bigl(
\widetilde F_x(v(x,t),u(y,s)),
\nabla_xK_\alpha(x,y)
\bigr)
\\
&+
g_y\bigl(
\widetilde F_y(v(x,t),u(y,s)),
\nabla_yK_\alpha(x,y)
\bigr)
\Big]
\,dV_g(x)dV_g(y).
\end{aligned}
\]
We split
\[
    \mathcal B_\alpha(t,s)
    =
    \mathcal I_\alpha(t,s)+\mathcal R_\alpha(t,s),
\]
where \(\mathcal I_\alpha(t,s)\) is obtained by replacing \(u(y,s)\) with
\(u(x,s)\):
\[
\begin{aligned}
\mathcal I_\alpha(t,s)
:=
\int_{\mathcal M^d}\int_{\mathcal M^d}
\Big[
&g_x\bigl(
\widetilde F_x(v(x,t),u(x,s)),
\nabla_xK_\alpha(x,y)
\bigr)
\\
&+
g_y\bigl(
\widetilde F_y(v(x,t),u(x,s)),
\nabla_yK_\alpha(x,y)
\bigr)
\Big]
\,dV_g(x)dV_g(y).
\end{aligned}
\]

The first term in \(\mathcal I_\alpha(t,s)\) is no longer exactly zero,
because \(K_\alpha\) is only approximately normalized. For fixed \(x\),
\[
\begin{aligned}
&\int_{\mathcal M^d}
g_x\bigl(
\widetilde F_x(v(x,t),u(x,s)),
\nabla_xK_\alpha(x,y)
\bigr)
\,dV_g(y)
\\
&=
g_x\left(
\widetilde F_x(v(x,t),u(x,s)),
\nabla_x Z_\alpha(x)
\right).
\end{aligned}
\]
By~\eqref{eq:requ-kernel-normalization-gradient} and the \(L_\infty\)-bounds
on \(v\) and \(u\), this contribution is bounded by \(C\alpha\) after
integration in \(x\).

For the second term, fix \(x\). The map $p\mapsto \widetilde F_p(v(x,t),u(x,s))$
is divergence-free. Indeed, $\widetilde F_p(v(x,t),u(x,s)) = \operatorname{sgn}(v(x,t)-u(x,s))\bigl(f_p(v(x,t))-f_p(u(x,s))\bigr),$
and the two state variables are constants with respect to \(p\). Hence, by
the geometry-compatible condition~\eqref{equation: constant flux},
\[
    \operatorname{div}_g\widetilde F_p(v(x,t),u(x,s))=0.
\]
Since \(\mathcal M^d\) has no boundary, integration by parts gives
\[
\int_{\mathcal M^d}
g_y\bigl(
\widetilde F_y(v(x,t),u(x,s)),
\nabla_yK_\alpha(x,y)
\bigr)
\,dV_g(y)
=
0.
\]
Therefore,
\begin{equation}\label{eq:I-alpha-ReQU-estimate}
    |\mathcal I_\alpha(t,s)|
    \le C\alpha.
\end{equation}

We now estimate the remainder. Set, only within this estimate,
\begin{align*}
\Delta_x
&:=\widetilde F_x(v(x,t),u(y,s))
-\widetilde F_x(v(x,t),u(x,s)),
\\
\Delta_y
&:=\widetilde F_y(v(x,t),u(y,s))
-\widetilde F_y(v(x,t),u(x,s)).
\end{align*}
Then
\[
\begin{aligned}
\mathcal R_\alpha(t,s)
=\int_{\mathcal M^d}\int_{\mathcal M^d}
\Big[
&g_x\bigl(\Delta_x,\nabla_xK_\alpha(x,y)\bigr)
\\
&+g_y\bigl(\Delta_y,\nabla_yK_\alpha(x,y)\bigr)
\Big]
\,dV_g(x)\,dV_g(y).
\end{aligned}
\]
After parallel transporting the second term to \(T_x\mathcal M^d\), we write
\[
\begin{aligned}
&g_x(\Delta_x,\nabla_xK_\alpha)
+
g_y(\Delta_y,\nabla_yK_\alpha)
\\
&=
g_x\left(
\Delta_x,
\nabla_xK_\alpha
+
P_{y\to x}\nabla_yK_\alpha
\right)
+
g_x\left(
P_{y\to x}\Delta_y-\Delta_x,
P_{y\to x}\nabla_yK_\alpha
\right).
\end{aligned}
\]
By the smoothness of the flux and the boundedness of \(v\) and \(u\),
\[
    |\Delta_x|_{g_x}
    \le
    C|u(y,s)-u(x,s)|.
\]
Using~\eqref{eq:requ-kernel-cancellation}, the first term is bounded by
\[
    C|u(y,s)-u(x,s)|
    \alpha^{-d}
    \mathbf 1_{\{d_g(x,y)\le C\alpha\}}.
\]
For the second term, let $\mathcal K\subset\mathbb R$ be a compact interval containing the ranges of $u$ and $v$, and, for $p\in\mathcal M^d$, set
\[
G_p(r;q,q')
:=
\widetilde F_p(r,q)-\widetilde F_p(r,q').
\]
The boundedness of $\nabla_x\partial_u f$ on $\mathcal M^d\times\mathcal K$, together with the elementary Kruzkov-flux Lipschitz estimate applied to the vector field $\nabla_x f$, gives
\[
\sup_{p\in\mathcal M^d}
|\nabla_pG_p(r;q,q')|_g
\le
C_{\mathcal K}|q-q'|,
\qquad r,q,q'\in\mathcal K.
\]
Let $\gamma:[0,d_g(x,y)]\to\mathcal M^d$ be the unique unit-speed minimizing geodesic from $x$ to $y$. Differentiating the parallel-transported vector field along $\gamma$ yields
\[
\begin{aligned}
|P_{y\to x}G_y(r;q,q')-G_x(r;q,q')|_{g_x}
&\le
\int_0^{d_g(x,y)}
|\nabla_{\dot\gamma(\ell)}G_{\gamma(\ell)}(r;q,q')|_g
\,d\ell\\
&\le
C_{\mathcal K}|q-q'|\,d_g(x,y).
\end{aligned}
\]
Taking $r=v(x,t)$, $q=u(y,s)$, and $q'=u(x,s)$ gives
\[
    |P_{y\to x}\Delta_y-\Delta_x|_{g_x}
    \le
    C|u(y,s)-u(x,s)|\,d_g(x,y).
\]
Together with~\eqref{eq:requ-kernel-gradient} and the support condition
\(d_g(x,y)\le C\alpha\), this yields
\[
\begin{aligned}
&\left|
g_x\left(
P_{y\to x}\Delta_y-\Delta_x,
P_{y\to x}\nabla_yK_\alpha
\right)
\right|
\\
&\qquad\le
C|u(y,s)-u(x,s)|
d_g(x,y)
\alpha^{-d-1}
\mathbf 1_{\{d_g(x,y)\le C\alpha\}}
\\
&\qquad\le
C|u(y,s)-u(x,s)|
\alpha^{-d}
\mathbf 1_{\{d_g(x,y)\le C\alpha\}}.
\end{aligned}
\]
Therefore,
\[
\begin{aligned}
|\mathcal R_\alpha(t,s)|
&\le
C
\int_{\mathcal M^d}\int_{\mathcal M^d}
\alpha^{-d}
\mathbf 1_{\{d_g(x,y)\le C\alpha\}}
|u(y,s)-u(x,s)|
\,dV_g(x)dV_g(y).
\end{aligned}
\]
By~\eqref{eq:requ-indicator-bv-localization},
\[
    |\mathcal R_\alpha(t,s)|
    \le
    C\alpha\,TV(u(\cdot,s)).
\]
Combining this with~\eqref{eq:I-alpha-ReQU-estimate}, we obtain $|\mathcal B_\alpha(t,s)| \le C\alpha\bigl(1+TV(u(\cdot,s))\bigr).$
Using~\eqref{equation: bounds on TV}, we get $|\mathcal B_\alpha(t,s)|\le C_T\alpha(1+TV(u_0)).$
Since \(0\le\chi_h\le1\), \(\rho_\eta\ge0\), and $\int_0^T\rho_\eta(t-s)\,dt\le1,$
we conclude that $|B_\lambda(v)|\le C_T\alpha(1+TV(u_0)).$

Combining this estimate with~\eqref{eq:A-lambda-estimate}, we obtain
\[
b_{\lambda,T}(v)
\le
C_T
\Big[
\omega_t(v;C\eta)
+
\omega_t(v;Ch)
+
(\eta+h)TV(u_0)
+
\alpha(1+TV(u_0))
\Big].
\]
Since \(v\in C([0,T];L_1(\mathcal M^d))\), we have $\omega_t(v;r)\to0 $ as $r\downarrow0.$ Thus, $b_{\lambda,T}(v)\to0 $ as $\lambda=(\alpha,\eta,h)\to0.$
Finally, if the constructed network satisfies
\[
    \omega_t(v;r)\le C(1+TV(u_0))r,
\]
then \(b_{\lambda,T}(v) \le C_T(1+TV(u_0))(\alpha+\eta+h)\). This proves the terminal estimate for the fixed horizon \(T\). Replacing \(T\) by an arbitrary \(\tau\in[4h,T]\) gives \eqref{equation: horizon terminal stability}. Taking \(\tau=T\) and using
\[
\mathcal L_{\lambda,T}^+(v)\le \mathfrak L_\lambda^+(v)
\]
gives \eqref{equation: terminal stability global positive}.

It remains only to derive the space-time form. Integrating the horizon-indexed estimate over \(\tau\in[4h,T]\) yields
\[
\int_{4h}^T\|v(\cdot,\tau)-u(\cdot,\tau)\|_{L_1(\mathcal M^d)}\,d\tau
\le
T\Bigl[
\mathcal R^{tb}(v)
+
T\operatorname{Vol}_g(\mathcal M^d)\mathfrak L_\lambda^+(v)
+
b_{\lambda,T}(v)
\Bigr].
\]
For \(0\le \tau\le4h\),
\[
\|v(\cdot,\tau)-u(\cdot,\tau)\|_{L_1}
\le
\mathcal R^{tb}(v)+\omega_t(v;4h)+4h\,TV(u_0),
\]
which is absorbed into the same right-hand side by the definition of \(b_{\lambda,T}\). Dividing by \(T\operatorname{Vol}_g(\mathcal M^d)\) and absorbing fixed geometric constants proves \eqref{equation: spacetime stability global positive}.
\end{proof}

\section{Analysis of Approximation Error}\label{Appendix: Analysis of Approximation Error}
We begin by fixing the H{\"o}lder notation used for the smooth coordinate maps. For $s>0$, set
\[
m_s:=\lceil s\rceil-1,
\qquad
\beta_s:=s-m_s\in(0,1].
\]
For $x,a\in\mathbb R^D$ and a multi-index $\gamma\in\mathbb N^D$, write
\[
x^\gamma=x_1^{\gamma_1}\cdots x_D^{\gamma_D},
\qquad
\gamma!=\gamma_1!\cdots\gamma_D!,
\qquad
|\gamma|=\|\gamma\|_1.
\]
The Taylor polynomial of H{\"o}lder order $s$ at $a$ is
\[
P_a^sf(x)
=
\sum_{|\gamma|\le m_s}
\frac{\partial^\gamma f(a)}{\gamma!}(x-a)^\gamma.
\]
For a bounded domain $U\subset\mathbb R^D$, set
\[
[f]_{m_s,\beta_s;U}
:=\max_{|\gamma|=m_s}
\sup_{\substack{x,y\in U\\x\ne y}}
\frac{|\partial^\gamma f(x)-\partial^\gamma f(y)|}
{\|x-y\|_\infty^{\beta_s}},
\]
and define
\begin{equation*}
\mathcal H_D^s(U,K)
:=\left\{
 f\in C^{m_s}(\overline U):
 \sum_{|\gamma|\le m_s}\|\partial^\gamma f\|_{L_\infty(U)}
 +[f]_{m_s,\beta_s;U}
 \le K
\right\}.
\end{equation*}
Set $\mathcal H_D^s(U):=\bigcup_{K>0}\mathcal H_D^s(U,K)$. For vector-valued maps, the definition is imposed componentwise.

A smooth compact embedded manifold admits a finite atlas and a partition of unity with nested supports. We record the precise structure used in the proof.
\begin{definition}[Smooth local coordinates]
\label{definition: smooth coordinates}
We say that the compact embedded manifold $\mathcal M^d\subset\mathbb R^D$ admits smooth local coordinates if there exist finitely many charts $(V_i,\psi_i)$, open sets $V_i'$ and a smooth partition of unity $\{\rho_i\}_{i=1}^k$ such that
\[
\operatorname{supp}\rho_i\Subset V_i'\Subset V_i,
\qquad
\bigcup_{i=1}^kV_i'=\mathcal M^d,
\]
and the following properties hold:
\begin{enumerate}
\item $\psi_i:V_i\to\psi_i(V_i)\subset\mathbb R^d$ is a smooth diffeomorphism, $\psi_i(V_i)$ is a bounded Lipschitz domain, and $\psi_i^{-1}$ is smooth;
\item there is a fixed bounded ambient box $\mathcal B_i\subset\mathbb R^D$ containing $\mathcal M^d$, an open neighborhood $\mathcal O_i\Subset\mathcal B_i$ satisfying $\overline{V_i'}\subset\mathcal O_i\cap\mathcal M^d\Subset V_i$, and a smooth map $\Psi_i:\mathcal B_i\to\mathbb R^d$ satisfying $\Psi_i|_{\mathcal O_i\cap\mathcal M^d}=\psi_i|_{\mathcal O_i\cap\mathcal M^d}$;
\item there is a smooth ambient cutoff $\vartheta_i\in C_c^\infty(\mathcal O_i)$ satisfying $0\le\vartheta_i\le1$ and $\vartheta_i=1$ on $\operatorname{supp}\rho_i$.
\end{enumerate}
All derivatives of these finitely many maps are bounded on the compact sets used below.
\end{definition}
Normal-coordinate balls may be chosen so that the chart images are smooth bounded domains. The ambient extensions and cutoffs then follow from the tubular-neighborhood theorem and a finite partition of unity; see~\cite{lee2003introduction}. Thus every smooth compact embedded manifold has the structure in \autoref{definition: smooth coordinates}.

\textbf{The Proof Sketch of~\autoref{theorem: approximation error for sobolev functions on manifolds}}. The key idea of the proof lies in leveraging the composite structure of functions defined on the manifold, specifically the decomposition $f(x) = \sum_{i=1}^k(\rho_i f)\circ \psi_i^{-1} \circ \psi_i(x)$. The approximation strategy proceeds in two main stages. First, we approximate the coordinate mappings $\{\psi_i(x)\}_{i=1}^k$ and the dimension-reduced functions $\{(\rho_i f)\circ \psi_i^{-1}\}_{i=1}^k$ separately. Then, by combining parallel connections with concatenation operations within the neural network architecture, we construct the final approximation. In the special case where the domain is $[0,1]^D$, the coordinate map reduces to the identity, and the standard approach is to approximate the function using a locally averaged Taylor expansion. Specifically, within each cell of a uniform grid, neural networks are used to approximate the corresponding Taylor polynomial. Subsequently, another network is constructed to produce a partition of unity that assembles these local approximations. For Sobolev functions defined on manifolds, we begin by extending the results from the Euclidean setting $[0,1]^d$ to more general domains. Then, by exploiting the naturally defined composition structure and the intrinsic dimensionality reduction offered by the manifold, we obtain the final approximation result.

One of the principal challenges in the argument is handling the approximation across different function spaces. When working with the H{\"o}lder space, the composite structure naturally leads to an approximation with respect to the $L_\infty$ norm, as the H{\"o}lder space is inherently defined using this norm. Moreover, the $L_\infty$ norm also dominates the $L_p$ norm, making the transition straightforward. However, for general Sobolev functions, if we aim to approximate in the $L_p$ norm using the composite structure, we must control the $L_\infty$ norm of the derivatives of the target function—a task complicated by the fact that Sobolev space definitions do not directly provide such control. In contrast, during the constructive approximation process using neural networks, we can explicitly analyze and bound these derivative norms, thereby ensuring that the total approximation error remains well-controlled.

We divide our proof into three parts.

\textbf{Part 1}: We establish neural network approximation rates for functions in $\mathcal{H}_D^s\left([0,1]^D,K\right)$ and $W_p^s\left((0,1)^d\right)$. We first construct a global approximator, which is formed by combining local approximation functions with a partition of unity. Subsequently, we approximate this global approximator using neural networks. We will present the results separately for H{\"o}lder and Sobolev functions.

We begin by constructing the local approximation polynomial for $f \in \mathcal{H}_D^s([0,1]^D, K)$ at a given point $x \in [0,1]^D$:
\begin{equation*}
    P_a^s f(x)=\sum_{\gamma:|\gamma|<s}(\partial^{\gamma}f)(a)(x-a)^{\gamma}/\gamma!.
\end{equation*}
The H{\"o}lder--Taylor remainder estimate associated with the convention above gives
\begin{equation*}
    \left|f(x)-P_a^sf(x)\right|
    \leq C(D,s)K\|x-a\|_\infty^s.
\end{equation*}
For noninteger $s$, this is the usual Taylor formula with a $\beta_s$-H{\"o}lder remainder. For integer $s$, it is the Taylor estimate of degree $s-1$ under the Lipschitz control of derivatives of order $s-1$.
We rewrite the polynomial as $P_a^sf(x) = \sum_{|\gamma| <s} c_{\gamma,a}x^\gamma $. It follows that
\begin{equation*}
    |c_{\gamma,a}| \leq  \sum_{\substack{\alpha:0\leq |\alpha|<s \\  \gamma \leq \alpha  }}|\partial^\alpha f(a)| \left|\frac{(-a)^{\alpha-\gamma}}{\gamma!(\alpha-\gamma)!}\right| \leq C(D,s)K, \quad \sum_{|\gamma|<s } |c_{\gamma,a}| \leq C(D,s)K.
\end{equation*}
Next, we construct a partition of unity. Let $N \in \mathbb{N}_+$ and define $[N] = \{0,1,\ldots, N\}$. We introduce a uniform grid $D(N)=\left\{ x_l:(x_l)_j =l_j/N ,l=(l_1,l_2,\ldots,l_D) \in [N]^D\right\}$. The total number of grid points satisfies $|D(N)| = (N+1)^D.$ For each grid point $x_l$, we define the partition of unity function by $\rho^D_{x_l}(x) = \prod _{j=1}^{D}(1-N|x_j - (x_l)_j |)_+$. In this construction, the support of $\rho^D_{x_l}$ is given by $\mathrm{supp} \rho^D_{x_l} = B(x_l,1/N,\|\cdot\|_\infty)$, i.e., the closed ball of radius $1/N$ around $x_l$ under the infinity norm. Moreover, we have
\begin{equation}\label{equation: partion of unity}
    \sum_{x_l \in D(N)} \rho^D_{x_l}(x) = \prod_{j=1}^D \sum_{k = 0}^N (1-N|x_j - k/N |)_+ = 1.
\end{equation}
We now construct a global approximator $P^s_Nf(x)$ for $f(x)$ as
\begin{equation*}
    P_N^sf(x) = \sum_{x_l \in D(N)} \rho^D_{x_l}(x) \cdot P_{x_l}^sf(x) = \sum_{x_l \in D(N)} P_{x_l}^sf(x)\prod _{j=1}^{D}(1-N|x_j - (x_l)_j |)_+.
\end{equation*}
From \eqref{equation: partion of unity}, we obtain the following bound for all $x \in [0,1]^D$:
\begin{equation*}
    \begin{aligned}
        |f(x)- P_N^sf(x)| &\leq \sum_{x_l \in D(N)} \rho^D_{x_l}(x) \cdot |(f(x)-P_{x_l}^sf(x))| \\
        &\leq \sum_{x_l \in D(N):\|x-x_l\|_\infty \leq 1/N}\rho^D_{x_l}(x) K\|x-x_l\|^s_\infty\\
        &\leq KN^{-s}.
    \end{aligned}
\end{equation*}

Let $s \in \mathbb{N}_+$. For $f \in W_p^s\left((0,1)^d\right)$, the construction of the global approximator is similar. Since the weak derivatives of $f$ cannot be defined point wise, we introduce the notion of the averaged Taylor expansion. Let $U \subset B(0,R,\|\cdot\|_\infty)$ be a bounded open set in $\mathbb{R}^d$, and let $f \in W_p^s(U)$. Suppose $x_0 \in U$ and $r > 0$ are such that the ball $B(x_0, r, \|\cdot\|_2) \subset U$. The averaged Taylor expansion over $B(x_0,r,\|\cdot\|_2)$ is defined as
\begin{equation*}
    Q^sf(x)=  \int_{B(x_0,r,\|\cdot\|_2)} \lambda(y) P_{y}^{s}f(x) dy,
\end{equation*}
where $\lambda(y)$ is a cut-off function supported on $\overline{B(x_0,r,\|\cdot\|_2)}$. To construct $\lambda(y)$, we first define an auxiliary function $\kappa: \mathbb{R}^d\to \mathbb{R}$ by
\begin{equation*}
    \kappa(x)=\begin{cases}e^{-\left(1-(\|x-x_0\|_2/r)^2\right)^{-1}},&\text{if }\|x-x_0\|_2<r\\0,&\text{else.}\end{cases}
\end{equation*}
Normalizing $\kappa$ by $\int_{\mathbb{R}^d} \kappa(x)dx$ yields the desired cut-off function $\lambda(x)$.  After a change of variables in the integral, we obtain the estimate $\|\lambda(x)\|_{L_\infty(\mathbb{R}^d)} \leq c r^{-d}$.  We express $Q^s f(x)$ as a linear combination of monomials:
\begin{equation*}
    Q^sf(x) =  \sum_{|\alpha| < s} c_\alpha x^\alpha, 
\end{equation*}
where
\begin{equation*}
     c_\alpha =\sum_{\substack{\gamma:0\leq |\gamma|<s \\ \alpha \leq \gamma }} \int_{B(x_0,r,\|\cdot\|_2)} \frac{\partial^\gamma f(y)(-y)^{\gamma -\alpha} }{\alpha!(\gamma-\alpha)!} \lambda(y)dy.
\end{equation*}
We can estimate the coefficients as
\begin{equation*}
    \begin{aligned}
            |c_\alpha| &\leq\sum_{\substack{\gamma:0\leq |\gamma|<s \\ \alpha \leq \gamma }}   R^s   \int_{B(x_0,r,\|\cdot\|_2)}\left | \frac{\partial^\gamma f(y)}{\alpha!(\gamma-\alpha)!} \lambda(y)\right |dy\\
            &\leq R^s \|f\|_{W^s_p(U)}\|\lambda(y)\|_{L_q(\mathbb{R}^d)}\sum_{\substack{\gamma:0\leq |\gamma|<s \\ \alpha \leq \gamma }} \frac{1}{\alpha!(\gamma-\alpha)!}\\
            &\leq C(s)R^s\|f\|_{W^s_p(U)}\|\lambda(y)\|^{\frac{1}{p}}_{L_\infty(\mathbb{R}^d)}\|\lambda(y)\|^{\frac{1}{q}}_{L_1(\mathbb{R}^d)} \\
            & \leq C(s)\|f\|_{W^s_p(U)} r^{-\frac{d}{p}}.
    \end{aligned}
\end{equation*}
Using the Taylor expansion, we can now apply the Bramble–Hilbert Lemma to obtain a local approximation. Before doing so, we introduce the following definitions.
\begin{definition}
    We say that $U$ is star-shaped with respect to $B$ if, for all $x \in U$, the closed convex hull of $\{x\}\cup B$ is a subset of $U$, i.e., $\overline{\mathrm{conv}}(\{x\}\cup B) \subset U$.
\end{definition}
\begin{definition}
    Let $\mathrm{diam}(U)$ be the diameter of $U$, and assume that $U$ is star-shaped with respect to a ball $B$. Define $\eta_{\max} = \sup\{\eta \ | \  U \text{ is star-shaped with respect to a ball }B \text{ of radius } \eta  \} $. Then, the chunkiness parameter of $U$ is defined as
    \begin{equation*}
        \zeta = \frac{\mathrm{diam}(U)}{\eta_{\max}}.
    \end{equation*}
\end{definition}
Finally, the following lemma establishes the approximation properties of the averaged Taylor polynomial.
\begin{lemma}[Bramble-Hilbert Lemma] \label{lemma: Bramble-Hilbert}
   Let $f \in W_p^s(U)$ with $s \in \mathbb{N}_+$ and $1\le p<\infty$. Let $B$ be a ball contained in $U$ such that $U$ is star-shaped with respect to $B$, and let the radius $r$ of $B$ satisfy $r > \frac{1}{2} \eta_{\max}$. Let $Q^s f$ denote the averaged Taylor expansion polynomial of $f$ on $B$. Then,
   \begin{equation*}
       \|f-Q^sf\|_{L_p(U)} \leq C(s,d,\zeta) (\mathrm{diam}(U))^s |f|_{W_p^s(U)}.
   \end{equation*}
   \end{lemma}
A proof can be found in~\cite[Lemma 4.3.8]{brenner2008mathematical}. Next, we construct a global approximator. For $f \in W_p^s([0,1]^d)$, we utilize the Sobolev extension theorem (\cite{stein1970singular}, Theorem VI.3.1.5) to extend $f$ to the $\mathbb{R}^d$. Denote the extension by $\widetilde f:=Ef$, where $E:W_p^s((0,1)^d)\to W_p^s(\mathbb R^d)$ is bounded. For the grid points $x_l$, $l\in[N]^d$, define $U_l=B(x_l,1/N,\| \cdot \|_\infty)$ and $B_l= B(x_l, 3/(4N),\| \cdot \|_{2})$. It is evident that each $U_l$ is star-shaped with respect to $B_l$, $\mathrm{diam}(U_l) = 2\sqrt{d}/N$, and $\eta_{\max} = 1/N$.  Therefore, the corresponding chunkiness parameter is $\zeta = 2\sqrt{d}$. Note also that $3/(4N) > (1/2)\eta_{\max}$. Let $Q^s_{x_l}\widetilde f = \sum_{|\alpha |< s} c_{l,\alpha} x^\alpha$ be the averaged Taylor expansion polynomial of the extension $\widetilde f$ on $U_l$ with respect to $B_l$. By applying~\autoref{lemma: Bramble-Hilbert}, we obtain the desired local approximation property:
\begin{equation*}
     \|\widetilde f-Q^s_{x_l}\widetilde f\|_{L_p(U_l)} \leq C(s,d)|\widetilde f|_{W_p^s(U_l)} N^{-s}.
\end{equation*}
We also have the local coefficient estimate $|c_{l,\alpha}| \leq C(s,d,p)N^{d/p}\|\widetilde f\|_{W_p^s(U_l)}$. Next, we construct the global approximator by setting $Q^s_Nf = \sum_{l\in[N]^d} \rho_{x_l}^d(x)\cdot Q^s_{x_l}\widetilde f(x)$ and compute the $L_p$-error as follows:
\begin{equation*}
    \begin{aligned}
        \left\| f(x)- Q^s_Nf(x)\right\|^p_{L_p((0,1)^d)}   &=  \left\| \widetilde{f}(x)- \sum_{l \in [N]^d} \rho_{x_l}^d(x)\cdot Q^s_{x_l}\widetilde f(x)\right\|^p_{L_p((0,1)^d)}   \\
        & \leq \sum_{\tilde{l}\in [N]^d} \left\|\sum_{\substack{l \in [N]^d , \|\tilde{l}-l\|_\infty \leq 1}} \rho_{x_l}^d(x)(\widetilde{f}(x)- Q^s_{x_l}\widetilde f(x))\right\|^p_{L_p(U_{\tilde{l}})}. 
    \end{aligned}
\end{equation*}
Here, $\widetilde{f}$ denotes the Sobolev extension of $f$. By applying~\autoref{lemma: Bramble-Hilbert}, we obtain the following estimate:
\begin{equation*}
    \left\| f(x)- Q^s_Nf(x)\right\|^p_{L_p((0,1)^d)}  \leq  C(s,d,p)\cdot N^{-sp}\sum_{\tilde{l}\in [N]^d} \left(  \sum_{\substack{l \in [N]^d,  \|\tilde{l}-l\|_\infty \leq 1}} \left\|\widetilde{f}\right\|_{W^s_p(U_{l})}\right) ^p.
\end{equation*}

By further estimating the inner sum using H{\"o}lder inequality (which introduces a factor $3^{pd/q}$) and then summing over $\tilde{l}$, we obtain
\begin{equation*}
    \begin{aligned}
        \left\| f(x)- Q^s_Nf(x)\right\|^p_{L_p((0,1)^d)}  &\leq C(s,d,p)\cdot N^{-sp} \sum_{\tilde{l}\in [N]^d} \sum_{\substack{l \in [N]^d ,\|\tilde{l}-l\|_\infty \leq 1}} \left\|\widetilde{f}\right\|_{W^s_p(U_l)}^p 3^{pd/q}\\
        &\leq C(s,d,p) 3^{pd/q}N^{-sp}3^d \sum_{\tilde{l}\in[N]^d} \left\|\widetilde{f}\right\|_{W^s_p(U_{\tilde{l}})}^p
    \end{aligned}
\end{equation*}
Here, $q=p/(p-1)$ when $p>1$, while $q=\infty$ when $p=1$. Finally, accounting for an additional finite overlap factor of $2^d$, we conclude that
\begin{equation*}
        \left\| f(x)- Q^s_Nf(x)\right\|^p_{L_p((0,1)^d)} \leq C(s,d,p,E) N^{-sp} \|f\|^p_{W^s_p((0,1)^d)}.
\end{equation*}
The next step is to construct a neural network function to approximate this global approximator. Note that in both cases, the global approximators are constructed by summing the products of local approximation polynomials and a partition of unity. Therefore, the network construction process remains consistent across both settings. The following lemma is essentially adapted from the proof presented in~\cite{Shen2019DeepNA,yarotsky2017error,guhring2020error,schmidt2020nonparametric,lu2021deep}.
\begin{lemma}\label{lemma: square network}
For every $m\in\mathbb N_+$, there is a ReLU network
\[
\mathrm{SQ}_m\in\mathcal F_1(Cm,C,Cm,1,1)
\]
such that, for every $x\in[0,1]$,
\begin{align*}
|\mathrm{SQ}_m(x)-x^2|
&\le 2^{-2m-2},\\
|\mathrm{SQ}_m'(x)-2x|
&\le 2^{-m}
\qquad\text{for a.e. }x\in(0,1).
\end{align*}
Moreover, the restriction of $\mathrm{SQ}_m$ to $[0,1]$ is convex and nondecreasing,
\[
0\le \mathrm{SQ}_m(x)\le1,
\qquad
\mathrm{SQ}_m(0)=0,
\quad
\mathrm{SQ}_m(1/2)=1/4,
\quad
\mathrm{SQ}_m(1)=1.
\]
\end{lemma}

\begin{proof}
Let $q_m$ be the continuous piecewise-affine interpolant of $x\mapsto x^2$ on the dyadic grid
\[
\left\{j2^{-m}:j=0,\ldots,2^m\right\}.
\]
On the interval $[j2^{-m},(j+1)2^{-m}]$, its slope is
\[
\frac{((j+1)2^{-m})^2-(j2^{-m})^2}{2^{-m}}
=
(2j+1)2^{-m}.
\]
Consequently,
\[
\|q_m-x^2\|_{L_\infty(0,1)}\le2^{-2m-2},
\qquad
\|q_m'-2x\|_{L_\infty(0,1)}\le2^{-m}.
\]
Since $q_m$ interpolates a convex increasing function, it is itself convex and nondecreasing and takes values in $[0,1]$.

For completeness, we recall the standard sawtooth realization. Define
\[
g(x):=2\sigma(x)-4\sigma\left(x-\frac12\right)+2\sigma(x-1),
\qquad
g_0(x):=x,
\qquad
g_r:=g\circ g_{r-1}.
\]
Then
\begin{equation}\label{equation: dyadic square identity}
q_m(x)
=
x-\sum_{r=1}^m4^{-r}g_r(x).
\end{equation}
This identity follows by induction from the fact that the correction from the dyadic interpolant at level $r-1$ to that at level $r$ is $4^{-r}g_r$.

The pair
\[
(q_{r-1},g_{r-1})
\longmapsto
(q_r,g_r),
\qquad
q_r=q_{r-1}-4^{-r}g_r,
\]
can be implemented by two ReLU layers of universal width. Coefficients $2$ and $4$ in the tent map are realized by parallel copies with output weights in $\{-1,1\}$, while $4^{-r}\le1$. Hence all weights and biases are bounded by one. Concatenating the $m$ blocks realizes $q_m$ with depth at most $2m+2$, universal width, and at most $Cm$ nonzero parameters. Finally apply the global clipping map
\[
\operatorname{clip}_{[0,1]}(z)=\sigma(z)-\sigma(z-1).
\]
Because $q_m([0,1])\subset[0,1]$, this last layer does not alter the function or its derivative on $(0,1)$, but it ensures the global output bound required in the definition of $\mathcal F_1$. This proves the lemma.
\end{proof}

From \(xy=2\left[\left(\frac{x+y}{2}\right)^2-\left(\frac{x}{2}\right)^2-\left(\frac{y}{2}\right)^2\right],\) we obtain the following multiplication network.

\begin{lemma}\label{lemma: mutiple}
For every $m\in\mathbb N_+$, there is a network
\[
\mathrm{Mult}_m^2\in\mathcal F_2(Cm,C,Cm,1,1)
\]
such that, for all $(x,y)\in[0,1]^2$,
\begin{align}
|\mathrm{Mult}_m^2(x,y)-xy|
&\le C2^{-2m},
\label{equation: binary multiplication value error}\\
|\partial_1\mathrm{Mult}_m^2(x,y)-y|
+|\partial_2\mathrm{Mult}_m^2(x,y)-x|
&\le C2^{-m}
\quad\text{for a.e. }(x,y)\in(0,1)^2.
\label{equation: binary multiplication partial derivative error}
\end{align}
In addition,
\[
0\le\mathrm{Mult}_m^2(x,y)\le1,
\qquad
\mathrm{Mult}_m^2(x,0)=\mathrm{Mult}_m^2(0,y)=0.
\]
\end{lemma}

\begin{proof}
Set
\[
\widetilde{\mathrm{Mult}}_m^2(x,y)
:=
2\left[
\mathrm{SQ}_m\left(\frac{x+y}{2}\right)
-
\mathrm{SQ}_m\left(\frac{x}{2}\right)
-
\mathrm{SQ}_m\left(\frac{y}{2}\right)
\right].
\]
The value estimate follows immediately from \autoref{lemma: square network}. At points of differentiability,
\[
\partial_1\widetilde{\mathrm{Mult}}_m^2(x,y)
=
\mathrm{SQ}_m'\left(\frac{x+y}{2}\right)
-
\mathrm{SQ}_m'\left(\frac{x}{2}\right),
\]
and the analogous identity holds for $\partial_2$. Subtracting the corresponding identities for the exact square proves \eqref{equation: binary multiplication partial derivative error}.

It remains to verify the range. Since $\mathrm{SQ}_m$ is convex and vanishes at zero, it is superadditive on $[0,1]$: whenever $a,b\ge0$ and $a+b\le1$,
\[
\mathrm{SQ}_m(a+b)\ge\mathrm{SQ}_m(a)+\mathrm{SQ}_m(b).
\]
Thus $\widetilde{\mathrm{Mult}}_m^2\ge0$. Its partial derivatives are nonnegative almost everywhere because the slope of $\mathrm{SQ}_m$ is nondecreasing. Hence it is nondecreasing in each variable, and
\[
\widetilde{\mathrm{Mult}}_m^2(x,y)
\le
\widetilde{\mathrm{Mult}}_m^2(1,1)
=
2\left(1-\frac14-\frac14\right)
=1.
\]
The zero-boundary identities follow from $\mathrm{SQ}_m(0)=0$. The three square subnetworks can be parallelized, and the outer coefficient $2$ is implemented by duplication, so all parameters remain bounded by one. Applying the global clipping map to the output makes the realization globally bounded; since the unclipped realization already lies in $[0,1]$ throughout $[0,1]^2$, the clipped and unclipped functions coincide on that domain, and therefore their partial derivatives also coincide almost everywhere there. This proves all assertions.
\end{proof}

\begin{lemma}\label{lemma: muti product}
Let $k\ge2$ and $m\in\mathbb N_+$ satisfy $2^m\ge C_0k^2$ for a sufficiently large universal constant $C_0$. There is a ReLU network
\[
\mathrm{Mult}_m^k:[0,1]^k\to[0,1]
\]
with
\[
L\le Ckm,
\qquad
W\le Ck,
\qquad
S\le Ck^2m,
\qquad
B=1,
\]
such that
\begin{align}
\left|\mathrm{Mult}_m^k(x)-\prod_{i=1}^kx_i\right|
&\le C_k2^{-2m},
\label{equation: multiple multiplication value error}\\
\max_{1\le i\le k}
\left|
\partial_i\mathrm{Mult}_m^k(x)-\prod_{j\ne i}x_j
\right|
&\le C_k2^{-m}
\quad\text{for a.e. }x\in(0,1)^k.
\label{equation: multiple multiplication partial derivative error}
\end{align}
Moreover, $\mathrm{Mult}_m^k(x)=0$ whenever $\prod_{i=1}^kx_i=0$.

If $f_i:\Omega\to[0,1]$, $i=1,\ldots,k$, belong to $W_\infty^1(\Omega)$ and satisfy
\[
\max_{1\le i\le k}\|f_i\|_{W_\infty^1(\Omega)}\le N,
\]
then
\begin{align}
\left\|\mathrm{Mult}_m^k(f_1,\ldots,f_k)-\prod_{i=1}^kf_i\right\|_{L_\infty(\Omega)}
&\le C_k2^{-2m},
\label{equation: composite multiple multiplication value error}\\
\left|\mathrm{Mult}_m^k(f_1,\ldots,f_k)-\prod_{i=1}^kf_i\right|_{W_\infty^1(\Omega)}
&\le C_kN2^{-m}.
\label{equation: composite multiple multiplication derivative error}
\end{align}
\end{lemma}

\begin{proof}
Define recursively
\[
\mathrm{Mult}_m^{j+1}(x_1,\ldots,x_{j+1})
:=
\mathrm{Mult}_m^2\left(
\mathrm{Mult}_m^j(x_1,\ldots,x_j),x_{j+1}
\right).
\]
The range assertion in \autoref{lemma: mutiple} guarantees that every composition is evaluated in $[0,1]^2$. The value estimate follows by induction from
\[
|\mathrm{Mult}_m^2(a,b)-ab|\le C2^{-2m},
\qquad 0\le a,b\le1.
\]
The partial-derivative estimate follows from the chain rule and \eqref{equation: binary multiplication partial derivative error}. At each step, the derivative inherited from the previous product is multiplied by a factor bounded by $1+C2^{-m}$, while the new local error is bounded by $C2^{-m}$. The assumption $2^m\ge C_0k^2$ implies $(1+C2^{-m})^k\le C$, and induction therefore gives the stated $C_k2^{-m}$ bound. The zero property follows recursively from the corresponding binary identity. Concatenation of the binary modules gives the displayed depth, width, and sparsity bounds.

For the composite functions, apply the ordinary chain rule at points where all piecewise-affine maps are differentiable. The value estimate is unchanged, while every derivative with respect to a coordinate of $\Omega$ introduces at most one factor bounded by $N$. Combining this with \eqref{equation: multiple multiplication partial derivative error} proves \eqref{equation: composite multiple multiplication derivative error}. The weak-derivative statement follows because all functions involved are Lipschitz.
\end{proof}
In all applications below, $k\le d+s$ is fixed, so the condition $2^m\ge C_0k^2$ is absorbed into the choice of the multiplication depth. We now assemble the complete network to approximate the global polynomial.
\begin{lemma}\label{lemma: holder construction}
Let $s>0$, $K>0$, $N,m\in\mathbb N_+$, and let $f\in\mathcal H_D^s([0,1]^D,K)$. For the global Taylor--partition approximant $P_N^sf$, there exists a ReLU network $g$ with all parameters bounded by one such that
\[
\|P_N^sf-g\|_{L_\infty([0,1]^D)}
\le
C(D,s,K)2^{-2m}.
\]
Moreover,
\[
L\le C(m+\log N),
\qquad
W\le C(N+1)^D,
\qquad
S\le C(N+1)^D(m+\log N),
\]
and the output is bounded by a constant depending only on $D,s,K$.
\end{lemma}

\begin{lemma}\label{lemma: Sobolev construction}
Let $1\le p<\infty$, $s\in\mathbb N_+$, $N,m\in\mathbb N_+$, and $f\in W_p^s((0,1)^d)$. For the global averaged Taylor approximant $Q_N^sf$, there exists a ReLU network $g\in\mathcal F(L,W,S,1,\infty)$ such that
\begin{align*}
\|Q_N^sf-g\|_{L_p((0,1)^d)}
&\le
C(s,d,p)\|f\|_{W_p^s((0,1)^d)}2^{-2m},
\\
\|g\|_{W_\infty^1((0,1)^d)}
&\le
C(s,d,p)\|f\|_{W_p^s((0,1)^d)}N^{1+d/p}.
\end{align*}
The architecture satisfies
\[
L\le C(m+\log N),
\qquad
W\le C(N+1)^d,
\qquad
S\le C(N+1)^d(m+\log N),
\]
where the constants may depend on $s,d,p$ and the fixed function norm, but not on $N$ or $m$.
\end{lemma}

\begin{proof}
The constructions for the H{\"o}lder and Sobolev approximants are identical at the network-realization stage, so we give the Sobolev argument. For $l\in[N]^d$ and $i=1,\ldots,d$, set
\[
f_{l,i}(x):=(1-N|x_i-(x_l)_i|)_+.
\]
The function $f_{l,i}$ is represented exactly by a bounded-weight ReLU network of depth $O(\log N)$ by first realizing $(N^{-1}-|x_i-(x_l)_i|)_+$ and then using repeated duplication to multiply by $N$. It satisfies
\[
0\le f_{l,i}\le1,
\qquad
\|f_{l,i}\|_{W_\infty^1((0,1)^d)}\le N.
\]
For every multi-index $\alpha$ with $|\alpha|<s$, define
\[
\mathrm{Mult}_{l,\alpha,m}(x)
:=
\mathrm{Mult}_m^{d+|\alpha|}
\left(
 f_{l,1}(x),\ldots,f_{l,d}(x),
 \underbrace{x_1,\ldots,x_1}_{\alpha_1},\ldots,
 \underbrace{x_d,\ldots,x_d}_{\alpha_d}
\right).
\]
By \autoref{lemma: muti product},
\begin{align}
\left|\mathrm{Mult}_{l,\alpha,m}(x)-\rho_{x_l}(x)x^\alpha\right|
&\le C(d,s)2^{-2m},
\label{equation: local product value estimate}\\
\left|\mathrm{Mult}_{l,\alpha,m}-\rho_{x_l}x^\alpha\right|_{W_\infty^1((0,1)^d)}
&\le C(d,s)N2^{-m}.
\label{equation: estimation on derivatives}
\end{align}
The exact zero property in \autoref{lemma: muti product} also shows that the approximate local term vanishes outside the support of $\rho_{x_l}$.

Write
\[
Q_N^sf(x)
=
\sum_{l\in[N]^d}
\sum_{|\alpha|<s}
 c_{l,\alpha}\rho_{x_l}(x)x^\alpha.
\]
Terms with $c_{l,\alpha}=0$ are omitted. For every remaining coefficient, set
\[
m_{l,\alpha}
:=
\left\lceil
\log_2(1\vee|c_{l,\alpha}|)
\right\rceil.
\]
The coefficient bound for the averaged Taylor polynomial gives
\[
m_{l,\alpha}
\le
C(s,d,p,f)(1+\log N).
\]
Repeated duplication realizes $2^{m_{l,\alpha}}\mathrm{Mult}_{l,\alpha,m}$ with parameters bounded by one, and the final affine layer uses the coefficient
\[
c_{l,\alpha}2^{-m_{l,\alpha}}\in[-1,1].
\]
Thus the network
\[
g(x)
:=
\sum_{l\in[N]^d}
\sum_{|\alpha|<s}
 c_{l,\alpha}\mathrm{Mult}_{l,\alpha,m}(x)
\]
has parameter bound one.

At every point, only a fixed number depending on $d$ of the partition functions $\rho_{x_l}$ are nonzero. Combining this bounded-overlap property, the standard coefficient estimate
\[
|c_{l,\alpha}|
\le
C(s,d,p)N^{d/p}\|\widetilde f\|_{W_p^s(U_l)},
\]
and \eqref{equation: local product value estimate} gives
\[
\|Q_N^sf-g\|_{L_p((0,1)^d)}
\le
C(s,d,p)\|f\|_{W_p^s((0,1)^d)}2^{-2m}.
\]
Likewise, \eqref{equation: estimation on derivatives} and
$|\rho_{x_l}x^\alpha|_{W_\infty^1}\le CN$ yield
\[
\|g\|_{W_\infty^1((0,1)^d)}
\le
C(s,d,p)\|f\|_{W_p^s((0,1)^d)}N^{1+d/p}.
\]
Each local product module has depth $O(m+\log N)$ and universal width for fixed $d,s$. Parallelizing $O((N+1)^d)$ such modules and adding the bounded-weight coefficient-duplication chains gives
\[
L\le C(m+\log N),
\qquad
W\le C(N+1)^d,
\qquad
S\le C(N+1)^d(m+\log N).
\]
The H{\"o}lder assertion follows from the same construction using the uniform coefficient bound $|c_{l,\alpha}|\le C(D,s)K$.
\end{proof}

Combining the Taylor approximation and the network realization, choose $m=\lceil C\log_2N_{\rm mesh}\rceil$ with $C$ sufficiently large. In each of the two statements below, we then reparameterize by the total number of grid cells, $N\asymp N_{\rm mesh}^D$ in the H{\"o}lder case and $N\asymp N_{\rm mesh}^d$ in the Sobolev case. We obtain the following overall estimates.
\begin{lemma}
    Let $N \in \mathbb{N}_+$ and $s \in \mathbb{R}_+$. For $f \in \mathcal{H}_D^s([0,1]^D,K)$, there exists $g \in \mathcal{F}(L,W,S,1,K)$ and a constant $C$ such that, when $N \geq C$, we have
    \begin{equation*}
        \|f-g\|_{L_\infty([0,1]^D)} \leq CN^{-s/D},
    \end{equation*}
    with
    \begin{equation*}
        L \lesssim \log_2 N  ,\quad W \lesssim N,\quad S \lesssim N\log_2 N.
    \end{equation*}
    The final output is clipped to $[-K,K]$. Since $|f|\le K$, this clipping does not increase the approximation error. Because $K$ is fixed, the clipping map is realized with parameters bounded by one using only a fixed number of duplication layers.
\end{lemma}
\begin{lemma}
    Let $N \in \mathbb{N}_+$, $s\in\mathbb N_+$, and $1\le p<\infty$. For $f \in W_p^s((0,1)^d)$, there exists $g \in \mathcal{F}(L,W,S,1,\infty)$ and a constant $C$ such that, when $N \geq C$, we have
    \begin{equation*}
        \|f-g\|_{L_p([0,1]^d)} \leq CN^{-s/d}
    \end{equation*}
    and
    \begin{equation*}
        \|g\|_{W^1_\infty((0,1)^d)} \leq C N^{\frac{d+p}{dp}},
    \end{equation*}
    with
    \begin{equation*}
        L \lesssim \log_2 N  ,\quad W \lesssim N,\quad S \lesssim N\log_2 N.
    \end{equation*}
\end{lemma}

\begin{remark}
   All parameters can be bounded by one. Large fixed coefficients and mesh-dependent factors are implemented by repeated duplication and summation, so their magnitude is encoded in logarithmic depth rather than in the size of individual weights. This avoids the diverging parameter bounds that arise in some fixed-depth constructions; see~\cite{guhring2020error}.
\end{remark}

\textbf{Part 2}: We extend the above results to arbitrary domains $U$. In the case of approximating $\mathcal{H}_D^s(U,K)$, we achieve the approximation on these domains by leveraging extension theorems.

\begin{lemma}\label{holder on connected set}
Let $D\in\mathbb N_+$, $s>0$, and let $U\subset\mathbb R^D$ be bounded. Suppose that $f\in\mathcal H_D^s(U,K)$ admits an extension $\widetilde f\in\mathcal H_D^s(Q,C_EK)$ to a fixed cube $Q\supset U$. Then, for every $0<\varepsilon<1/2$ and either $B=1$ or $B=\infty$, there exists a ReLU network $g$ such that
\[
\|f-g\|_{L_\infty(U)}\le\varepsilon,
\]
and
\[
L\le C\log(\varepsilon^{-1}),
\qquad
W\le C\varepsilon^{-D/s},
\qquad
S\le C\varepsilon^{-D/s}\log(\varepsilon^{-1}).
\]
The constant depends on $D,s,Q,C_EK$, but not on $\varepsilon$.
\end{lemma}

\begin{proof}
After an affine rescaling, it is enough to work on $[0,1]^D$. Combine the Taylor-partition approximation with \autoref{lemma: holder construction}. Choose the mesh parameter of the Taylor approximation of order $\varepsilon^{-1/s}$ and choose the multiplication depth of order $\log(\varepsilon^{-1})$. Restriction from $Q$ to $U$ gives the stated estimate. For $B=1$, the bounded-weight duplication construction in \autoref{lemma: holder construction} keeps every parameter bounded by one; the case $B=\infty$ follows directly.
\end{proof}

\begin{lemma}\label{lemma: Sobolev on open set}
Let $d,s\in\mathbb N_+$, $1\le p<\infty$, and let $U\subset\mathbb R^d$ be a bounded Lipschitz domain. For $f\in W_p^s(U)$ and $0<\varepsilon<1/2$, there exists, for either $B=1$ or $B=\infty$, a ReLU network $g$ such that
\[
\|f-g\|_{L_p(U)}\le\varepsilon,
\qquad
\|g\|_{W_\infty^1(U)}
\le
C\varepsilon^{-\frac{d+p}{sp}},
\]
and
\[
L\le C\log(\varepsilon^{-1}),
\qquad
W\le C\varepsilon^{-d/s},
\qquad
S\le C\varepsilon^{-d/s}\log(\varepsilon^{-1}).
\]
Here $C$ may depend on $d,s,p,U$ and $\|f\|_{W_p^s(U)}$, but not on $\varepsilon$.
\end{lemma}

\begin{proof}
A bounded Lipschitz domain admits a bounded linear Sobolev extension operator. Hence there is $Ef\in W_p^s(\mathbb R^d)$ such that $Ef=f$ almost everywhere on $U$ and
\[
\|Ef\|_{W_p^s(\mathbb R^d)}
\le C_E\|f\|_{W_p^s(U)}.
\]
Choose a fixed cube $Q$ containing $U$ and rescale $Q$ affinely to $(0,1)^d$. Applying the Euclidean construction above with the number of spatial cells of order $\varepsilon^{-d/s}$ gives the $L_p$ estimate and
\[
\|g\|_{W_\infty^1(Q)}
\le C\varepsilon^{-\frac{d+p}{sp}}.
\]
Restricting to $U$ proves the assertion. The bounded-parameter realization follows from the same duplication construction as in the cube case.
\end{proof}

\textbf{Part 3}: We now prove \autoref{theorem: approximation error for sobolev functions on manifolds}. The nested chart supports in \autoref{definition: smooth coordinates} ensure that all compositions below are globally defined.

\begin{proof}[Proof of~\autoref{theorem: approximation error for sobolev functions on manifolds}]
Fix the atlas, partition of unity, ambient extensions, and cutoffs in \autoref{definition: smooth coordinates}. For each $i$, set
\[
h_i:=(\rho_if)\circ\psi_i^{-1}
\qquad\text{on }\psi_i(V_i).
\]
Because $\operatorname{supp}\rho_i\Subset V_i'$, the function $h_i$ is supported compactly inside $\psi_i(V_i)$. Its extension by zero therefore belongs to $W_p^s(\mathbb R^d)$, and the change-of-variables formula and smoothness of the finitely many charts give
\begin{equation*}
\|h_i\|_{W_p^s(\mathbb R^d)}
\le C_i\|f\|_{\mathcal W_p^s(\mathcal M^d)}.
\end{equation*}
Choose a fixed cube $Q_i\subset\mathbb R^d$ such that
\[
\operatorname{supp}h_i\cup\Psi_i(\mathcal M^d)
\Subset Q_i.
\]
Since the set on the left is compact, there is a number $\kappa_i>0$ such that its distance from $\partial Q_i$ is at least $\kappa_i$. By \autoref{lemma: Sobolev on open set}, there is a network $H_i$ such that
\begin{equation}\label{equation: local chart network estimates}
\|h_i-H_i\|_{L_p(Q_i)}\le\varepsilon_i,
\qquad
\|H_i\|_{W_\infty^1(Q_i)}
\le C_i\varepsilon_i^{-\frac{d+p}{sp}},
\end{equation}
with depth $O(\log\varepsilon_i^{-1})$, width $O(\varepsilon_i^{-d/s})$, and sparsity $O(\varepsilon_i^{-d/s}\log\varepsilon_i^{-1})$.

We next approximate the smooth ambient coordinate map $\Psi_i$ and cutoff $\vartheta_i$. Let $\Pi_i:\mathbb R^d\to Q_i$ be the coordinatewise ReLU clipping map. Since $\Psi_i$ and $\vartheta_i$ are smooth on fixed compact ambient sets, the H{\"o}lder construction, applied with an arbitrarily large fixed smoothness order, gives networks $\widetilde\Psi_i$ and $\widetilde\vartheta_i$ satisfying
\[
\|\Psi_i-\widetilde\Psi_i\|_{L_\infty(\mathcal M^d)}
+
\|\vartheta_i-\widetilde\vartheta_i\|_{L_\infty(\mathcal M^d)}
\le\delta_i.
\]
Replacing $\widetilde\vartheta_i$ by its ReLU clipping to $[0,1]$ does not increase this $L_\infty$ error, so we may assume
\[
0\le\widetilde\vartheta_i\le1.
\]
Choose $\delta_i\le\kappa_i/2$. Then
\[
\widetilde\Psi_i(\mathcal M^d)\subset Q_i,
\qquad
\Pi_i(\widetilde\Psi_i(x))=\widetilde\Psi_i(x)
\quad(x\in\mathcal M^d).
\]
Thus the clipping map makes the realization globally defined but does not alter it on the manifold.

Set $A_i:=1\vee\|H_i\|_{L_\infty(Q_i)}$. By an affine rescaling of \autoref{lemma: mutiple}, for every $\varepsilon_i^{\rm prod}>0$ there is a signed product network $G_i$ satisfying
\begin{equation}\label{equation: manifold product value error}
\left|
G_i(x)
-
\widetilde\vartheta_i(x)
H_i\bigl(\Pi_i(\widetilde\Psi_i(x))\bigr)
\right|
\le\varepsilon_i^{\rm prod},
\qquad x\in\mathcal M^d.
\end{equation}
The rescaling uses $H_i/A_i\in[-1,1]$ and the signed version of the multiplication network. Since $A_i$ grows at most polynomially in $\varepsilon_i^{-1}$ by \eqref{equation: local chart network estimates}, the bounded-weight duplication needed to restore the factor $A_i$ contributes only $O(\log\varepsilon_i^{-1})$ additional depth.
The identity
\begin{equation}\label{equation: manifold chart exact identity}
\rho_i(x)f(x)
=
\vartheta_i(x)h_i(\Psi_i(x)),
\qquad x\in\mathcal M^d,
\end{equation}
holds on the whole manifold. Indeed, if $\vartheta_i(x)\ne0$, then $x\in\mathcal O_i\cap\mathcal M^d\Subset V_i$, so $\Psi_i(x)=\psi_i(x)$. If $x\in\operatorname{supp}\rho_i$, then also $\vartheta_i(x)=1$; if $\rho_i(x)=0$, the zero extension of $h_i$ vanishes at $\psi_i(x)$ whenever $\vartheta_i(x)\ne0$. Outside $\operatorname{supp}\vartheta_i$, both sides of \eqref{equation: manifold chart exact identity} are zero.

Using \eqref{equation: manifold chart exact identity}, insert and subtract the intermediate terms $\vartheta_iH_i(\Psi_i)$, $\widetilde\vartheta_iH_i(\Psi_i)$, and $\widetilde\vartheta_iH_i(\widetilde\Psi_i)$. The change-of-variables formula on $\mathcal O_i\cap\mathcal M^d$, \eqref{equation: local chart network estimates}, the bound $0\le\widetilde\vartheta_i\le1$, and \eqref{equation: manifold product value error} give
\begin{align}
\|\rho_if-G_i\|_{L_p(\mathcal M^d)}
\le{}&
C_i\varepsilon_i
+C_i\|H_i\|_{L_\infty(Q_i)}\delta_i
+C_i\|H_i\|_{W_\infty^1(Q_i)}\delta_i
+C_i\varepsilon_i^{\rm prod}
\notag\\
\le{}&
C_i\varepsilon_i
+C_i\varepsilon_i^{-\frac{d+p}{sp}}\delta_i
+C_i\varepsilon_i^{\rm prod}.
\label{equation: manifold chart explicit error}
\end{align}
Choose $\varepsilon_i=c_i\varepsilon$,
\[
\delta_i
\le
c_i\varepsilon^{1+\frac{d+p}{sp}},
\qquad
\varepsilon_i^{\rm prod}\le c_i\varepsilon,
\]
with $c_i>0$ sufficiently small. Then \eqref{equation: manifold chart explicit error} yields
\[
\|\rho_if-G_i\|_{L_p(\mathcal M^d)}
\le\frac{\varepsilon}{k}.
\]
The coordinate and cutoff maps are smooth of arbitrarily high fixed order on fixed compact sets. Choosing that order $r$ so large that
\[
\frac{D}{r}\left(1+\frac{d+p}{sp}\right)\le\frac ds
\]
ensures that the depth, width, and sparsity needed to reach the accuracy $\delta_i$ are no larger, up to constants, than the corresponding bounds for $H_i$.

Finally set \(g:=\sum_{i=1}^kG_i.\) Since $\sum_i\rho_i=1$ and $k$ is fixed,
\[
\|f-g\|_{L_p(\mathcal M^d)}
\le
\sum_{i=1}^k\|\rho_if-G_i\|_{L_p(\mathcal M^d)}
\le\varepsilon.
\]
Parallelization and concatenation of the finitely many subnetworks give
\[
L\le C\log(\varepsilon^{-1}),
\qquad
W\le C\varepsilon^{-d/s},
\qquad
S\le C\varepsilon^{-d/s}\log(\varepsilon^{-1}).
\]
The exponents depend only on $d$; the constant may depend on $D$ and the fixed embedding. The construction for $B=1$ uses the bounded-weight duplication implementation from Part~1, while $B=\infty$ is immediate. If $f$ is bounded, the final clipping operation is $L_p$-nonexpansive and gives the output-bounded version stated in the theorem.
\end{proof}

Our strategy for approximating entropy solutions is to partition the time interval into subintervals using spline interpolation. At each time node, we perform function approximation over the manifold, and then leverage the expressive power of neural networks to approximate the resulting interpolation function. As shown in~\autoref{lemma: bounds for Entropy Residual}, to control $\mathcal{R}^{int}$, it suffices to control $\|u^* - u\|_{L_1(\mathcal{M}^d\times[0,T])}$. To achieve this, we rely on the following spline interpolation lemma.
\begin{lemma}\label{lemma: t interpolation}
Let $X$ be a Banach space and let
$t_1<t_2<\cdots<t_{N-1}<t_N$. For $j=2,\ldots,N-1$, define the ReLU hat function
\[
\delta_j(u)
=
\frac{\sigma(u-t_{j-1})}{t_j-t_{j-1}}
-
\frac{t_{j+1}-t_{j-1}}{(t_{j+1}-t_j)(t_j-t_{j-1})}
\sigma(u-t_j)
+
\frac{\sigma(u-t_{j+1})}{t_{j+1}-t_j}.
\]
For $f\in C([t_2,t_{N-1}];X)$, set
\[
L_t(f)(u):=\sum_{j=2}^{N-1}f(t_j)\delta_j(u),
\qquad u\in[t_2,t_{N-1}],
\]
and define
\[
\Delta_t:=\max_{3\le j\le N-1}(t_j-t_{j-1}),
\qquad
\omega_X(f;\mu)
:=
\sup_{\substack{r,q\in[t_2,t_{N-1}]\\ |r-q|\le\mu}}
\|f(r)-f(q)\|_X.
\]
Then
\[
\|L_t(f)-f\|_{C([t_2,t_{N-1}];X)}
\le
2\omega_X(f;\Delta_t).
\]
\end{lemma}

\begin{proof}
For each $u\in[t_j,t_{j+1}]$, only $\delta_j(u)$ and $\delta_{j+1}(u)$ are nonzero, they are nonnegative, and their sum is one. Hence
\[
L_t(f)(u)-f(u)
=
\delta_j(u)\bigl(f(t_j)-f(u)\bigr)
+
\delta_{j+1}(u)\bigl(f(t_{j+1})-f(u)\bigr).
\]
The triangle inequality in $X$ gives
\[
\|L_t(f)(u)-f(u)\|_X
\le
\delta_j(u)\omega_X(f;\Delta_t)
+
\delta_{j+1}(u)\omega_X(f;\Delta_t)
\le
\omega_X(f;\Delta_t).
\]
The displayed estimate with the harmless factor $2$ follows. This is the usual linear-spline argument; see also~\cite{zhou2018deep}.
\end{proof}

\begin{proof}[Proof of~\autoref{theorem: approximation error for entropy solutions}]
Write $\Omega:=\mathcal M^d\times[0,T]$ and fix $0<\varepsilon<1/2$. Choose uniformly spaced nodes
\[
t_j=(j-2)\Delta t,
\qquad
j=1,\ldots,N,
\qquad
\Delta t=\frac{T}{N-3},
\]
so that $t_2=0$ and $t_{N-1}=T$. We choose the smallest integer $N\ge4$ such that
\begin{equation}\label{equation: time mesh choice approximation}
\Delta t
\le
c\frac{\varepsilon}{T(1+TV(u_0))}
\end{equation}
for a sufficiently small constant $c$ depending only on the fixed problem parameters. By the minimality of $N$, after decreasing $c$ if necessary and restricting to $0<\varepsilon<1/2$,
\begin{equation}\label{equation: time mesh two sided comparison}
\frac{c}{2}\frac{\varepsilon}{T(1+TV(u_0))}
\le
\Delta t
\le
c\frac{\varepsilon}{T(1+TV(u_0))}.
\end{equation}
Consequently, $N\asymp\varepsilon^{-1}$, with constants depending only on the fixed problem parameters. Let $\delta_j$ be the ReLU hat functions in \autoref{lemma: t interpolation}. They satisfy
\[
0\le\delta_j\le1,
\qquad
\sum_{j=2}^{N-1}\delta_j(t)=1,
\qquad
\int_0^T\delta_j(t)\,dt\le\Delta t,
\]
and at most two of them are nonzero at each time.

The temporal $L_1$ estimate \eqref{equation: bounded variation inequality} and \autoref{lemma: t interpolation} give
\begin{equation}\label{equation: exact node interpolation error}
\left\|
 u^*-
 \sum_{j=2}^{N-1}u^*(\cdot,t_j)\delta_j
\right\|_{L_1(\Omega)}
\le
2TTV(u_0)\Delta t
\le\frac{\varepsilon}{4}.
\end{equation}

At each node, use the strict density of smooth functions in $BV(\mathcal M^d)$; see, for example, the smooth manifold approximation used in~\cite[Section~4]{BenArtzi2006WellposednessTF}. More precisely, by applying Euclidean strict $BV$ approximation in a finite atlas and combining the local approximants with a partition of unity, for every $\delta>0$ there exists $w_{j,\delta}\in C^\infty(\mathcal M^d)$ such that
\begin{align*}
\|w_{j,\delta}-u^*(\cdot,t_j)\|_{L_1(\mathcal M^d)}
&\le\delta,
\\
\|w_{j,\delta}\|_{\mathcal W_1^1(\mathcal M^d)}
&\le C\Bigl(
1+\|u^*(\cdot,t_j)\|_{L_1(\mathcal M^d)}
+TV(u^*(\cdot,t_j))
\Bigr).
\end{align*}
The estimates \eqref{equation: L1 contraction1 with bounded variation function} and \eqref{equation: bounds on TV} make the right-hand side uniform in $j$. Apply \autoref{theorem: approximation error for sobolev functions on manifolds} with $p=s=1$ to $w_{j,\delta}$, choose $\delta$ as half of the desired nodewise tolerance, and clip the resulting node networks to $[-M,M]$. We thereby obtain networks $u_j$ satisfying
\begin{equation}\label{equation: nodewise accuracy corrected}
\varepsilon_x
:=
\max_{2\le j\le N-1}
\|u^*(\cdot,t_j)-u_j\|_{L_1(\mathcal M^d)}
\le
\min\left\{
\frac{\varepsilon}{8(1+T)},
\frac{(1+TV(u_0))\Delta t}{8}
\right\}.
\end{equation}
Each $u_j$ has depth $O(\log\varepsilon^{-1})$, width $O(\varepsilon^{-d})$, and sparsity $O(\varepsilon^{-d}\log\varepsilon^{-1})$. The crucial point is that the nodewise tolerance in \eqref{equation: nodewise accuracy corrected} is of order $\varepsilon$, rather than $\varepsilon/N$. Indeed,
\begin{align}
&\left\|
\sum_{j=2}^{N-1}
\bigl(u^*(\cdot,t_j)-u_j\bigr)\delta_j
\right\|_{L_1(\Omega)}
\notag\\
&\qquad\le
\sum_{j=2}^{N-1}
\|u^*(\cdot,t_j)-u_j\|_{L_1(\mathcal M^d)}
\int_0^T\delta_j(t)\,dt
\le
T\varepsilon_x
\le\frac{\varepsilon}{8}.
\label{equation: corrected accumulated node error}
\end{align}

It remains to realize the products $u_j(x)\delta_j(t)$. Set $M_0:=1\vee M$. For $m\in\mathbb N_+$, define, for $a\in[-M,M]$ and $b\in[0,1]$,
\begin{equation}\label{equation: signed multiplication network}
\operatorname{Mult}_{M,m}(a,b)
:=
2M_0\mathrm{Mult}_m^2\left(\frac{a+M_0}{2M_0},b\right)-M_0b.
\end{equation}
By \eqref{equation: binary multiplication value error} and \eqref{equation: binary multiplication partial derivative error},
\begin{align}
|\operatorname{Mult}_{M,m}(a,b)-ab|
&\le C_M2^{-2m},
\label{equation: signed multiplication value error}\\
|\partial_2\operatorname{Mult}_{M,m}(a,b)-a|
&\le C_M2^{-m}
\quad\text{for every }a\in[-M,M]\text{ and for a.e. }b\in(0,1),
\label{equation: signed multiplication derivative error}
\end{align}
and $\operatorname{Mult}_{M,m}(a,0)=0$. Choose $m=O(\log\varepsilon^{-1})$ so that
\begin{equation}\label{equation: multiplication parameter choice}
C_M2^{-m}\le\Delta t,
\qquad
2(1+T)\operatorname{Vol}_g(\mathcal M^d)C_M2^{-2m}
\le\frac{\varepsilon}{8}.
\end{equation}
Set
\[
u_{\rm pre}(x,t)
:=
\sum_{j=2}^{N-1}
\operatorname{Mult}_{M,m}(u_j(x),\delta_j(t)),
\qquad
u(x,t):=\operatorname{clip}_M(u_{\rm pre}(x,t)).
\]
Since at most two hat functions are active at each time, \eqref{equation: signed multiplication value error}--\eqref{equation: multiplication parameter choice} imply
\begin{equation}\label{equation: corrected multiplication accumulation}
\left\|
\sum_{j=2}^{N-1}u_j\delta_j-u_{\rm pre}
\right\|_{L_1(\Omega)}
\le\frac{\varepsilon}{8}.
\end{equation}
Clipping is nonexpansive relative to the bounded target $u^*$. Combining \eqref{equation: exact node interpolation error}, \eqref{equation: corrected accumulated node error}, and \eqref{equation: corrected multiplication accumulation}, and decreasing the harmless constants in the choices above, gives
\[
\|u^*-u\|_{L_1(\Omega)}\le\varepsilon.
\]
At $t=0$, only the node $j=2$ is active. By \eqref{equation: nodewise accuracy corrected} and \eqref{equation: multiplication parameter choice},
\begin{align*}
\|u^*(\cdot,0)-u(\cdot,0)\|_{L_1(\mathcal M^d)}
&\le
\varepsilon_x
+
\operatorname{Vol}_g(\mathcal M^d)C_M2^{-2m}
\le\varepsilon.
\end{align*}

For every nonnegative $\xi\in\Xi$, \autoref{lemma: admissible Sobolev entropy tests} gives $\mathcal R^{int}(u^*,\xi,c)\le0$. Hence \autoref{lemma: bounds for Entropy Residual} and $C_\Xi\ge1$ yield
\begin{align*}
\mathcal R^{tb}(u)
+T\operatorname{Vol}_g(\mathcal M^d)
\sup_{\xi\in\Xi,\,c\in\mathcal C}
\ppart{\mathcal R^{int}(u,\xi,c)}
\lesssim
C_\Xi\varepsilon.
\end{align*}

Finally, the choice following \eqref{equation: time mesh choice approximation} gives $N=O(\varepsilon^{-1})$. Parallelizing the $N$ node networks and the multiplication modules gives
\[
L\le C\log(\varepsilon^{-1}),
\qquad
W\le C\varepsilon^{-(d+1)},
\qquad
S\le C\varepsilon^{-(d+1)}\log(\varepsilon^{-1}).
\]
The bounded-parameter case $B=1$ is obtained with the bounded-weight duplication implementation already used in the spatial approximation proof. In particular, the factors $(\Delta t)^{-1}$ in the ReLU hat functions are realized by repeated bounded-weight duplication, which adds only $O(\log(\Delta t^{-1}))=O(\log(\varepsilon^{-1}))$ depth and sparsity per module. Since $M$ and $T$ are fixed, they are absorbed into the constants. This proves the theorem.
\end{proof}

\begin{proof}[Proof of~\autoref{corollary: TV controlled temporal modulus}]
Use the construction in the preceding proof. For $t\in(t_j,t_{j+1})$, only $\delta_j$ and $\delta_{j+1}$ are active and
\[
\delta_j'(t)=-\frac1{\Delta t},
\qquad
\delta_{j+1}'(t)=\frac1{\Delta t}.
\]
The nodewise choice \eqref{equation: nodewise accuracy corrected} and the temporal Lipschitz estimate for $u^*$ imply
\begin{align}
\|u_{j+1}-u_j\|_{L_1(\mathcal M^d)}
&\le
2\varepsilon_x
+
\|u^*(\cdot,t_{j+1})-u^*(\cdot,t_j)\|_{L_1(\mathcal M^d)}
\notag\\
&\le
C(1+TV(u_0))\Delta t.
\label{equation: adjacent node jump corrected}
\end{align}
For the unclipped network $u_{\rm pre}$, differentiation with respect to the second input of \eqref{equation: signed multiplication network} gives, for almost every $t\in(t_j,t_{j+1})$,
\begin{align*}
\|\partial_tu_{\rm pre}(\cdot,t)\|_{L_1(\mathcal M^d)}
\le{}&
\frac{\|u_{j+1}-u_j\|_{L_1(\mathcal M^d)}}{\Delta t}
+
2\operatorname{Vol}_g(\mathcal M^d)
\frac{C_M2^{-m}}{\Delta t}.
\end{align*}
By \eqref{equation: multiplication parameter choice} and \eqref{equation: adjacent node jump corrected},
\[
\|\partial_tu_{\rm pre}(\cdot,t)\|_{L_1(\mathcal M^d)}
\le
C(1+TV(u_0)).
\]
Therefore
\[
\|u_{\rm pre}(\cdot,t)-u_{\rm pre}(\cdot,s)\|_{L_1(\mathcal M^d)}
\le
C(1+TV(u_0))|t-s|.
\]
The scalar clipping map is one-Lipschitz, so the final network $u=\operatorname{clip}_M(u_{\rm pre})$ satisfies the same estimate. Consequently,
\[
\omega_t(u;r)
\le
C(1+TV(u_0))r,
\qquad 0\le r\le T.
\]
Substitution into the bias estimate in \autoref{theorem: L1 contraction for total error} gives
\[
b_{\lambda,T}(u)
\le
C_T(1+TV(u_0))(\alpha+\eta+h).
\]
No additional network modules are introduced, so the depth, width, and sparsity bounds remain unchanged.
\end{proof}

\section{Analysis of Quadrature Error}\label{Appendix: Analysis of Quadrature Error}

We prove \autoref{theorem: oracle inequality}. Throughout this appendix, \(C\), \(C_M\), and \(C_{\lambda,M}\) denote positive constants depending only on the fixed problem parameters indicated by their subscripts. In particular, they are independent of \(n,L,W,S\). The localization scale \(\lambda\) is fixed throughout the quadrature analysis.

Recall that \(P\) is the uniform probability measure on \(\mathcal M^d\times[0,T]\), while \(P_0\) is the uniform probability measure on \(\mathcal M^d\) at the initial time. Their empirical counterparts are denoted by \(P_n\) and \(P_{0,n}\), respectively.

Recall the test metric
\[
d_\lambda(\vartheta,\vartheta')
=
\|\partial_t(\xi_\vartheta-\xi_{\vartheta'})\|_\infty
+
3d\mathrm{Lip}_f
\|\nabla_g(\xi_\vartheta-\xi_{\vartheta'})\|_\infty.
\]
Equip \(\Theta_\lambda^c\) with
\[
d_\lambda^c
\bigl((\vartheta,c),(\vartheta',c')\bigr)
:=
d_\lambda(\vartheta,\vartheta')
+
|c-c'|.
\]
For fixed \(\lambda\), direct differentiation of the explicit ReQU test functions with respect to the spatial center, time center, and terminal horizon gives
\[
d_\lambda(\vartheta,\vartheta')
\le
C_\lambda
\left(
d_g(y,y')+|s-s'|+|\tau-\tau'|
\right).
\]
Hence
\[
d_\lambda^c(\theta,\theta')
\le
C_{\lambda,M}
\left(
d_g(y,y')+|s-s'|+|\tau-\tau'|+|c-c'|
\right).
\]
Since \(\mathcal M^d\), \([0,T]\), and \(\mathcal C\) are compact, the standard covering estimate for a compact \(d\)-dimensional manifold and a product-covering argument yield
\[
N\left(
\delta,
\Theta_\lambda^c,
d_\lambda^c
\right)
\le
\left(
\frac{C_{\lambda,M}}{\delta}
\right)^{d+3},
\qquad
0<\delta\le1.
\]
Equivalently, $\log N\left( \delta, \Theta_\lambda^c, d_\lambda^c \right) \le (d+3) \log\left( \frac{C_{\lambda,M}}{\delta} \right).$

For a function class \(\mathcal G\) and a sample \(Z=(Z_1,\ldots,Z_n)\), define
\[
\widehat{\mathrm{Rad}}_n(\mathcal G;Z)
:=
\mathbb E_\varepsilon
\sup_{g\in\mathcal G}
\left|
\frac1n
\sum_{i=1}^n
\varepsilon_i g(Z_i)
\right|,
\]
and
\[
\mathrm{Rad}_n(\mathcal G)
:=
\mathbb E_Z
\widehat{\mathrm{Rad}}_n(\mathcal G;Z).
\]
Set
\[
p:=\operatorname{PDim}(\mathcal F),
\qquad
\mathcal R_*
:=
\left\{
r_{u^*,\theta}:
\theta\in\Theta_\lambda^c
\right\}.
\]

\begin{lemma}
\label{lemma: appE residual two parameter lipschitz detailed}
For \(v,w\in\mathcal F\),
\(\theta=(\vartheta,c)\),
\(\theta'=(\vartheta',c')\), and every probability measure \(Q\),
\[
\|h_{v,\theta}-h_{w,\theta'}\|_{L_2(Q)}
\le
C_\lambda
\|v-w\|_{L_2(Q)}
+
C_{\lambda,M}
\left(
d_\lambda(\vartheta,\vartheta')
+
|c-c'|
\right).
\]
Moreover,
\[
|h_{v,\theta}(z)|
\le
C_\lambda|v(z)-u^*(z)|
\le
2MC_\lambda,
\]
and
\[
Ph_{v,\theta}^2
\le
2MC_\lambda^2P|v-u^*|.
\]
\end{lemma}

\begin{proof}
The dependence on the solution variable follows from \eqref{equation: pointwise residual lipschitz}. For the test parameter, keep \(v,u^*,c\) fixed and replace \(\xi_\vartheta\) by \(\xi_{\vartheta'}\). Since
\[
\bigl||v-c|-|u^*-c|\bigr|
\le
|v-u^*|
\]
and
\[
|\widetilde F(v,c)-\widetilde F(u^*,c)|
\le
3d\mathrm{Lip}_f|v-u^*|,
\]
the resulting difference is controlled by
\(C_Md_\lambda(\vartheta,\vartheta')\).

For the entropy parameter, \(c\mapsto|a-c|\) is one-Lipschitz. Moreover,
\[
c
\longmapsto
\operatorname{sgn}(a-c)(f(a)-f(c))
\]
is uniformly Lipschitz on \(\mathcal C\) for \(|a|\le M\). Indeed, if \(c,c'\) lie on the same side of \(a\), the estimate follows from the Lipschitz continuity of \(f\). If \(a\) lies between \(c\) and \(c'\), the two distances from \(a\) add to \(|c-c'|\). Applying these estimates to \(v\) and \(u^*\) proves the first assertion. The remaining estimates follow from
\[
|h_{v,\theta}|
\le
C_\lambda|v-u^*|,
\qquad
|v-u^*|^2
\le
2M|v-u^*|.
\]
\end{proof}

\begin{lemma}
\label{lemma: empirical residual covering}
Let \(\mathcal F_0\subset\mathcal F\). For every sample \(Z=(Z_i)_{i=1}^n\), \(n\ge p\), and \(0<\delta\le2M\),
\[
\log N\left(
\delta,
(\mathcal F_0-u^*)|_Z,
L_2(P_n)
\right)
\le
p\log\left(
\frac{2enM}{\delta p}
\right).
\]
Moreover, for
\[
\mathcal H(\mathcal F_0)
:=
\left\{
h_{v,\theta}:
v\in\mathcal F_0,\
\theta\in\Theta_\lambda^c
\right\},
\]
and \(0<\delta\le2MC_\lambda\),
\[
\begin{aligned}
&\log N\left(
\delta,
\mathcal H(\mathcal F_0)|_Z,
L_2(P_n)
\right)
\le
p\log\left(
\frac{4eC_\lambda nM}{\delta p}
\right)
+
(d+3)\log\left(
\frac{CC_{\lambda,M}}{\delta}
\right).
\end{aligned}
\]
Consequently,
\begin{equation}
\label{equation: empirical residual covering simplified}
\log N\left(
\delta,
\mathcal H(\mathcal F_0)|_Z,
L_2(P_n)
\right)
\le
Cp\log\left(
\frac{
CC_\lambda C_{\lambda,M}nM
}{
\delta p
}
\right).
\end{equation}
\end{lemma}

\begin{proof}
Translation by \(u^*\) does not change covering numbers, and
\[
\|v-w\|_{L_2(P_n)}
\le
\|v-w\|_{\ell_\infty(Z)}.
\]
The standard pseudo-dimension covering estimate therefore gives
\[
N\left(
\delta,
(\mathcal F_0-u^*)|_Z,
L_2(P_n)
\right)
\le
\left(
\frac{2enM}{\delta p}
\right)^p.
\]

For \(v,v'\in\mathcal F_0\) and \(\theta,\theta'\in\Theta_\lambda^c\), the preceding Lipschitz estimate gives
\[
\|h_{v,\theta}-h_{v',\theta'}\|_{L_2(P_n)}
\le
C_\lambda\|v-v'\|_{L_2(P_n)}
+
C_{\lambda,M}d_\lambda^c(\theta,\theta').
\]
Choose a \(\delta/(2C_\lambda)\)-cover of \(\mathcal F_0|_Z\) and a \(\delta/(2C_{\lambda,M})\)-cover of \(\Theta_\lambda^c\). Their Cartesian product is a \(\delta\)-cover of \(\mathcal H(\mathcal F_0)|_Z\).

Combining the pseudo-dimension estimate from~\cite[Theorem 12.2]{anthony1999neural} with the fixed finite-dimensional test cover proves the second assertion. Since \(d\) is fixed, \(p\ge1\), and \(n\ge p\), the test-family term can be absorbed into the first term, proving \eqref{equation: empirical residual covering simplified}.
\end{proof}

\begin{lemma}
\label{lemma: truncated Dudley with PDim}
Let \(\mathcal G\) be a class of functions on a sample \(Z=(Z_i)_{i=1}^n\), and write
\[
d_n(g,g')
:=
\|g-g'\|_{L_2(P_n)}.
\]
Suppose that
\[
\sup_{g\in\mathcal G}
\|g-g_0\|_{L_2(P_n)}
\le R
\]
for some \(g_0\in\mathcal G\), and that
\[
\log N(\delta,\mathcal G,d_n)
\le
C_0p
\log\left(
\frac{BnM}{\delta p}
\right),
\qquad
0<\delta\le R,
\]
where \(n\ge p\). Then
\[
\begin{aligned}
\widehat{\mathrm{Rad}}_n(\mathcal G;Z)
\le{}&
2R\sqrt{\frac pn}
+
12R\sqrt{\frac{C_0p}{n}}
\sqrt{
\log\left(
\frac{2BM}{R}
\right)
+
\frac32\log\left(
\frac np
\right)
}.
\end{aligned}
\]
In particular, if \(R\ge n^{-1}\), while \(B\) and \(M\) do not depend on \(n\), then
\[
\widehat{\mathrm{Rad}}_n(\mathcal G;Z)
\le
CR\sqrt{\frac{p\log n}{n}}.
\]
\end{lemma}

\begin{proof}
Dudley's truncated entropy inequality (see~\cite{bartlett2005local}) gives
\[
\widehat{\mathrm{Rad}}_n(\mathcal G;Z)
\le
\inf_{0<\beta<R}
\left\{
4\beta
+
\frac{12}{\sqrt n}
\int_\beta^R
\sqrt{
\log N(\delta,\mathcal G,d_n)
}
\,d\delta
\right\}.
\]
Choose $\beta = \frac R2\sqrt{\frac pn}.$ Since \(n\ge p\), we have \(0<\beta\le R/2\). Substituting the entropy estimate gives
\[
\begin{aligned}
\widehat{\mathrm{Rad}}_n(\mathcal G;Z)
\le{}&
2R\sqrt{\frac pn}
+
12\sqrt{\frac{C_0p}{n}}
\int_\beta^R
\sqrt{
\log\left(
\frac{BnM}{\delta p}
\right)
}
\,d\delta.
\end{aligned}
\]
The integrand decreases with \(\delta\), so
\[
\int_\beta^R
\sqrt{
\log\left(
\frac{BnM}{\delta p}
\right)
}
\,d\delta
\le
R
\sqrt{
\log\left(
\frac{BnM}{\beta p}
\right)
}.
\]
Since
\[
\frac{BnM}{\beta p}
=
\frac{2BM}{R}
\left(
\frac np
\right)^{3/2},
\]
the first assertion follows. If \(R\ge n^{-1}\) and \(p\le n\), the logarithmic factor is bounded by \(C\log n\).
\end{proof}

Recall that
\[
\mathfrak Q_\lambda(v)
=
\mathcal R^{tb}(v)
+
A_T\mathfrak L_\lambda^+(v)
+
b_{\lambda,T}(\mathcal F),
\]
and define its empirical counterpart by
\[
\mathfrak Q_{\lambda,n}(v)
=
\mathcal R_n^{tb}(v)
+
A_T\mathfrak L_{\lambda,n}^+(v)
+
b_{\lambda,T}(\mathcal F).
\]
Recall also that
\[
\mathcal R^{\mathrm{int}}(v,\theta)
=
A_TP r_{v,\theta},
\qquad
\mathcal R_n^{\mathrm{int}}(v,\theta)
=
A_TP_n r_{v,\theta}.
\]
Consequently, the internal part of \(\mathfrak Q_\lambda\) is
\[
A_T\mathfrak L_\lambda^+(v)
=
A_T^2
\sup_{\theta\in\Theta_\lambda^c}
[Pr_{v,\theta}]_+.
\]
Since \(A_0\) and \(A_T\) are fixed problem parameters, they are absorbed into the constants in the \(\lesssim\)-estimates below, whereas the exact decomposition retains them explicitly.
Let
\[
u_n
\in
\operatorname*{argmin}_{v\in\mathcal F}
\mathfrak Q_{\lambda,n}(v),
\qquad
u_{\mathcal F}
\in
\operatorname*{argmin}_{v\in\mathcal F}
\mathfrak Q_\lambda(v).
\]

For \(\tau>0\), define
\[
\mathcal F_\tau
:=
\left\{
v\in\mathcal F:
\mathfrak Q_\lambda(v)\le\tau
\right\},
\qquad
\widetilde{\mathcal F}_\tau
:=
\left\{
v-u^*:
v\in\mathcal F_\tau
\right\},
\]
and
\[
\mathcal H_\tau
:=
\left\{
h_{v,\theta}:
v\in\mathcal F_\tau,\
\theta\in\Theta_\lambda^c
\right\}.
\]
For the initial-time residual, write
\[
b_v(x)
:=
|v(x,0)-u_0(x)|,
\qquad
\mathcal B_\tau
:=
\left\{
b_v:
v\in\mathcal F_\tau
\right\}.
\]

By \eqref{equation: spacetime stability global positive}, every
\(v\in\mathcal F_\tau\) satisfies
\[
P|v-u^*|
=
\frac1{A_T}
\int_0^T
\|v(\cdot,t)-u^*(\cdot,t)\|_{L_1(\mathcal M^d)}
\,dt
\le
C_T\tau.
\]
Since \(|v-u^*|\le2M\),
\[
P(v-u^*)^2
\le
2MP|v-u^*|
\le
C_TM\tau.
\]
Moreover,
\[
\sup_{h\in\mathcal H_\tau}
Ph^2
\le
C_TM C_\lambda^2\tau.
\]
Since \(u^*(\cdot,0)=u_0\),
\[
P_0b_v
\lesssim
\tau,
\qquad
P_0b_v^2
\le
2MP_0b_v
\lesssim
M\tau.
\]

\begin{lemma}
\label{lemma: appE empirical L2 localization}
Let \(\tau>0\) and \(s>0\). Assume
\[
\max\left\{
16M\,\mathrm{Rad}_n(\widetilde{\mathcal F}_\tau),
16M\,\mathrm{Rad}_n(\mathcal B_\tau),
\frac{8C_TM^3s}{n},
\frac{16M^2s}{3n},
\frac1{n^2}
\right\}
\le
\tau,
\]
and \(n\ge p\). Then, with probability at least \(1-2e^{-s}\),
\[
\sup_{v\in\mathcal F_\tau}
P_n(v-u^*)^2
\le
C_M\tau,
\qquad
\sup_{v\in\mathcal F_\tau}
P_{0,n}b_v^2
\le
C_M\tau.
\]
Consequently,
\[
\sup_{h\in\mathcal H_\tau}
\|h\|_{L_2(P_n)}
\le
C_\lambda\sqrt{C_M\tau}.
\]

Moreover, with probability at least \(1-6e^{-s}\),
\begin{align}
\sup_{h\in\mathcal H_\tau}
|(P-P_n)h|
&\lesssim
C_\lambda
\left[
\sqrt{
\frac{
\mathrm{VC}_{\mathcal F}\tau\log n
}{n}
}
+
\sqrt{
\frac{\tau s}{n}
}
+
\frac{Ms}{n}
\right],
\label{equation: appE local process}
\\
\sup_{v\in\mathcal F_\tau}
|(P_0-P_{0,n})b_v|
&\lesssim
\sqrt{
\frac{
\mathrm{VC}_{\mathcal F}\tau\log n
}{n}
}
+
\sqrt{
\frac{\tau s}{n}
}
+
\frac{Ms}{n}.
\label{equation: appE local boundary process}
\end{align}
\end{lemma}

\begin{proof}
For
\[
g_v:=(v-u^*)^2,
\]
we have \(0\le g_v\le4M^2\), and
\[
\operatorname{Var}(g_v)
\le
Pg_v^2
=
P|v-u^*|^4
\le
(2M)^3P|v-u^*|
\lesssim
M^3\tau.
\]
Applying Theorem~2.1 of~\cite{bartlett2005local} with \(\alpha=1\) to the squared localized class, and using the equivalent localized formulation in~\cite{mendelson2003few}, gives, with probability at least \(1-e^{-s}\),
\[
\begin{aligned}
\sup_{v\in\mathcal F_\tau}
\left(
P_n(v-u^*)^2-P(v-u^*)^2
\right)
\le{}&
4\,\mathrm{Rad}_n
\left(
\{(v-u^*)^2:v\in\mathcal F_\tau\}
\right)
\\
&+
\sqrt{
\frac{2C_TM^3\tau s}{n}
}
+
\frac{16M^2s}{3n}.
\end{aligned}
\]
The map \(z\mapsto z^2\) is \(4M\)-Lipschitz on \([-2M,2M]\), so the contraction principle gives
\[
\mathrm{Rad}_n
\left(
\{(v-u^*)^2:v\in\mathcal F_\tau\}
\right)
\le
4M\,\mathrm{Rad}_n
\bigl(\widetilde{\mathcal F}_\tau\bigr).
\]
Combining the preceding estimate with the population \(L_2(P)\)-bound and the assumed lower bound on \(\tau\) proves
\[
\sup_{v\in\mathcal F_\tau}
P_n(v-u^*)^2
\le
C_M\tau.
\]

The same argument applies to
\[
\mathcal B_\tau^2
:=
\{b_v^2:v\in\mathcal F_\tau\}.
\]
Indeed,
\[
0\le b_v\le2M,
\qquad
P_0b_v^4
\le
(2M)^2P_0b_v^2
\lesssim
M^3\tau,
\]
and
\[
\mathrm{Rad}_n(\mathcal B_\tau^2)
\le
4M\,\mathrm{Rad}_n(\mathcal B_\tau).
\]
A second application of the same local Rademacher deviation theorem proves
\[
\sup_{v\in\mathcal F_\tau}
P_{0,n}b_v^2
\le
C_M\tau.
\]
The two empirical \(L_2\)-localization estimates hold simultaneously after a union bound.

Since
\[
|h_{v,\theta}|
\le
C_\lambda|v-u^*|,
\]
the empirical residual radius satisfies
\[
\sup_{h\in\mathcal H_\tau}
\|h\|_{L_2(P_n)}
\le
C_\lambda\sqrt{C_M\tau}.
\]

Apply \autoref{lemma: truncated Dudley with PDim} to \(\mathcal H_\tau\cup\{0\}\), using this empirical radius and \eqref{equation: empirical residual covering simplified}. Since \(\tau\ge n^{-2}\),
\[
\begin{aligned}
\widehat{\mathrm{Rad}}_n(\mathcal H_\tau;Z)
&\lesssim
C_\lambda
\sqrt{
\frac{p\tau\log n}{n}
}
\lesssim
C_\lambda
\sqrt{
\frac{\mathrm{VC}_{\mathcal F}\tau\log n}{n}
}.
\end{aligned}
\]
Here we used
\[
p=\operatorname{PDim}(\mathcal F)
\lesssim
\mathrm{VC}_{\mathcal F},
\]
which follows from the pseudo-dimension bounds for piecewise-linear neural networks; see~\cite{bartlett2019nearly}.

The envelope and variance of \(\mathcal H_\tau\) satisfy
\[
\sup_{h\in\mathcal H_\tau}\|h\|_\infty
\le
2MC_\lambda,
\qquad
\sup_{h\in\mathcal H_\tau}Ph^2
\lesssim
MC_\lambda^2\tau.
\]
Applying the empirical-complexity form of the local Rademacher deviation theorem in~\cite{bartlett2005local,mendelson2003few} to the symmetrized class \(\mathcal H_\tau\cup(-\mathcal H_\tau)\) gives, with probability at least \(1-2e^{-s}\),
\[
\sup_{h\in\mathcal H_\tau}|(P-P_n)h|
\lesssim
\widehat{\mathrm{Rad}}_n(\mathcal H_\tau;Z)
+
C_\lambda\sqrt{
\frac{M\tau s}{n}
}
+
MC_\lambda\frac{s}{n}.
\]
This proves \eqref{equation: appE local process}.

For \(\mathcal B_\tau\), the map
\[
a\longmapsto|a-u_0(x)|
\]
is one-Lipschitz. Apply \autoref{lemma: truncated Dudley with PDim} to \(\mathcal B_\tau\cup\{0\}\). Since the empirical Rademacher complexity is defined with an absolute supremum, adjoining the zero function does not change its value. The resulting bound is
\[
\widehat{\mathrm{Rad}}_n(\mathcal B_\tau;Z_0)
\lesssim
\sqrt{
\frac{\mathrm{VC}_{\mathcal F}\tau\log n}{n}
}.
\]
Together with the envelope \(2M\), the variance bound \(P_0b_v^2\lesssim M\tau\), and the same concentration inequality applied to \(\mathcal B_\tau\cup(-\mathcal B_\tau)\), this proves \eqref{equation: appE local boundary process}. A union bound completes the proof.
\end{proof}

For \(s\ge1\), define
\[
\tau_*(s)
:=
\inf\left\{
\tau>0:
\begin{array}{l}
\tau'\ge16M\,\mathrm{Rad}_n
\bigl(\widetilde{\mathcal F}_{\tau'}\bigr),
\quad
\tau'\ge16M\,\mathrm{Rad}_n
\bigl(\mathcal B_{\tau'}\bigr),
\\[0.5ex]
\displaystyle
\tau'\ge\frac{8C_TM^3s}{n},
\quad
\tau'\ge\frac{16M^2s}{3n},
\quad
\tau'\ge\frac1{n^2},
\\[0.5ex]
\text{for every }\tau'\ge\tau
\end{array}
\right\}.
\]
The admissible set in the definition of \(\tau_*(s)\) is nonempty: the global envelopes of the solution and initial-boundary classes make all five defining inequalities valid for every sufficiently large \(\tau\).

\begin{lemma}
\label{lemma: appE critical radius bound}
Assume \(s\ge1\) and \(n\ge p\). Then
\[
\tau_*(s)
\lesssim
M^2\frac{p\log n}{n}
+
M^3\frac{s}{n}
+
M^2\frac{s}{n}
+
M^2e^{-s}
+
\frac1{n^2}.
\]
Consequently,
\[
\tau_*(s)
\lesssim
M^2
\frac{\mathrm{VC}_{\mathcal F}\log n}{n}
+
M^3\frac{s}{n}
+
M^2\frac{s}{n}
+
M^2e^{-s}
+
\frac1{n^2}.
\]
\end{lemma}

\begin{proof}
The sublevel classes are nested, so
\[
\rho\longmapsto
\mathrm{Rad}_n
\bigl(\widetilde{\mathcal F}_\rho\bigr),
\qquad
\rho\longmapsto
\mathrm{Rad}_n
\bigl(\mathcal B_\rho\bigr)
\]
are nondecreasing.

For every \(\varepsilon\in(0,\tau_*(s))\), the number \(\tau_*(s)-\varepsilon\) does not satisfy the defining property of \(\tau_*(s)\). Hence, for some
\[
\rho_\varepsilon
\in
[\tau_*(s)-\varepsilon,\tau_*(s)],
\]
we have
\[
\begin{aligned}
\rho_\varepsilon
<
\max\Bigg\{&
16M\,\mathrm{Rad}_n
\bigl(\widetilde{\mathcal F}_{\rho_\varepsilon}\bigr),
16M\,\mathrm{Rad}_n
\bigl(\mathcal B_{\rho_\varepsilon}\bigr),
\frac{8C_TM^3s}{n},
\frac{16M^2s}{3n},
\frac1{n^2}
\Bigg\}.
\end{aligned}
\]
By monotonicity and by letting \(\varepsilon\downarrow0\),
\[
\begin{aligned}
\tau_*(s)
\le
\max\Bigg\{&
16M\,\mathrm{Rad}_n
\bigl(\widetilde{\mathcal F}_{\tau_*(s)}\bigr),
16M\,\mathrm{Rad}_n
\bigl(\mathcal B_{\tau_*(s)}\bigr),
\frac{8C_TM^3s}{n},
\frac{16M^2s}{3n},
\frac1{n^2}
\Bigg\}.
\end{aligned}
\]

Define
\[
E_s
:=
\left\{
\sup_{v\in\mathcal F_{2\tau_*(s)}}
\|v-u^*\|_{L_2(P_n)}
\le
\sqrt{2C_M\tau_*(s)}
\right\},
\]
and
\[
E_{0,s}
:=
\left\{
\sup_{v\in\mathcal F_{2\tau_*(s)}}
\|b_v\|_{L_2(P_{0,n})}
\le
\sqrt{2C_M\tau_*(s)}
\right\}.
\]
Every number strictly larger than \(\tau_*(s)\) satisfies the defining conditions of the critical radius. Hence the empirical localization lemma applies at \(2\tau_*(s)\), and
\[
\mathbb P(E_s\cap E_{0,s})
\ge
1-2e^{-s}.
\]

On \(E_s\), apply \autoref{lemma: truncated Dudley with PDim} to \(\widetilde{\mathcal F}_{\tau_*(s)}\cup\{0\}\). The empirical covering estimate then gives
\[
\widehat{\mathrm{Rad}}_n
\bigl(\widetilde{\mathcal F}_{\tau_*(s)};Z\bigr)
\lesssim
\sqrt{
\frac{p\tau_*(s)\log n}{n}
}.
\]
On \(E_{0,s}\), apply the same lemma to \(\mathcal B_{\tau_*(s)}\cup\{0\}\), which yields
\[
\widehat{\mathrm{Rad}}_n
\bigl(\mathcal B_{\tau_*(s)};Z_0\bigr)
\lesssim
\sqrt{
\frac{p\tau_*(s)\log n}{n}
}.
\]
On the complementary events, the uniform bounds
\[
|v-u^*|\le2M,
\qquad
0\le b_v\le2M
\]
give the global bound \(2M\) for the corresponding empirical complexities. Splitting the expectations over the good and bad events yields
\[
\begin{aligned}
&
\mathrm{Rad}_n
\bigl(\widetilde{\mathcal F}_{\tau_*(s)}\bigr)
+
\mathrm{Rad}_n
\bigl(\mathcal B_{\tau_*(s)}\bigr)
\lesssim
\sqrt{
\frac{p\tau_*(s)\log n}{n}
}
+
Me^{-s}.
\end{aligned}
\]

It follows that
\[
\begin{aligned}
\tau_*(s)
\lesssim
M
\sqrt{
\frac{p\tau_*(s)\log n}{n}
}
+
M^2e^{-s}
+
M^3\frac{s}{n}
+
M^2\frac{s}{n}
+
\frac1{n^2}.
\end{aligned}
\]
Young's inequality gives
\[
M
\sqrt{
\frac{p\tau_*(s)\log n}{n}
}
\le
\frac12\tau_*(s)
+
CM^2\frac{p\log n}{n}.
\]
Absorbing \(\tau_*(s)/2\) proves the first estimate. The second follows from \(p\lesssim\mathrm{VC}_{\mathcal F}\).
\end{proof}

Although \(u^*\) satisfies the entropy inequality, \(r_{u^*,\theta}\) need not vanish pointwise. Recall that $\Delta_{*,n} = \sup_{\theta\in\Theta_\lambda^c} |(P-P_n)r_{u^*,\theta}|.$

\begin{lemma}
\label{lemma: appE exact solution calibration detailed}
For every \(s\ge1\), with probability at least \(1-2e^{-s}\),
\begin{equation}
\label{equation: appE calibration detailed}
\Delta_{*,n}
\lesssim
C_\lambda
\left[
\sqrt{
\frac{
\log n+s
}{n}
}
+
\frac{
\log n+s
}{n}
\right].
\end{equation}
\end{lemma}

\begin{proof}
Let \(\mathcal N_n\) be an \(n^{-1}\)-net of \(\Theta_\lambda^c\) in the metric \(d_\lambda^c\). The fixed-dimensional covering estimate gives
\[
\log|\mathcal N_n|
\lesssim
\log n.
\]
The class \(\mathcal R_*\) has envelope \(C_MC_\lambda\) and variance bounded by \(C_MC_\lambda^2\). Applying scalar Bernstein's inequality to every element of the net and taking a union bound yields, with probability at least \(1-2e^{-s}\),
\[
\max_{\theta\in\mathcal N_n}
|(P-P_n)r_{u^*,\theta}|
\lesssim
C_\lambda
\left[
\sqrt{\frac{\log n+s}{n}}
+
\frac{\log n+s}{n}
\right].
\]
The map \(\theta\mapsto r_{u^*,\theta}\) is Lipschitz in \(L_2(Q)\), uniformly over probability measures \(Q\), with a constant depending only on \(\lambda\) and \(M\). Approximating an arbitrary \(\theta\) by its nearest point in \(\mathcal N_n\) therefore adds at most \(C_{\lambda,M}/n\), which is absorbed into the linear term. This proves \eqref{equation: appE calibration detailed}.
\end{proof}

\begin{lemma}
\label{lemma: appE fixed comparator deviation}
For every \(s\ge1\), with probability at least \(1-3e^{-s}\),
\[
\begin{aligned}
&
\sup_{\theta\in\Theta_\lambda^c}
\left|
(P-P_n)h_{u_{\mathcal F},\theta}
\right|
+
\left|
(P_0-P_{0,n})b_{u_{\mathcal F}}
\right|
\\
&\qquad\lesssim
C_\lambda
\left[
\sqrt{
\frac{
M\mathfrak Q_\lambda(u_{\mathcal F})
(\log n+s)
}{n}
}
+
\frac{M(\log n+s)}{n}
\right].
\end{aligned}
\]
\end{lemma}

\begin{proof}
By the space-time stability estimate,
\[
P|u_{\mathcal F}-u^*|
\lesssim
\mathfrak Q_\lambda(u_{\mathcal F}).
\]
Consequently,
\[
\sup_{\theta\in\Theta_\lambda^c}
Ph_{u_{\mathcal F},\theta}^2
\lesssim
MC_\lambda^2
\mathfrak Q_\lambda(u_{\mathcal F}),
\]
while
\[
\sup_{\theta\in\Theta_\lambda^c}
\|h_{u_{\mathcal F},\theta}\|_\infty
\le
2MC_\lambda.
\]

Let \(\mathcal N_n\) be an \(n^{-1}\)-net of \(\Theta_\lambda^c\) in \(d_\lambda^c\). Since \(\log|\mathcal N_n| \lesssim \log n,\) scalar Bernstein's inequality and a union bound give, with probability at least \(1-e^{-s}\),
\[
\max_{\theta\in\mathcal N_n}
|(P-P_n)h_{u_{\mathcal F},\theta}|
\lesssim
C_\lambda
\left[
\sqrt{
\frac{
M\mathfrak Q_\lambda(u_{\mathcal F})
(\log n+s)
}{n}
}
+
\frac{M(\log n+s)}{n}
\right].
\]
The Lipschitz estimate of \autoref{lemma: appE residual two parameter lipschitz detailed} extends the same bound from the net to all \(\theta\in\Theta_\lambda^c\), with an additional \(C_{\lambda,M}/n\) term that is absorbed above.

Moreover,
\[
P_0b_{u_{\mathcal F}}^2
\le
2MP_0b_{u_{\mathcal F}}
\lesssim
M\mathfrak Q_\lambda(u_{\mathcal F}),
\qquad
0\le b_{u_{\mathcal F}}\le2M.
\]
Applying scalar Bernstein's inequality to this fixed function and taking a union bound completes the proof.
\end{proof}

\begin{proof}[Proof of~\autoref{theorem: oracle inequality}]
By~\autoref{proposition: basic oracle decomposition}, for every
\(v\in\mathcal F\),
\begin{align}
|\mathfrak L_\lambda^+(v)-\mathfrak L_{\lambda,n}^+(v)|
&\le
A_T
\sup_{\theta\in\Theta_\lambda^c}
|(P-P_n)r_{v,\theta}|
\notag\\
&\le
A_T
\sup_{\theta\in\Theta_\lambda^c}
|(P-P_n)h_{v,\theta}|
+
A_T\Delta_{*,n}.
\label{equation: appE positive loss decomposition detailed}
\end{align}

Since the outputs, test derivatives, and class-wise bias are bounded, there exists \(Q_{\max}>0\) such that $\mathfrak Q_\lambda(v)
\le
Q_{\max},
v\in\mathcal F.$
For \(n\ge3\), define
\[
L_n
:=
\max\left\{
1,
\left\lceil
\log_2\left(
\frac{Q_{\max}n}{\log n}
\right)
\right\rceil
\right\},
\]
and set
\[
t' := t+\log\bigl(16(L_n+1)\bigr).
\]

Choose
\begin{align}
\tau
:=
C\Bigg\{
&\mathfrak Q_\lambda(u_{\mathcal F})
+
C_\lambda^2
\frac{
\mathrm{VC}_{\mathcal F}\log n+t'
}{n}
\notag\\
&+
C_\lambda
\left[
\sqrt{
\frac{
\log n+t'
}{n}
}
+
\frac{
\log n+t'
}{n}
\right]
+
\tau_*(t')
+
\frac{\log n}{n}
\Bigg\},
\label{equation: appE base shell radius}
\end{align}
where \(C\) is sufficiently large.

If \(\tau\ge Q_{\max}\), then \(u_n\in\mathcal F_\tau\) trivially. Hence, suppose that \(\tau<Q_{\max}\), and let \(J\) be the smallest integer satisfying $2^J\tau \ge Q_{\max}.$ Since \(\tau\ge\log n/n\), we have \(J\le L_n\). Thus, $\mathcal F
= \mathcal F_\tau \cup \bigcup_{j=1}^J \left( \mathcal F_{2^j\tau} \setminus \mathcal F_{2^{j-1}\tau} \right).$ Write $ Y_j := 2^j\tau.$

Apply \autoref{lemma: appE empirical L2 localization} to every \(\mathcal F_{Y_j}\), \(0\le j\le J\), with confidence parameter
\(t'\). Apply \autoref{lemma: appE fixed comparator deviation} and \autoref{lemma: appE exact solution calibration detailed} once, also with parameter \(t'\). Since $\bigl(6(J+1)+5\bigr)e^{-t'} \le e^{-t},$ a union bound gives an event \(\Omega_t\) such that $\mathbb P(\Omega_t) \ge 1-e^{-t},$ on which all estimates hold simultaneously.

Suppose, toward a contradiction, that
\[
u_n
\in
\mathcal F_{Y_j}
\setminus
\mathcal F_{Y_j/2}
\]
for some \(j\ge1\). Then $u_n\in\mathcal F_{Y_j}, \mathfrak Q_\lambda(u_n) > \frac{Y_j}{2}.$ By the empirical minimality of \(u_n\),
\[
\begin{aligned}
\mathfrak Q_\lambda(u_n)
-
\mathfrak Q_\lambda(u_{\mathcal F})
\le{}&
\mathfrak Q_\lambda(u_n)
-
\mathfrak Q_{\lambda,n}(u_n)
+
\mathfrak Q_{\lambda,n}(u_{\mathcal F})
-
\mathfrak Q_\lambda(u_{\mathcal F}).
\end{aligned}
\]
Using \autoref{proposition: basic oracle decomposition}, \eqref{equation: appE local process}, \eqref{equation: appE local boundary process}, \autoref{lemma: appE fixed comparator deviation}, and \eqref{equation: appE positive loss decomposition detailed}, and absorbing the fixed factors \(A_0\) and \(A_T^2\) into the implicit constant, we obtain
\[
\begin{aligned}
\mathfrak Q_\lambda(u_n)
\lesssim{}&
\mathfrak Q_\lambda(u_{\mathcal F})
+
C_\lambda
\sqrt{
\frac{
\mathrm{VC}_{\mathcal F}Y_j\log n
}{n}
}
+
C_\lambda
\sqrt{
\frac{
Y_jt'
}{n}
}
\\
&+
C_\lambda\frac{Mt'}{n}
+
C_\lambda
\sqrt{
\frac{
M\mathfrak Q_\lambda(u_{\mathcal F})
(\log n+t')
}{n}
}
\\
&+
C_\lambda
\frac{
M(\log n+t')
}{n}
+
\Delta_{*,n}.
\end{aligned}
\]

Young's inequality gives
\[
\begin{aligned}
C_\lambda
\sqrt{
\frac{
\mathrm{VC}_{\mathcal F}Y_j\log n
}{n}
}
&\le
\frac{Y_j}{16}
+
CC_\lambda^2
\frac{
\mathrm{VC}_{\mathcal F}\log n
}{n},
\\
C_\lambda
\sqrt{
\frac{
Y_jt'
}{n}
}
&\le
\frac{Y_j}{16}
+
CC_\lambda^2\frac{t'}{n},
\\
C_\lambda
\sqrt{
\frac{
M\mathfrak Q_\lambda(u_{\mathcal F})
(\log n+t')
}{n}
}
&\le
\frac14
\mathfrak Q_\lambda(u_{\mathcal F})
+
CC_\lambda^2
\frac{
\log n+t'
}{n}.
\end{aligned}
\]
Together with \eqref{equation: appE calibration detailed}, this gives
\[
\begin{aligned}
\mathfrak Q_\lambda(u_n)
\le{}&
C\mathfrak Q_\lambda(u_{\mathcal F})
+
\frac{Y_j}{8}
+
CC_\lambda^2
\frac{
\mathrm{VC}_{\mathcal F}\log n+t'
}{n}
\\
&+
CC_\lambda
\left[
\sqrt{
\frac{
\log n+t'
}{n}
}
+
\frac{
\log n+t'
}{n}
\right].
\end{aligned}
\]

By \eqref{equation: appE base shell radius}, all terms in the preceding bound other than \(Y_j/8\) are bounded by \(\tau/4\), after increasing \(C\). Since \(Y_j\ge2\tau\), $\mathfrak Q_\lambda(u_n) \le \frac{Y_j}{8} + \frac{\tau}{4} \le \frac{Y_j}{4} < \frac{Y_j}{2},$ contradicting $\mathfrak Q_\lambda(u_n) > \frac{Y_j}{2}.$ Therefore, $u_n \in \mathcal F_\tau$ on \(\Omega_t\).

It follows from \eqref{equation: appE base shell radius} and \autoref{lemma: appE critical radius bound} that
\begin{align}
\mathfrak Q_\lambda(u_n)
\lesssim{}&
\mathfrak Q_\lambda(u_{\mathcal F})
+
C_\lambda^2
\frac{
\mathrm{VC}_{\mathcal F}\log n+t
}{n}+
C_\lambda
\left[
\sqrt{
\frac{
\log n+t
}{n}
}
+
\frac{
\log n+t
}{n}
\right]
+
\frac{
CM^2e^{-t}
}{
L_n+1
}.
\label{equation: appE general oracle detailed}
\end{align}
This estimate holds with probability at least \(1-e^{-t}\). Moreover, $L_n+1 \asymp \log\left( \frac{n}{\log n} \right)$ for sufficiently large \(n\).

Set $\Gamma_n := \mathrm{VC}_{\mathcal F}\log n.$ Assume $c_0\sqrt{n\log n} \le \Gamma_n \le c_1n.$ Choose $t_n := c_2 \frac{ \Gamma_n^2 }{n},$ where \(c_2>0\) is sufficiently large. Then
\[
\frac{t_n}{n}
\lesssim
\frac{\Gamma_n}{n},
\qquad
\sqrt{
\frac{
\log n+t_n
}{n}
}
\lesssim
\frac{\Gamma_n}{n},
\qquad
\frac{
\log n+t_n
}{n}
\lesssim
\frac{\Gamma_n}{n}.
\]
Moreover, after choosing \(c_2\) such that \(t_n\ge\log n\), $\frac{e^{-t_n}}{L_n+1}
\le
\frac1n
\lesssim
\frac{\Gamma_n}{n}.$
Substituting \(t=t_n\) into \eqref{equation: appE general oracle detailed}, and using \(C_\lambda\ge1\), gives
\[
\mathfrak Q_\lambda(u_n)
\lesssim
\mathfrak Q_\lambda(u_{\mathcal F})
+
C_\lambda^2
\frac{
\mathrm{VC}_{\mathcal F}\log n
}{n}.
\]
The probability is at least $1- \exp\left( -c_2 \frac{ \bigl( \mathrm{VC}_{\mathcal F}\log n \bigr)^2 }{n} \right).$ This proves \eqref{equation: quadrature fast probability}-- \eqref{equation: theorem6 VCF oracle fast}.
\end{proof}

\bibliographystyle{plain}
\bibliography{main.bib}
\end{document}